\documentclass{article}

\usepackage{microtype}
\usepackage[nottoc]{tocbibind}
\usepackage[utf8]{inputenc}
\usepackage{graphicx}
\usepackage{float}
\usepackage{scalerel}
\usepackage{bm}
\usepackage{amsthm}
\usepackage{amsmath}
\usepackage{physics}
\usepackage{amssymb}
\usepackage{mathrsfs}
\usepackage{dsfont}
\usepackage{wasysym}
\usepackage[all]{xy}
\usepackage{tikz-cd}
\usepackage{stmaryrd}
\usepackage{enumitem}
\usepackage{tabularx}
\usepackage{placeins}
\usepackage{hyperref}
\usepackage{tikz}
\usetikzlibrary{shapes.geometric,calc}
\usetikzlibrary{shapes,arrows,chains}
\usetikzlibrary{decorations.markings}
\usetikzlibrary{decorations.pathmorphing}
\usetikzlibrary{shapes.multipart}
\tikzset{snake it/.style={decorate, decoration=snake}}
\tikzstyle{GraphNode}=[circle, draw=black, fill=black, inner sep=2pt, minimum size=5pt]
\tikzstyle{GraphEdge}=[black]

\theoremstyle{definition}
\newtheorem{Definition}{Definition}[subsection]

\theoremstyle{plain}
\newtheorem{Theorem}[Definition]{Theorem}

\theoremstyle{plain}

\theoremstyle{plain}
\newtheorem{Proposition}[Definition]{Proposition}

\theoremstyle{plain}
\newtheorem{Lemma}[Definition]{Lemma}

\theoremstyle{plain}
\newtheorem{Corollary}[Definition]{Corollary}

\theoremstyle{plain}

\newcounter{mainthm}

\newtheorem{maintheorem}[mainthm]{Theorem}

\theoremstyle{plain}

\theoremstyle{plain}

\theoremstyle{definition}

\theoremstyle{definition}
\newtheorem{Example}[Definition]{Example}

\theoremstyle{definition}

\theoremstyle{remark}
\newtheorem{Remark}[Definition]{Remark}

\theoremstyle{plain}
\newcommand{\thistheoremname}{}
\newtheorem*{genericthm*}{\thistheoremname}
\newenvironment{namedthm*}[1]
  {\renewcommand{\thistheoremname}{#1}%
   \begin{genericthm*}}
  {\end{genericthm*}}
  
\newcommand\cat[1]{\mathbf{#1}}

\newcommand\Mod{\cat{Mod}}
\newcommand\rRep{\cat{Rep}}

\newcommand\Bimod{\cat{Bimod}}

\newcommand\Vect{\cat{Vect}}
\newcommand\supp{\mathrm{supp}}
\newcommand\Z{\mathbb{Z}}
\newcommand\cG{\mathcal{G}}
\newcommand\fC{\mathfrak{C}}

\title{Fiber 2-Functors and Tambara-Yamagami Fusion 2-Categories}
\author{Thibault D. Décoppet and Matthew Yu}
\date{December 2024}

\begin{document}

\bibliographystyle{alpha}

\maketitle
    \hspace{1cm}
    \begin{abstract}
        We introduce group-theoretical fusion 2-categories, a categorification of the notion of a group-theoretical fusion 1-category. Physically speaking, such fusion 2-categories arise by gauging subgroups of a global symmetry. We show that group-theoretical fusion 2-categories are completely characterized by the property that the braided fusion 1-category of endomorphisms of the monoidal unit is Tannakian. Then, we describe the underlying finite semisimple 2-category of group-theoretical fusion 2-categories, and, more generally, of certain 2-categories of bimodules. We also partially describe the fusion rules of group-theoretical fusion 2-categories. Using our previous results, we classify fusion 2-categories admitting a fiber 2-functor. Next, we study fusion 2-categories with a Tambara-Yamagami defect, that is $\mathbb{Z}/2$-graded fusion 2-categories whose non-trivially graded factor is $\mathbf{2Vect}$. We classify these fusion 2-categories, and examine more closely the more restrictive notion of Tambara-Yamagami fusion 2-categories. Throughout, we give many examples to illustrate our various results.
    \end{abstract}

\tableofcontents

\section*{Introduction}

Perhaps the most well-studied class of fusion 1-categories consist of those fusion 1-categories that admit a fiber functor. Namely, it follows from Tannaka-Krein duality \cite{EGNO} that such fusion 1-categories are equivalent to the 1-categories of representation of a finite semisimple Hopf algebra (see also \cite{JS2} for a historical reference). Further, finite dimensional Hopf algebras, and especially semisimple ones, have been extensively studied (see, for instance, \cite{LR, EG}). Categorifying the notion of a Hopf algebra, Hopf 1-categories were originally introduced in \cite{CF}. However, only few examples of Hopf 1-categories have been constructed thus far \cite{Pfe}. More recently, fusion 2-categories were introduced in \cite{DR}. The corresponding version of Tannaka-Krein duality was established in \cite{Green}, i.e.\ it is shown that finite semisimple Hopf 1-categories correspond to fusion 2-categories equipped with a fiber 2-functor. This provides a strong mathematical motivation for the study of such fusion 2-categories. From a physical perspective, fiber functors for fusion 1-categories are useful for identifying in which fusion 1-categories it is possible to condense out all the objects. In particular, fiber functors were used to discuss anomalies for categorical symmetries for theories in $(1+1)$-dimensions in \cite{Thorngren:2019iar,KNZZ23}. We expect that that fiber 2-functors for fusion 2-categories can be used to carry out a similar analysis for theories in $(2+1)$-dimensions.

Classifying finite semisimple Hopf algebras remains an open problem. One of the main difficulty is the existence of non group-theoretical finite semisimple Hopf algebras. More precisely, group-theoretical fusion 1-categories as in \cite{ENO1} are constructed using exclusively finite groups and cocycles. Further, fiber functors on such fusion 1-categories are completely understood \cite{O2}, and it is natural to ask whether every fusion 1-category admitting a fiber functor is group-theoretical \cite{ENO1}. This question was answered in the negative in \cite{N2}, and thereby provides examples of non group-theoretical semisimple Hopf algebras. We will see that the categorified version of this question has a very different answer.

One context in which fiber functors for fusion 1-categories naturally arise is the study of Tambara-Yamagami 1-categories. Namely, recall that Tambara-Yamagami 1-categories are $\mathbb{Z}/2$-graded fusion 1-categories $\mathcal{C} = \mathcal{C}_+\oplus \mathcal{C}_-$ such that $\mathcal{C}_+$ is group-like, also called pointed, and $\mathcal{C}_-=\mathbf{Vect}$. Then, $\mathcal{C}_-$ provides a fiber functor for $\mathcal{C}_+$. Moreover, $\mathcal{C}_-$ implements a Morita autoequivalence, i.e.\ a self-duality, of $\mathcal{C}_+$. This observation was used in \cite{ENO2} to recover the classification of Tambara-Yamagami 1-categories obtained in \cite{TY} using homotopy-theoretic methods. On the other hand, the original approach in \cite{TY} proceeds instead by direct computations. We also wish to point out that Tambara-Yamagami 1-categories are not always group-theoretical, and that the question of the existence of fiber functors for Tambara-Yamagami 1-categories was considered in \cite{Tam}.

Tambara-Yamagami 1-categories have also appeared in Physics as they provide the mathematical framework to describe the algebra of extended operators of $(1+1)$-dimensional theories that are self-dual and admit a \textit{duality defect} \cite{JSW19,LS}. The most famous example of such a theory is the Ising model at critial temperature  \cite{FFRS,FFRS09}. The switching of correlators of spins with the correlators of disorder, or twist, operators at the critical temperature yields a Morita autoequivalence. More precisely, gauging the 0-form $\Z/2$ symmetry on half of a $(1+1)$-dimensional spacetime leads to the Kramers-Wanier duality defect \cite{CCHLS}. This procedure provides a boundary for the theory and imposing Dirichlet conditions makes the defect topological, which therefore gives the desired self-duality. Many other applications of duality defects have been studied in Physics such as in \cite{Thorngren:2019iar,Burbano:2021loy,CLSWY:2018iay,LS,Choi,KOZ21,KOZ22}.

In this work, we are interested in Tambara-Yamagami fusion 2-categories, a categorification of the notion of a Tambara-Yamagami fusion 1-category. Slightly more generally, a fusion 2-category with a Tambara-Yamagami defect is a $\Z/2$-graded fusion 2-category $\mathfrak{C}=\mathfrak{C}_+\boxplus\mathfrak{C}_-$ such that $\mathfrak{C}_-=\mathbf{2Vect}$. In particular, $\mathfrak{C}_-$ provides a fiber 2-functor for $\mathfrak{C}_+$. But, unlike in the 1-categorical case, we will see that every fusion 2-category that admits a fiber 2-functor is group-theoretical. This implies that every fusion 2-category with a Tambara-Yamagami defect is automatically group-theoretical. As a consequence, we obtain a classification of Tambara-Yamagami fusion 2-categories from a general classification result for group-theoretical fusion 2-categories. In addition, our analysis of group-theoretical fusion 2-categories allows us to completely describe fusion 2-categories admitting a fiber 2-functor.

\subsection*{Results}

Throughout, we will work over an algebraically closed field $\mathds{k}$ of characteristic zero. We begin by introducing a categorification of the notion of a group-theoretical fusion 1-category using work of the first author on fusion 2-categories \cite{D7, D8}. 
Fusion 2-categories can be used to describe the interactions of topological operators (surfaces and lines) for $(2+1)$-dimensional theories and we note that group-theoretical fusion 2-categories have already appeared in the Physics literature \cite{BBSNT1} (see also \cite{BBFP, BBSNT}). Such fusion 2-categories are exactly the ones that arise by considering a finite group $G$, potentially equipped with an anomaly $\pi$, and gauging a subgroup $H$, equipped with a trivialization $\psi$ of the anomaly. 
We will use the notation $\mathfrak{C}(G,H,\pi,\psi)$ to denote the corresponding fusion 2-category;\ the precise definition being given at the beginning of \S\ref{section:groupF2C}. Examples include the fusion 2-categories $\mathbf{2Vect}_{\cG}$ of 2-vector spaces graded by $\cG$, and $\mathbf{2Rep}(\cG)$ of 2-representations of $\cG$, where $\cG$ is a finite 2-group. 
We begin by establishing the following detection criterion for group-theoretical fusion 2-categories.

\begin{maintheorem}[Theorem \ref{thm:grouptheoreticalrecognition}]
A fusion 2-category $\mathfrak{C}$ is group-theoretical if and only if $\Omega\mathfrak{C}$ is a Tannakian fusion 1-category, i.e.\ is equivalent as a braided fusion 1-category to the 1-category of representations of a finite group.
\end{maintheorem}

\noindent In particular, this implies that every fusion 2-category admitting a fiber 2-functor is group-theoretical, so that every finite semisimple Hopf 1-category is group-theoretical. This is in sharp contrast with the decategorified setting recalled above. Then, we obtain a description of the underlying finite semisimple 2-category of any group theoretical fusion 2-category, and, crucially, of certain finite semisimple bimodule 2-categories between them. This generalizes results of \cite{El} and \cite{KTZ}. We also give a partial description of the fusion rules of group-theoretical fusion 2-categories, and investigate when two group-theoretical fusion 2-categories are monoidally equivalent.

Next, we consider fusion 2-categories admitting a fiber 2-functor, that is fusion 2-categories $\mathfrak{C}$ for which there exists a monoidal 2-functor $\mathfrak{C}\rightarrow \mathbf{2Vect}$. We obtain the following categorification of corollary 3.1 of \cite{O2}. To state our result, recall that an exact factorization of a finite group $G$ is the data of two subgroups $H$ and $K$ of $G$ such that $H\cap K=\{e\}$ and $HK=G$.

\begin{maintheorem}[Theorem \ref{thm:fiber2functor}]
A fusion 2-category admits a fiber 2-functor if and only if it is equivalent to a fusion 2-category of the form
$\mathfrak{C}(G,H,\pi,\psi)$, with $G$ a finite group admitting an exact factorization by two subgroups $H$ and $K$, $\pi\in Z^4_{gp}(G; \mathds{k}^\times)$, $\psi\in C^3_{gp}(H;\mathds{k}^\times)$ is such that $\pi|_{H} = d\psi$, and $\pi|_K$ is a coboundary.
\end{maintheorem}

\noindent Furthermore, we describe the set of equivalence classes of fiber 2-functors on any such fusion 2-category. Moreover, we exhibit the dual fusion 2-category in the sense of \cite{D8} to a fusion 2-category with respect to the module 2-category associated to a fiber 2-functor. We also use the above theorem to argue that some of the examples of group-theoretical fusion 2-categories that we have previously introduced are not equivalent to the fusion 2-category of 2-representations of a finite 2-group.

The next part of this work examines duality surface defects. In the Physics literature, duality operators are only discussed in even dimensions for the reason that it is the only case in which gauging a $p$-form symmetry may result in a dual $p$-form symmetry. Given our particular focus on $(2+1)$-dimensions, gauging a $0$- or $1$-form symmetry with the structure of a finite group results in a dual $1$- or $0$-form symmetry respectively. In this sense, one might not have expected interesting self-dualities to exist as in the cases of $(1+1)$- or $(3+1)$-dimensions. In spite of this, for a specific (non-split) finite 2-group $\cG$, we do observe the existence of an exotic fiber functor for $\mathbf{2Rep}(\cG)$, for which the corresponding dual fusion 2-category is again $\mathbf{2Rep}(\cG)$ (and not $\mathbf{2Vect}_{\cG}$). This predicts an interesting self-duality in $(2+1)$-dimensions that is not visible via the usual framework of gauging a subgroup \`{a} la \cite{Tachikawa17}, and requires a more subtle analysis of the category.

As a first step towards our mathematical analysis of self-duality in $(2+1)$-dimensions, we study the general properties of group-graded fusion 2-categories. We also give many interesting examples of group-theoretical fusion 2-categories, for which we completely describe both the underlying finite semisimple 2-category and the fusion rules. Then, in order to specifically capture the concept of self-duality, we introduce the notion of a fusion 2-category with a Tambara-Yama\-gami defect, that is a $\Z/2$-graded fusion 2-category $\mathfrak{C}=\mathfrak{C}_+\boxplus\mathfrak{C}_-$ such that $\mathfrak{C}_-=\mathbf{2Vect}$. As is familiar from the theory of Tambara-Yamagami 1-categories recalled above, the factor $\mathfrak{C}_-$ implements a self-dual duality, i.e.\ a Morita autoequivalence, of the fusion 2-category $\mathfrak{C}_+$. In particular, $\mathfrak{C}_+$ admits a fiber 2-functor, so that such a fusion 2-category $\mathfrak{C}$ is always group-theoretical. We can therefore use our previous results to classify fusion 2-categories with a Tamabara-Yamagami defect using finite groups and cocycles. We also show that the fusion rules of the Tamabara-Yamagami defect resembles those of the duality line in Tambara-Yamagami 1-categories. 

We move on to study further the categorified version of the Tamabara-Yamagami 1-categories introduced in \cite{TY}. More precisely, we define a Tambara-Yamagami 2-category as a fusion 2-category with a Tamabara-Yamagami defect $\mathfrak{C}=\mathfrak{C}_+\boxplus\mathfrak{C}_-$ such that $\mathfrak{C}_+$ is $\textbf{2}\Vect_{\cG}$ for some finite 2-group $\mathcal{G}$. In particular, the relations between the concepts that we have introduced are summarized by the following inclusions:
\[
    \text{TY F2Cs} \subset \text{F2Cs with a TY defect} \subset \text{group theoretical F2Cs}\,.
\]
\noindent Then, we obtain the following 2-categorical analogue of the classification results of \cite{TY}.

\begin{maintheorem}[Proposition \ref{prop:2TYclassification}, Proposition \ref{prop:fusionD}]
For any finite abelian group $A$, let us write $A\wr \Z/2 = (A \oplus A) \rtimes \Z/2$ for the wreath product.
\begin{enumerate}
    \item Every Tambara-Yamagami 2-category can be constructed, up to monoidal equivalence, from the data of a finite abelian group $A$ together with a class $\pi$ in $H^4_{gp}(A\wr\mathbb{Z}/2;\mathds{k}^\times)$ whose restriction to $A\oplus A$ is trivial. We write $\mathbf{2TY}(A,\pi)$ for the corresponding (equivalence class of) fusion 2-category.
    \item Two Tambara-Yamagami 2-categories $\mathbf{2TY}(A_1,\pi_1)$ and $\mathbf{2TY}(A_2,\pi_2)$ are equivalent if and only if there exists a group isomorphism $f:A_1\wr \mathbb{Z}/2\cong A_2\wr \mathbb{Z}/2$ such that $f(A_1\oplus 0) = A_2\oplus 0$, and $f^*\pi_2/\pi_1$ is the trivial class.
    \item The trivially graded factor of $\mathbf{2TY}(A,\pi)$ is $$\mathbf{2TY}(A,\pi)_+ = \mathbf{2Vect}_{A\lbrack 1\rbrack\times A\lbrack 0\rbrack}= \underbrace{\mathbf{Mod}(\mathbf{Vect}_A)\boxplus ... \boxplus \mathbf{Mod}(\mathbf{Vect}_A)}_{|A|\ \mathrm{times}},$$ and, under this decomposition, the fusion of the TY defect $D$ is given by $$D\Box D = \boxplus_{a\in A} \mathbf{Vect}_{a},$$ where $\mathbf{Vect}_{a}$ is a condensation surface defect in the $a$-th copy of $\mathbf{Mod}(\mathbf{Vect}_A)$.
\end{enumerate}
 
\end{maintheorem}

\noindent We note that the Tambara-Yamagami 2-category $\mathbf{2TY}(\mathbb{Z}/2,triv)$ has first appeared in \cite{BBSNT1}, and then in \cite{BBFP:I,BBFP,BBSNT}, which provided the original motivation for this work. Following \cite{Tam}, it is interesting to try to determine when a Tambara-Yamagami 2-category admits a fiber 2-functor. If $A$ has odd order, then there is a unique Tambara-Yamagami 2-category $\mathbf{2TY}(A,triv)$, and we show that this fusion 2-category admits a fiber 2-functor. However, we prove that $\mathbf{2TY}(A,triv)$ is not equivalent to the 2-category of 2-representations of any finite 2-group so long as $A$ is non-trivial. With $A = \mathbb{Z}/2$, we find that all the associated Tambara-Yamagami 2-categories admit a fiber 2-functor, and we identify two of them with the 2-categories of 2-representations of a finite 2-group.

\subsection*{Acknowledgments}

We would like to thank Lakshya Bhardwaj, Clement Delcamp, Andrea Ferrari, Theo Johnson-Freyd, Adrià Mar\'in Salvador, Sean Sanford, Sakura Sch\"afer-Nameki, and Ryan Thorngren for conversations. We would also like to thank the referees for their numerous thorough and constructive comments. In particular, we would like to express our most profound gratitude towards one of the referees for spotting a subtle inaccuracies in the previous versions of the statement of proposition \ref{prop:classificationgrouptheoretical} and corollary \ref{cor:fiber2functor}.

\section{Physical Intuition}\label{sec:physicalintuition}

We presently discuss the place that the mathematical objects that we will study inhabit in Physics. The most relevant notion is that of a categorical symmetry, that is a physical theory on which a category of symmetries acts. Manipulations that change the category result in different symmetries, and we say that the theory is in a different ``phase" when the category of symmetries that acts has changed.

This work relies on the construction that associates a 2-category of bimodules $\Bimod_{\fC}(\mathcal{A})$ to an algebra object $\mathcal{A}$ in a fusion 2-category $\fC$. Physically, we adopt the same perspective as \cite{BT17} where the category
$\Bimod_{\fC}(\mathcal{A})$
describes a particular phase that is 
the result of \textit{condensing} $\mathcal{A}$ in the original phase.
This is a type of topological manipulation that is done to our theory which changes its symmetries in a controlled manner, thereby preserving some structure. A concrete realization of such a manipulation is when one can tune a parameter adiabatically in time to change a given system in such a way that the original system can be recovered by slowly tuning the parameter in the opposite direction.
This adiabatic procedure creates a gapped interface in the theory and no information about the dynamics has been lost, only the symmetries have been modified.
For more physical examples of condensations in the case of 1-categories, we refer the reader to \cite{Burnell,Yu:2021zmu,Hung:2015hfa}.
Condensing is a type of topological manipulation and it is possible to recover the original category by reversing the manipulation, in the manner previously discussed. 
Furthermore, we will assume that the condensed phase only has a single ground state which puts additional constraints on the algebra.

Fusion 2-categories act as symmetries of (2+1)d theories, which means that that we can surface and line operators that implement the symmetry. As a specific example, associated to the finite group $G$, we can think of the objects of $\mathbf{2Vect}_G$ as surface operators that implement a $0$-form $G$ symmetry. For a given finite group $G=A[0]\rtimes H[0]$, we can consider the two rigid algebras $\Vect_{H}$ and $\Vect_{A}$. Taking bimodules with respect to these algebras then amounts to gauging a quotient or a subgroup.\footnote{There is a slight subtlety that deserves to be mentioned: Gauging $A$ results in an emergent dual symmetry that is 1-form, but we get condensation defects of the topological Wilson lines upon choosing to gauge $H$.} 
The resulting fusion 2-categories are
\begin{align}\label{eq:bim1}
    \Bimod_{\textbf{2}\Vect_{G}}(\Vect_{H}) & \simeq \mathbf{2Rep}(\cG)\,,  \\\label{eq:bim2}
    \Bimod_{\textbf{2}\Vect_{G}}(\Vect_{A}) & \simeq \textbf{2}\Vect_{\cG}\,, 
\end{align}
where $\cG$ is the finite 2-group $\widehat{A}[1]\rtimes H[0]$. 
Above and throughout, we use $\widehat{A}$ to denote the Pontryagin dual group of $A$. 
In the example above, we see that there is an emergent dual 1-form group. The above equivalences have already appeared in \cite{BBFP,BBSNT,DelT}, and we offer a rigorous proof in \S\ref{section:Fib2Functor}. While the two categories $\mathbf{2Rep}(\cG)$ and $\textbf{2}\Vect_{\cG}$ are in general not monoidally equivalent, they are Morita equivalent. In particular, they have the same centers \cite{D9}. More precisely, there is a sequence of gapped procedures that can be used to go between these two 2-categories. In fact, this can be done directly by gauging or ungauging the 2-group $\cG$.

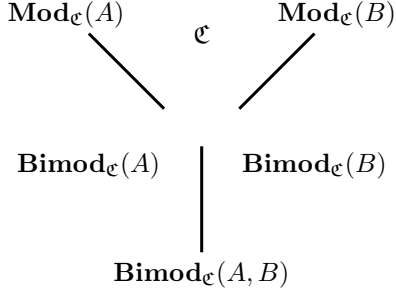
\begin{figure}
    \centering
   \begin{tikzpicture}
    \draw[line width = .4 mm] (0,0)--(1,1);
      \draw[line width = .4 mm] (-1,0)--(-2,1);
      \draw[line width = .4 mm] (-.5,-.5)--(-.5,-.5-1.41);
      \draw (1,-.75) node {$\Bimod_{\fC}{(\mathcal{B})}$};
       \draw (-2,-.75) node {$\Bimod_{\fC}{(\mathcal{A})}$};
       \draw (-2.3,1.25) node {$\Mod_{\fC}(\mathcal{A})$};
       \draw (1.3,1.25) node {$\Mod_{\fC}(\mathcal{B})$};
       \draw (-.5,-.5-1.7) node {$\Bimod_{\fC}(\mathcal{A},\mathcal{B})$};
       \draw (-.5,1) node {$\fC$};
   \end{tikzpicture}
    \caption{Starting with some fusion 2-category $\fC$, taking bimodules with respect to some algebras $\mathcal{A}$ and $\mathcal{B}$ moves between different bulks, namely $\Bimod_{\fC}(\mathcal{A})$ and $\Bimod_{\fC}(\mathcal{B})$. The solid lines denote boundaries. In particular, $\Bimod_{\fC}(\mathcal{A},\mathcal{B})$ is a boundary between $\Bimod_{\fC}(\mathcal{A})$ and $\Bimod_{\fC}(\mathcal{B})$.}
    \label{fig:bimodAB}
\end{figure}

When studying condensations between phases described by fusion 2-categories, it is natural to ask whether one can condense to the vacuum, in this case $\textbf{2}\Vect$. Fiber 2-functors are maps from a fusion 2-category $\fC$ to $\textbf{2}\Vect$, or, equivalently, $\fC$-module structures on $\textbf{2}\Vect$.
The existence of a fiber functor therefore corresponds to the existence of a \textit{gapped trivial boundary}.
These gapped trivial boundaries are important in condensed matter settings and are related to Lagrangian algebras \cite{Lan,Kong}.
For our purposes, starting with a fusion 2-category $\fC$, fiber 2-functors for the fusion 2-category $\Bimod_{\fC}(\mathcal{A})$ will arise through an equivalence $\Bimod_{\fC}(\mathcal{A},\mathcal{B})\simeq \textbf{2}\Vect$ for some algebras $\mathcal{A}$ and $\mathcal{B}$. The dual fusion 2-category to $\Bimod_{\fC}(\mathcal{A})$ with respect to this fiber 2-functor is then $\Bimod_{\fC}(\mathcal{B})$. Physically, the process is described in figure \ref{fig:bimodAB} where one maps from a theory with symmetry given by $\fC$ to a condensed phase with symmetry given by either $\Bimod_{\fC}(\mathcal{A})$ or $\Bimod_{\fC}(\mathcal{B})$.

Any given fusion 2-category $\fC$ need not admit a single fiber 2-functor, but can in fact admit multiple. For fusion 1-categories, this behaviour is well-known. For instance, $\rRep(D_8)$ admits three distinct fiber functors \cite{Tam}. In section 3.2.5 of \cite{Thorngren:2019iar}, the authors show that the boundary between them must support gapless edge modes. These edge modes give some amount of noninvertibility to the boundary-changing operator between two fiber functors. Without the modes, two fiber functors would have equivalent associated module categories, contradicting the fact that the fiber functors are inequivalent. We therefore expect that there exists gapless edge modes protected by fusion 2-category symmetry in $(3+1)$-dimensions associated to the multiple fiber 2-functors that will be discussed in \S\ref{section:Fib2Functor}.

We now give another physical  perspective on fiber 2-functors. Namely, the fusion 2-category $\fC$ admits a fiber 2-functor precisely if there exists an algebra such that the objects of $\mathfrak{C}$ are either condensates of the algebra, or are screened by the algebra. This viewpoint will be particularly useful for the intuition behind the classification of fusion 2-categories that admit a fiber 2-functor in \S\ref{section:Fib2Functor}.

\section{Mathematical Preliminaries}\label{section:prelim}
\subsection{Fusion 2-Categories associated to Groups}\label{sub:F2Cgroup}

Let us fix $\mathds{k}$ an algebraically closed field of characteristic zero. When working with monoidal 2-categories, we will use $\Box$ to denote the monoidal product, $I$ to denote the monoidal unit, and otherwise follow the notations of \cite{GPS} (see also section 2.3 of \cite{SP}). We will focus our attention on fusion 2-categories as introduced in \cite{DR}, which are finite semisimple monoidal 2-categories with duals and simple monoidal unit.

For our purposes, the most important class of fusion 2-categories will be the fusion 2-categories $\mathbf{2Vect}_G$ of $G$-graded 2-vector spaces for some finite group $G$, with monoidal structure given by the grading preserving Deligne tensor product $\boxtimes$. More generally, given $\pi\in Z^4_{gp}(G, \mathds{k}^\times)$, a 4-cocycle for the group $G$ with coefficients in $\mathds{k}^{\times}$ (see e.g.\ section 1.7 of \cite{EGNO} for the precise definition),\footnote{We will always assume that the group cochains that we consider are normalized. Thanks to \cite{EML}, this is not a loss of generality.} we may twist the pentagonator of $\textbf{2}\Vect_G$. More precisely, given for $i=1,2,3,4$ a 2-vector space $V_i$ with grading $g_i\in G$, then the invertible modification $\pi$ witnessing the coherence $$\begin{tikzcd}[sep=small]
 {(V_1\boxtimes(V_2\boxtimes V_3))\boxtimes V_4} \arrow[d, "\mathfrak{a}"'] &                       &     ((V_1\boxtimes V_2)\boxtimes V_3)\boxtimes V_4 \arrow[ll, "\mathfrak{a} 1"'] \arrow[d, "\mathfrak{a}"]\\
V_1\boxtimes((V_2\boxtimes V_3)\boxtimes V_4) \arrow[rd, "1\mathfrak{a}"']{}\arrow[rr, "\pi", Rightarrow, shorten >=10ex, shorten <=10ex ]                                                                               &       & (V_1\boxtimes V_2)\boxtimes (V_3\boxtimes V_4) \arrow[ld, " \mathfrak{a}"] \\
 & V_1\boxtimes(V_2\boxtimes (V_3\boxtimes V_4))  &  
\end{tikzcd}$$ is given by multiplication by $\pi(g_1,g_2,g_3,g_4)$. This produces the fusion 2-category $\textbf{2}\Vect^\pi_G$. Even though this construction needs a 4-cocycle as input, the resulting fusion 2-category depends up to monoidal equivalence only on the class of $\pi$ in the cohomology group $H^4_{gp}(G, \mathds{k}^\times)$. As was shown in \cite{JFY}, such fusion 2-categories can be completely characterized via an elementary property:\ A bosonic strongly fusion 2-category is a fusion 2-category $\mathfrak{C}$ such that $\Omega\mathfrak{C}$, the braided fusion 1-category of endomorphism of the monoidal unit, is $\mathbf{Vect}$. Every bosonic strongly fusion 2-category is of the form $\textbf{2}\Vect^\pi_G$ for some finite group $G$ and 4-cocycle $\pi$.

More generally, we will also be interested in the fusion 2-categories associated to a finite 2-group $\mathcal{G}$. The data of a finite 2-group consists of a finite group $H[0]$ and a finite abelian group $A[1]$, where the brackets denote the degree in which each group lives and gives the ``form number" of the symmetry that the group describes.
Said differently, for any non-negative integer $n$ and group $B$, abelian if $n\geq 1$, we use $B[n]$ to denote the higher group whose classifying space is the Eilenberg-MacLane space $K(B,n+1)$.
As part of the data of the finite 2-group $\mathcal{G}$, there is also an action of the degree-zero group on the degree-one group by a map $\rho: H \to \mathrm{Aut}(A)$. Moreover, there is a $k$-invariant given by a Postnikov class $\beta\in H^3_{gp}(H;A_{\rho})$. We denote the finite 2-group corresponding to this data by $A[1]\boldsymbol{\cdot}H[0]$. When the Postnikov class $\beta$ is trivial this gives a split 2-group, which we write as $A[1]\rtimes H[0]$. In this case, it is possible to specify local projective actions of the symmetry $H$ on each topological operator that is labeled by $A$. In particular, there could be a natural transformation from acting with $\rho_g \circ \rho_h$ and acting by $\rho_{g \circ h}$, where $\rho: H \rightarrow Aut(A)$. The equivalence classes of choices are classified by a class in $H^2_{gp}(H;A_{\rho})$. This class is referred to in the Physics literature as the \textit{symmetry fractionalization} \cite{BBCW,DGHK,HS19,Yu20}. If, in addition, the action by $H$ on $A$ is trivial, then the 2-group is a direct product $A[1]\times H[0]$. For applications of 2-groups to field theory and string theory beyond our categorical setup, see \cite{CDI,Baez05,Baez2:05,Benini:2018reh,Debray:2023rlx,Sharpe15}.

Given a finite 2-group $\mathcal{G} = A[1]\boldsymbol{\cdot}H[0]$, we can consider the fusion 2-category $\mathbf{2Vect}_{\cG}$ of $\mathcal{G}$-graded 2-vector spaces. The underlying finite semisimple 2-category has $|H|$ connected components, all of which are equivalent to $\mathbf{Mod}(\mathbf{Vect}_A)$, the 2-category of finite semisimple $\mathbf{Vect}_A$-module 1-categories. Further, the fusion rules of $\textbf{2}\Vect_{\cG}$ are completely understood. As was explained in construction 2.1.16 of \cite{DR}, the structure constraints of the fusion 2-category $\textbf{2}\Vect_{\cG}$ can be twisted using a 4-cocycle $\pi \in Z^4(\mathrm{B}\cG, \mathds{k}^\times)$, the group of 4-cocycles for the space $\mathrm{B}\cG$ with coefficients in $\mathds{k}^\times$. Broadly speaking, the fusion 2-category $\mathbf{2Vect}^\pi_{\cG}$ is defined as the Cauchy completion of the linearization of $\cG$ twisted by $\pi$. We note that if $\cG = G[0]$ is a finite group, then $H^4(\mathrm{B}G; \mathds{k}^\times)\cong H^4_{gp}(G; \mathds{k}^\times)$, and we recover the twisting procedure that we have previously reviewed.

Given a finite 2-group $\cG$, there is also a fusion 2-category $\mathbf{2Rep}(\cG)$ of 2-representations of $\cG$, i.e.\ 2-functors from $\mathrm{B}\mathcal{G}$ to $\textbf{2}\Vect$. The underlying 2-category was extensively studied in \cite{El}. On the other hand, there is yet no general description of its fusion rules. If $\cG = G[0]$ is a finite group, then it was shown in \cite{D2} that there is an equivalence of symmetric fusion 2-categories between $\mathbf{2Rep}(G)$ and $\mathbf{Mod}(\mathbf{Rep}(G))$, the fusion 2-category of finite semisimple $\mathbf{Rep}(G)$-module 1-categories.

\subsection{Morita Equivalences between Fusion 2-Categories}\label{sub:MoritaEquivalence}

We now give a general mathematical account of the constructions considered above in section \ref{sec:physicalintuition}. In the next sections, we will specialize these constructions to the bosonic strongly fusion 2-categories introduced at the beginning of \S\ref{sub:F2Cgroup}. We presently consider the general case of a fixed arbitrary fusion 2-category $\mathfrak{C}$. A \textit{rigid algebra} $\mathcal{A}$ in $\mathfrak{C}$ is an algebra whose multiplication 1-morphism has a right adjoint as an $\mathcal{A}$-$\mathcal{A}$-bimodule 1-morphism. These objects have been extensively studied in \cite{D7}, and should be thought of as generalizing the notion of a multifusion 1-category, which is recovered by taking $\mathfrak{C}=\mathbf{2Vect}$. Most of the algebras that we will encounter are also going to be \textit{connected}; by this we mean that their unit 1-morphism is a simple 1-morphism in $\mathfrak{C}$. When $\mathfrak{C}=\mathbf{2Vect}$, connected rigid algebras are precisely fusion 1-categories.

Given an algebra $\mathcal{A}$ in the fusion 2-category $\mathfrak{C}$, one can consider the 2-category $\mathbf{Bimod}_{\mathfrak{C}}(\mathcal{A})$ of $\mathcal{A}$-$\mathcal{A}$-bimodules in $\mathfrak{C}$. If $\mathcal{A}$ is rigid, it follows from work of the first author that $\mathbf{Bimod}_{\mathfrak{C}}(\mathcal{A})$ is a multifusion 2-category. More precisely, this last statement is a combination of theorem 3.1.6 of \cite{D7}, theorem 5.3.2 of \cite{D8}, and theorem 5.1.1 of \cite{D9}. In order to ensure that $\mathbf{Bimod}_{\mathfrak{C}}(\mathcal{A})$ is a fusion 2-category, it is sufficient to require that $\mathcal{A}$ be connected in addition of being rigid. By definition, the fusion 2-category $\mathbf{Bimod}_{\mathfrak{C}}(\mathcal{A})$ is Morita equivalent to $\mathfrak{C}$.\footnote{More generally, one can define the notion of Morita equivalence between multifusion 2-categories using arbitrary rigid algebras. We refer the reader to the above references for detail.} This Morita equivalence is also witnessed by the finite semisimple $\mathfrak{C}$-module 2-category $\mathbf{Mod}_{\mathfrak{C}}(\mathcal{A})$ (the definition of a module 2-category is given in section 2.1 of \cite{D4}). For later use, it will be practical to have another formulation of this concept. In order to do so, we introduce the following definition (see \cite{ENO2} for the decategorified analogue).

\begin{Definition}
Let $\mathfrak{C}$ and $\mathfrak{D}$ be two fusion 2-categories. A finite semisimple $\mathfrak{C}$-$\mathfrak{D}$-bimodule 2-category $\mathfrak{M}$ is called invertible provided that the canonical monoidal 2-functor from $\mathfrak{D}^{mop}$ to $\mathbf{End}_{\mathfrak{C}}(\mathfrak{M})$, the monoidal 2-category of left $\mathfrak{C}$-module endo-2-functors on $\mathfrak{M}$, is an equivalence.
\end{Definition}

\noindent It follows from theorem 5.4.3 of \cite{D8} that a $\mathfrak{C}$-$\mathfrak{D}$-bimodule 2-category $\mathfrak{M}$ is invertible if and only if it induces a Morita equivalence between $\mathfrak{C}$ and $\mathfrak{D}$. In particular, given two connected rigid algebras $\mathcal{A}$ and $\mathcal{B}$ in $\mathfrak{C}$, $\mathbf{Bimod}_{\mathfrak{C}}(\mathcal{A}, \mathcal{B})$ is an invertible $\mathbf{Bimod}_{\mathfrak{C}}(\mathcal{A})$-$\mathbf{Bimod}_{\mathfrak{C}}(\mathcal{B})$-bimodule 2-category.

Further, it follows from theorem 5.4.3 of \cite{D8} that invertible bimodule 2-categories can be composed, i.e.\ stacked.\footnote{More generally, it follows from \cite{D9} that the relative 2-Deligne tensor product of any finite semisimple bimodule 2-categories exists.} For instance, any monoidal autoequivalence $\mathbf{E}:\mathfrak{D}\simeq \mathfrak{D}$ induces an invertible bimodule 2-category: We can consider the $\mathfrak{D}$-$\mathfrak{D}$-bimodule 2-category $\mathfrak{D}_{\mathbf{E}}$, whose right $\mathfrak{D}$-module structure is twisted by $\mathbf{E}$. Such invertible bimodule 2-categories compose in the obvious way. Moreover, the composite of $\mathfrak{M}$ with $\mathfrak{D}_{\mathbf{E}}$ is given by $\mathfrak{M}_{\mathbf{E}}$, the bimodule 2-category obtained by twisting the right action of $\mathfrak{D}$ by $\mathbf{E}$.

\subsection{Algebras and Bimodules in \texorpdfstring{$\mathbf{2Vect}_G^{\pi}$}{Twisted Group-Graded 2-Vector Spaces}}\label{ap:detaileddefinitions}

We fix a finite group $G$ with unit $e$, and $\pi$ a (normalized) 4-cocycle for $G$ with coefficients in $\mathds{k}^{\times}$. We will spell out the definitions of algebras, modules, and bimodules in the fusion 2-category $\mathbf{2Vect}_G^{\pi}$. More precisely, we specialize the definitions given in sections 3.1 and 3.2 of \cite{D4} as well as those of sections 1.2 and 1.3 of \cite{D7} for a general monoidal 2-category to the case of $\mathbf{2Vect}_G^{\pi}$. This will be needed in order to establish the proof of our main proposition in the next section.

\begin{Definition}
An algebra in $\mathbf{2Vect}_G^{\pi}$ is a $G$-graded finite semisimple 1-category $\mathcal{C} = \oplus_{g\in G} \mathcal{C}_g$ equipped with a grading preserving product $\otimes : \mathcal{C}\boxtimes\mathcal{C}\rightarrow\mathcal{C}$, a distinguished object $I\in\mathcal{C}_e$, and, for every objects $X,Y,Z$ in $\mathcal{C}$, natural isomorphisms $$\alpha_{X,Y,Z}:X\otimes (Y\otimes Z)\xrightarrow{\sim}(X\otimes Y)\otimes Z,$$ $$\lambda_Y:I\otimes Y\xrightarrow{\sim} Y,\ \ \ \rho_X:X\otimes I\xrightarrow{\sim} X.$$ For any $g_1,g_2,g_3,g_4\in G$, the above natural isomorphisms have to satisfy for every $W$ in $\mathcal{C}_{g_1}$, $X$ in $\mathcal{C}_{g_2}$, $Y$ in $\mathcal{C}_{g_3}$, and $Z$ in $\mathcal{C}_{g_4}$, the following equations: 
$$\begin{tikzcd}[sep=small]
W\otimes(X\otimes (Y\otimes Z)) \arrow[rr, "{\overset{\alpha_{W,X,Y\otimes Z}}{\phantom{a}}}"] \arrow[dd, "{W\otimes\alpha_{X,Y,Z}}"'] &  & (W\otimes X)\otimes (Y\otimes Z) \arrow[rr, "{\overset{\alpha_{W\otimes X, Y,Z}}{\phantom{a}}}"] &  & ((W\otimes X)\otimes Y)\otimes Z \arrow[dd, "{\pi(g_1,g_2,g_3,g_4)}"] \\
 &  &  &  &  \\
W\otimes ((X\otimes Y)\otimes Z) \arrow[rr, "{\overset{\phantom{a}}{\alpha_{W,X\otimes Y,Z}}}"']                                       &  & (W\otimes (X\otimes Y))\otimes Z \arrow[rr, "{\overset{\phantom{a}}{{\alpha_{W,X,Y}\otimes Z}}}"'] &  & ((W\otimes X)\otimes Y)\otimes Z,                                     
\end{tikzcd}$$

$$\begin{tikzcd}
X\otimes (I\otimes Y) \arrow[rr, "{\alpha_{X,I,Y}}"] \arrow[rd, "X\otimes\lambda_Y"'] &            & (X\otimes I)\otimes Y \arrow[ld, "\rho_X\otimes Y"] \\
 & X\otimes Y. & 
\end{tikzcd}$$
\end{Definition}

\begin{Remark}
When $\pi$ is trivial, algebras in $\mathbf{2Vect}_G$ are exactly $G$-graded monoidal finite semisimple 1-categories. In general, we could use the term $\pi$-twisted $G$-graded monoidal finite semisimple 1-categories to talk about the algebras in $\mathbf{2Vect}_G^{\pi}$ .
\end{Remark}

\begin{Example}
Let $H\subseteq G$ be a subgroup. Given a 3-cochain $\psi$ for $H$ with coefficients in $\mathds{k}^{\times}$ such that $d\psi = \pi|_H$, we can then consider the algebra $\mathbf{Vect}_H^{\psi}$ in $\mathbf{2Vect}_G^{\pi}$ with the tautological $G$-grading, and with $\alpha_{\mathds{k}_{h_1}, \mathds{k}_{h_2}, \mathds{k}_{h_3}} = \psi(h_1,h_2,h_3)$ for every $h_1,h_2,h_3\in H$. If $\psi$ is normalized, which we can always assume, then we can take $\lambda$ and $\rho$ to be the identity natural isomorphisms.
\end{Example}

Let us now fix an algebra $\mathcal{C}$ in $\mathbf{2Vect}_G^{\pi}$.

\begin{Definition}
A left $\mathcal{C}$-module in $\mathbf{2Vect}_G^{\pi}$ is a $G$-graded finite semisimple 1-category $\mathcal{M} = \oplus_{g\in G} \mathcal{M}_g$ equipped with a grading preserving action $\otimes : \mathcal{C}\boxtimes\mathcal{M}\rightarrow\mathcal{M}$, and for every object $X,Y$ in $\mathcal{C}$ and $M$ in $\mathcal{M}$, a natural isomorphism $$\kappa_{X,Y,M}:X\otimes (Y\otimes M)\xrightarrow{\sim}(X\otimes Y)\otimes M,\ \ \ \lambda_M:I\otimes M\xrightarrow{\sim} M.$$ For any $g_1,g_2,g_3,g_4\in G$, the above natural isomorphisms have to satisfy for every $X$ in $\mathcal{C}_{g_1}$, $Y$ in $\mathcal{C}_{g_2}$, $Z$ in $\mathcal{C}_{g_3}$, and $M$ in $\mathcal{M}_{g_4}$, the following equations:
$$\begin{tikzcd}[sep=small]
X\otimes(Y\otimes (Z\otimes M)) \arrow[rr, "{\overset{\kappa_{X,Y,Z\otimes M}}{\phantom{a}}}"] \arrow[dd, "{X\otimes\kappa_{Y,Z,M}}"'] &  & (X\otimes Y)\otimes (Z\otimes M) \arrow[rr, "{\overset{\kappa_{X\otimes Y, Z,M}}{\phantom{a}}}"] &  & ((X\otimes Y)\otimes Z)\otimes M \arrow[dd, "{\pi(g_1,g_2,g_3,g_4)}"] \\
 &  &  &  &  \\
X\otimes ((Y\otimes Z)\otimes M) \arrow[rr, "{\overset{\phantom{a}}{\kappa_{X,Y\otimes Z,M}}}"']                                       &  & (X\otimes (Y\otimes Z))\otimes M \arrow[rr, "{\overset{\phantom{a}}{\alpha_{X,Y,Z}\otimes M}}"'] &  & ((X\otimes Y)\otimes Z)\otimes M,              
\end{tikzcd}$$

$$\begin{tikzcd}
X\otimes (I\otimes M) \arrow[rr, "{\kappa_{X,I,M}}"] \arrow[rd, "X\otimes\lambda_M"'] &            & (X\otimes I)\otimes M \arrow[ld, "\rho_X\otimes M"] \\
 & X\otimes M. & 
\end{tikzcd}$$
\end{Definition}

\begin{Definition}
Let $\mathcal{M}$ and $\mathcal{N}$ be two left $\mathcal{C}$-modules in $\mathbf{2Vect}_G^{\pi}$. A left $\mathcal{C}$-module 1-morphism $\mathcal{M}\rightarrow \mathcal{N}$ is a grading preserving functor $F:\mathcal{M}\rightarrow \mathcal{N}$ equipped with a natural isomorphism given on the objects $X$ in $\mathcal{C}$ and $M$ in $\mathcal{M}$ by $$\chi_{X,M}:X\otimes F(M)\xrightarrow{\sim} F(X\otimes M).$$ This natural isomorphism satisfies the usual equation, that is, for every $X,Y$ in $\mathcal{C}$, and $M$ in $\mathcal{M}$, we have:
$$\begin{tikzcd}[sep=small]
X\otimes (Y\otimes F(M)) \arrow[d, "{X\otimes \chi_{Y,M}}"'] \arrow[rr, "{\kappa_{X,Y,F(M)}}"] &                                                            & (X\otimes Y)\otimes F(M) \arrow[d, "{\phantom{a}\chi_{X\otimes Y,M}}"] \\
X\otimes F(Y\otimes M) \arrow[rd, "{\chi_{X,Y\otimes M}}"']                                    &                                                            & F((X\otimes Y)\otimes M)                                    \\
& F(X\otimes(Y\otimes M)). \arrow[ru, "{F(\kappa_{X,Y,M})}"'] &                     
\end{tikzcd}$$

A left $\mathcal{C}$-module 2-morphism in $\mathbf{2Vect}_G^{\pi}$ between two left $\mathcal{C}$-module 1-morphisms $F_1,F_2:\mathcal{M}\rightarrow\mathcal{N}$ is a natural transformation between the functors $F_1$ and $F_2$ that is compatible with the left $\mathcal{C}$-module structures in the obvious way.
\end{Definition}

Let us now fix another algebra $\mathcal{D}$ in $\mathbf{2Vect}_G^{\pi}$.

\begin{Definition}
A right $\mathcal{D}$-module in $\mathbf{2Vect}_G^{\pi}$ is a $G$-graded finite semisimple 1-category $\mathcal{M} = \oplus_{g\in G} \mathcal{M}_g$ equipped with a grading preserving action $\otimes : \mathcal{M}\boxtimes\mathcal{D}\rightarrow\mathcal{M}$, and for every object  $M$ in $\mathcal{M}$ and $X,Y$ in $\mathcal{D}$, a natural isomorphism $$\nu_{M,X,Y}:M\otimes (X\otimes Y)\xrightarrow{\sim}(M\otimes X)\otimes Y,\ \ \ \rho_M:M\otimes I\xrightarrow{\sim} M.$$ For any $g_1,g_2,g_3,g_4\in G$, the above natural isomorphisms have to satisfy for every $M$ in $\mathcal{M}_{g_1}$, $X$ in $\mathcal{D}_{g_2}$, $Y$ in $\mathcal{D}_{g_3}$, and $Z$ in $\mathcal{D}_{g_4}$, the following equations: 
$$\begin{tikzcd}[sep=small]
M\otimes(X\otimes (Y\otimes Z)) \arrow[rr, "{\overset{\nu_{M,X,Y\otimes Z}}{\phantom{a}}}"] \arrow[dd, "{M\otimes\alpha_{X,Y,Z}}"'] &  & (M\otimes X)\otimes (Y\otimes Z) \arrow[rr, "{\overset{\nu_{M\otimes X, Y,Z}}{\phantom{a}}}"] &  & ((M\otimes X)\otimes Y)\otimes Z \arrow[dd, "{\pi(g_1,g_2,g_3,g_4)}"] \\
 &  &  &  &  \\
M\otimes ((X\otimes Y)\otimes Z) \arrow[rr, "{\overset{\phantom{a}}{\nu_{M,X\otimes Y,Z}}}"']                                       &  & (M\otimes (X\otimes Y))\otimes Z \arrow[rr, "{\overset{\phantom{a}}{\nu_{M,X,Y}\otimes Z}}"'] &  & ((M\otimes X)\otimes Y)\otimes Z,              
\end{tikzcd}$$

$$\begin{tikzcd}
M\otimes (I\otimes Y) \arrow[rr, "{\nu_{M,I,Y}}"] \arrow[rd, "M\otimes\lambda_Y"'] &            & (M\otimes I)\otimes Y \arrow[ld, "\rho_M\otimes Y"] \\
 & M\otimes Y. & 
\end{tikzcd}$$
\end{Definition}

\begin{Definition}
Let $\mathcal{M}$ and $\mathcal{N}$ be two right $\mathcal{D}$-modules in $\mathbf{2Vect}_G^{\pi}$. A right $\mathcal{D}$-module 1-morphism $\mathcal{M}\rightarrow \mathcal{N}$ is a grading preserving functor $F:\mathcal{M}\rightarrow \mathcal{N}$ equipped with a natural isomorphism given on the objects $M$ in $\mathcal{M}$ and $X$ in $\mathcal{D}$ by $$\phi_{M,X}:F(M)\otimes X\xrightarrow{\sim} F(M\otimes X).$$ This natural isomorphism satisfies the usual equation, that is, for every $M$ in $\mathcal{M}$, and $X,Y$ in $\mathcal{D}$, we have:
$$\begin{tikzcd}[sep=small]
F(M)\otimes (X\otimes Y) \arrow[d, "{\phi_{M,X\otimes Y}}"'] \arrow[rr, "{\nu_{F(M),X,Y}}"] &                           & (F(M)\otimes X)\otimes Y \arrow[d, "{\phi_{M,X}\otimes Y}"] \\
F(M\otimes (X\otimes Y)) \arrow[rd, "{F(\nu_{M,X,Y})}"']                                    &                           & F(M\otimes X)\otimes Y \arrow[ld, "{\phi_{M\otimes X, Y}}"] \\
 & F((M\otimes Y) \otimes M). &        
\end{tikzcd}$$

A right $\mathcal{D}$-module 2-morphism in $\mathbf{2Vect}_G^{\pi}$ between two right $\mathcal{D}$-module 1-morphisms $F_1,F_2:\mathcal{M}\rightarrow\mathcal{N}$ is a natural transformation between the functors $F_1$ and $F_2$ that is compatible with the right $\mathcal{D}$-module structures in the obvious way.
\end{Definition}

\begin{Definition}\label{def:bimodule}
A $\mathcal{C}$-$\mathcal{D}$-bimodule in $\mathbf{2Vect}_G^{\pi}$ is a $G$-graded finite semisimple 1-category $\mathcal{P} = \oplus_{g\in G} \mathcal{P}_g$ equipped with both a left $\mathcal{C}$-module structure and a right $\mathcal{D}$-module structure, and for every object $X$ in $\mathcal{C}$, $M$ in $\mathcal{M}$, and $Y\in\mathcal{D}$ a natural isomorphism $$\beta_{X,M,Y}:X\otimes (M\otimes Y)\xrightarrow{\sim}(X\otimes M)\otimes Y.$$ Further, for any $g_1,g_2,g_3,g_4,g_5\in G$, the above natural isomorphism has to satisfy for every $W$ in $\mathcal{C}_{g_1}$, $X$ in $\mathcal{C}_{g_2}$, $M$ in $\mathcal{M}_{g_3}$, $Y$ in $\mathcal{D}_{g_4}$, and $Z$ in $\mathcal{D}_{g_5}$, the following two equations: 
$$\begin{tikzcd}[sep=small]
W\otimes(X\otimes (M\otimes Y)) \arrow[rr, "{\overset{\kappa_{W,X,M\otimes Y}}{\phantom{a}}}"] \arrow[dd, "{W\otimes\beta_{X,M,Y}}"'] &  & (W\otimes X)\otimes (M\otimes Y) \arrow[rr, "{\overset{\beta_{W\otimes X, M,Y}}{\phantom{a}}}"] &  & ((W\otimes X)\otimes M)\otimes Y \arrow[dd, "{\pi(g_1,g_2,g_3,g_4)}"] \\
 &  &  &  &  \\
W\otimes ((X\otimes M)\otimes Y) \arrow[rr, "{\overset{\phantom{a}}{\beta_{W,X\otimes M,Y}}}"']                                       &  & (W\otimes (X\otimes M))\otimes Y \arrow[rr, "{\overset{\phantom{a}}{\kappa_{W,X,M}\otimes Y}}"'] &  & ((W\otimes X)\otimes M)\otimes Y,              
\end{tikzcd}$$
\vspace{3mm}
$$\begin{tikzcd}[sep=small]
X\otimes(M\otimes (Y\otimes Z)) \arrow[rr, "{\overset{\beta_{X,M,Y\otimes Z}}{\phantom{a}}}"] \arrow[dd, "{X\otimes\nu_{M,Y,Z}}"'] &  & (X\otimes M)\otimes (Y\otimes Z) \arrow[rr, "{\overset{\nu_{X\otimes M, Y,Z}}{\phantom{a}}}"] &  & ((X\otimes M)\otimes Y)\otimes Z \arrow[dd, "{\pi(g_2,g_3,g_4,g_5)}"] \\
 &  &  &  &  \\
X\otimes ((M\otimes Y)\otimes Z) \arrow[rr, "{\overset{\phantom{a}}{\beta_{X,M\otimes Y,Z}}}"']                                       &  & (X\otimes (M\otimes Y))\otimes Z \arrow[rr, "{\overset{\phantom{a}}{\beta_{X,M,Y}\otimes Z}}"'] &  & ((X\otimes M)\otimes Y)\otimes Z.              
\end{tikzcd}$$
\end{Definition}

\begin{Definition}\label{def:bimodulemorphism}
Let $\mathcal{P}$ and $\mathcal{Q}$ be two $\mathcal{C}$-$\mathcal{D}$-bimodules in $\mathbf{2Vect}_G^{\pi}$. A $\mathcal{C}$-$\mathcal{D}$-bimodule 1-morphism $\mathcal{P}\rightarrow \mathcal{Q}$ is a grading preserving functor $F:\mathcal{P}\rightarrow \mathcal{Q}$ equipped with both a left $\mathcal{C}$-module structure and a right $\mathcal{D}$-module structure, which are compatible in the obvious way.

A $\mathcal{C}$-$\mathcal{D}$-bimodule 2-morphism in $\mathbf{2Vect}_G^{\pi}$ between two right $\mathcal{C}$-$\mathcal{D}$-bimodule 1-morphisms $F_1,F_2:\mathcal{P}\rightarrow\mathcal{Q}$ is a natural transformation between the functors $F_1$ and $F_2$ that is both compatible with the left $\mathcal{C}$-module structures, and the right $\mathcal{D}$-module structures.
\end{Definition}

\section{Group-Theoretical Fusion 2-Categories}\label{section:groupF2C}

\subsection{Definition and Characterization}

Group-theoretical fusion 1-categories were implicitly introduced in \cite{O1}, and later explicitly in \cite{ENO1}. We study a strong categorification of this notion.
The first specific examples of such fusion 2-categories were considered in the Physics literature 
\cite{BBSNT1, BBFP, BBSNT}. In this context, the term group-theoretical is used to describe symmetry topological field theories that are gauge equivalent to a Dijkgraaf-Witten theory \cite[Section 1.2]{SunZheng}.

\begin{Definition}
Let $G$ be a finite group, $\pi\in Z^4_{gp}(G;\mathds{k}^{\times})$, $H\subseteq G$ a subgroup, and $\psi\in C^3_{gp}(H;\mathds{k}^{\times})$ a 3-cochain such that $d\psi = \pi|_H$. Then, we set $$\mathfrak{C}(G,H,\pi, \psi):=\mathbf{Bimod}_{\mathbf{2Vect}_G^{\pi}}(\mathbf{Vect}_H^{\psi}).$$ We say that a fusion 2-category is \textit{group-theoretical} if it is monoidally equivalent to a fusion 2-category of this form. 
\end{Definition}

\begin{Example}
Let $G$ be a finite group, and $\pi$ a 4-cocycle for $G$. We have $$\mathfrak{C}(G,\{e\},\pi, triv) = \mathbf{2Vect}_G^{\pi}.$$ Moreover, it follows from example 5.1.9 of \cite{D8} that $$\mathfrak{C}(G,G,triv, triv) \simeq \mathbf{2Rep}(G).$$
\end{Example}

Our objective is now to establish an intrinsic characterization of the property of being group-theoretical. We note that there is no similar criterion to detect group-theoretical fusion 1-categories.

\begin{Theorem}\label{thm:grouptheoreticalrecognition}
A fusion 2-category $\mathfrak{C}$ is group-theoretical if and only if there is an equivalence of braided fusion 1-categories $\Omega\mathfrak{C}\simeq \mathbf{Rep}(H)$ for some finite group $H$.
\end{Theorem}

The proof of the above theorem will be achieved using ideas from \cite{D9}, but, first, we prove some technical results.

\begin{Lemma}\label{lem:rigidalgebragrouptheoretical}
Let $\mathcal{A}$ be a connected rigid algebra in $\mathbf{2Vect}_G^{\pi}$. Then, the finite semisimple 2-category $\mathbf{Mod}_{\mathbf{2Vect}_G^{\pi}}(\mathcal{A})$ is of the form $\boxplus_{i=1}^n\mathbf{2Vect}$ if and only if $\mathcal{A}$ is equivalent to an algebra of the form $\mathbf{Vect}_H^{\psi}$ for some subgroup $H\subseteq G$ and $\psi\in C^3_{gp}(H;\mathds{k}^{\times})$ such that $d\psi = \pi|_H$.
\end{Lemma}
\begin{proof}
Let $\mathcal{A}$ be a connected rigid algebra in $\mathbf{2Vect}_G^{\pi}$ with multiplication 1-morphism $m:\mathcal{A}\boxtimes \mathcal{A}\rightarrow \mathcal{A}$. 
Let us assume that $\mathbf{Mod}_{\mathbf{2Vect}_G^{\pi}}(\mathcal{A})\simeq\boxplus_{i=1}^n\mathbf{2Vect}$. For any $g\in G$, let us write $\mathbf{Vect}_g$ for the corresponding simple object in $\mathbf{2Vect}_G^{\pi}$. Note that, as $\mathcal{A}$ is connected, the right $\mathcal{A}$-module $\mathbf{Vect}_{g}\boxtimes \mathcal{A}$ is simple for any $g$. Now, by lemma 3.2.13 of \cite{D4}, we have that \begin{equation}\label{eq:technicalgrouptheoretical}Hom_{\mathcal{A}}(\mathbf{Vect}_{g}\boxtimes \mathcal{A}, \mathcal{A})\simeq Hom_{\mathbf{2Vect}_G^{\pi}}(\mathbf{Vect}_{g}, \mathcal{A}).\end{equation} But, $\mathbf{Mod}_{\mathbf{2Vect}_G^{\pi}}(\mathcal{A})\simeq\boxplus_{i=1}^n\mathbf{2Vect}$ by hypothesis, so that the left hand side is either $0$ of $\mathbf{Vect}$. This shows that the underlying object of $\mathcal{A}$ in $\mathbf{2Vect}_G^{\pi}$ is of the form $\boxplus_{g\in H}\mathbf{Vect}_g$ for some subset $H\subseteq G$.

Let $g,h\in H$ with corresponding simple 1-morphisms $i_g:\mathbf{Vect}_g\rightarrow \mathcal{A}$ and $i_h:\mathbf{Vect}_h\rightarrow \mathcal{A}$. We claim that the composite $m\circ (i_g\boxtimes i_h):\mathbf{Vect}_g\boxtimes \mathbf{Vect}_h\rightarrow \mathcal{A}$ is non-zero. Namely, the two corresponding right $\mathcal{A}$-module 1-morphisms under the equivalence of equation \eqref{eq:technicalgrouptheoretical} are non-zero and given by $F_g:=m\circ (i_g\boxtimes \mathcal{A})$ and $F_h:=m\circ (i_h\boxtimes \mathcal{A})$. Then, by lemma 2.2.11 of \cite{D4}, we find that the composite $$F_g\circ (\mathbf{Vect}_g\boxtimes F_h): \mathbf{Vect}_g\boxtimes\mathbf{Vect}_h\boxtimes \mathcal{A}\rightarrow \mathcal{A}$$ is non-zero. But, the corresponding 1-morphism $\mathbf{Vect}_g\boxtimes \mathbf{Vect}_h\rightarrow \mathcal{A}$ under equation \eqref{eq:technicalgrouptheoretical} is $m\circ (i_g\boxtimes i_h)$, which is therefore non-zero as claimed.

The preceding paragraph not only shows that $H$ is a subgroup of $G$, but also that the multiplication 1-morphism $m:\mathcal{A}\boxtimes \mathcal{A}\rightarrow \mathcal{A}$ is given by the group structure of $H$, so that $\mathcal{A}\simeq \mathbf{Vect}_H^{\psi}$ for some subgroup $H\subseteq G$ and $\psi\in C^3_{gp}(H;\mathds{k}^{\times})$ such that $d\psi = \pi|_H$. This concludes the forward direction of proof. Conversely, the backward direction follows from \eqref{eq:technicalgrouptheoretical}, so that the proof of the lemma is complete.
\end{proof}

\begin{Corollary}\label{cor:grouptheoreticaltechnical}
A fusion 2-category is group-theoretical if and only if it is Morita equivalent to $\mathbf{2Vect}_G^{\pi}$ via a finite semisimple module 2-category $\mathfrak{M}$, whose underlying 2-category is $\boxplus_{i=1}^n\mathbf{2Vect}$.
\end{Corollary}
\begin{proof}
The forward direction is immediate by the previous lemma. Conversely, let $\mathfrak{C}$ be a fusion 2-category and $\mathfrak{M}$ be an invertible $\mathbf{2Vect}_G^{\pi}$-$\mathfrak{C}$-bimodule 2-category such that $\mathfrak{M}\simeq \boxplus_{i=1}^n\mathbf{2Vect}$ as finite semisimple 2-categories. As $\mathfrak{C}$ is fusion, $\mathfrak{M}$ is indecomposable as a left $\mathbf{2Vect}_G^{\pi}$-module 2-category by corollary 5.2.5 of \cite{D8}. Invoking theorem 5.2.9 of \cite{D4}, there exists a connected rigid algebra $\mathcal{A}$ in $\mathbf{2Vect}_G^{\pi}$ such that $\mathfrak{M}\simeq\mathbf{Mod}_{\mathbf{2Vect}_G^{\pi}}(\mathcal{A})$. It follows from the previous lemma that, as algebras in $\mathbf{2Vect}_G^{\pi}$, we must have $\mathcal{A}\simeq \mathbf{Vect}_H^{\psi}$ for some subgroup $H\subseteq G$ and 3-cochain $\psi$ for $H$ with $d\psi = \pi|_H$, so that $\mathfrak{C}$ is group-theoretical as claimed.
\end{proof}

\begin{proof}[Proof of Thm. \ref{thm:grouptheoreticalrecognition}]
Let us assume that $\mathfrak{C}$ is group-theoretical. It follows from the forward direction of corollary \ref{cor:grouptheoreticaltechnical} that $\mathbf{2Vect}$ is a module 2-category for $\mathfrak{C}^0$, the connected component of the identity of $\mathfrak{C}$. The data of such a module 2-category structure is equivalent to that of a braided monoidal functor $\Omega\mathfrak{C}\rightarrow \mathbf{Vect}$. By the main theorem of \cite{Del}, this proves that $\Omega\mathfrak{C}\simeq\mathbf{Rep}(H)$ for some finite group $H$, and therefore gives the desired result in this direction.

Conversely, let us assume that $\Omega\mathfrak{C}\simeq \mathbf{Rep}(H)$ as braided fusion 1-categories. We can consider $\mathbf{Vect}$ as a connected rigid algebra in $\mathbf{Mod}(\mathbf{Rep}(H))$ via the canonical symmetric monoidal functor $\mathbf{Rep}(H)\rightarrow \mathbf{Vect}$. As $\mathfrak{C}^0$ is equivalent to the fusion 2-category $\mathbf{Mod}(\mathbf{Rep}(H))$, we can therefore also view $\mathbf{Vect}$ as an algebra in $\mathfrak{C}$. Then, we have equivalences of braided fusion 1-categories $$\mathbf{Bimod}_{\mathfrak{C}}(\mathbf{Vect})^0\simeq \mathbf{Bimod}_{\mathbf{Mod}(\mathbf{Rep}(H))}(\mathbf{Vect})^0\simeq (\mathbf{2Vect}_H)^0\simeq\mathbf{2Vect},$$ where we have used the equivalence $\mathbf{Bimod}_{\mathbf{Mod}(\mathbf{Rep}(H))}(\mathbf{Vect})\simeq \mathbf{2Vect}_H$ of fusion 2-categories from example 5.1.9 of \cite{D8}. But, as $\mathbf{Bimod}_{\mathfrak{C}}(\mathbf{Vect})^0\simeq \mathbf{Vect}$, it follows from \cite{JFY} that there is an equivalence of fusion 2-categories $\mathbf{Bimod}_{\mathfrak{C}}(\mathbf{Vect})\simeq \mathbf{2Vect}_G^{\pi}$ for some finite group $G$ and 4-cocycle $\pi$. Moreover, by construction, the finite semisimple 2-category $\mathfrak{M}:=\mathbf{Mod}_{\mathfrak{C}}(\mathbf{Vect})$ contains $\mathbf{Mod}_{\mathbf{2Rep}(H)}(\mathbf{Vect})\simeq \mathbf{2Vect}$ as a connected summand. This implies that $\mathfrak{M}\simeq \boxplus_{i=1}^n\mathbf{2Vect}$ as $\mathfrak{M}$ is an indecomposable right module 2-category over $\mathbf{Bimod}_{\mathfrak{C}}(\mathbf{Vect})\simeq \mathbf{2Vect}_G^{\pi}$, and therefore all the connected summands of $\mathfrak{M}$ must be equivalent. Consequently, all the conditions of the backwards direction of corollary \ref{cor:grouptheoreticaltechnical} are satisfied by $\mathfrak{C}$ and $\mathfrak{M}$, so that $\mathfrak{C}$ is group-theoretical. This concludes the proof.
\end{proof}

\begin{Example}
Many known fusion 2-categories satisfy the criterion of theorem \ref{thm:grouptheoreticalrecognition}. For instance, given any finite 2-group $\mathcal{G}= A[1]\boldsymbol{\cdot}H[0]$, we have $\Omega\mathbf{2Rep}(\mathcal{G}) = \mathbf{Rep}(H)$ and $\Omega\mathbf{2Vect}_{\mathcal{G}} = \mathbf{Vect}_A$. In particular, both $\mathbf{2Rep}(\mathcal{G})$ and $\mathbf{2Vect}_{\mathcal{G}}$ are group-theoretical. Given $\pi\in Z^4(\mathrm{B}\cG;\mathds{k}^{\times})$, the fusion 2-category $\mathbf{2Vect}_{\cG}^{\pi}$ might not be group-theoretical. Namely, there is an equivalence of braided fusion 1-categories $\Omega\mathbf{2Vect}_{\cG}^{\pi}\simeq \mathbf{Vect}_A^{f^*\pi}$, where $f^*\pi$ is the class in $H^4(\mathrm{B}^2A;\mathds{k})$ given by the pullback of $\pi$ along the inclusion $f:\mathrm{B}^2A\rightarrow \mathrm{B}\cG$. We find that $\mathbf{2Vect}_{\cG}^{\pi}$ is group-theoretical if and only if the class $f^*\pi$ is trivial.
\end{Example}

\begin{Example}\label{ex:center}
Given a finite group $G$ and $\pi\in Z^4_{gp}(\cG;\mathds{k}^{\times})$, the fusion 2-category $\mathscr{Z}(\mathbf{2Vect}_G^{\pi})$, which is the Drinfeld center of $\mathbf{2Vect}_G^{\pi}$, was studied in \cite{KTZ}, and is group-theoretical. Namely, it was shown in \cite{KTZ} that $\Omega\mathscr{Z}(\mathbf{2Vect}_G^{\pi}) = \mathbf{Rep}(G)$. In fact, if we write $\Delta\hookrightarrow G\times G^{op}$ for the diagonal subgroup, it follows from lemma 2.1.1 of \cite{D9} that we have $$\mathfrak{C}(G\times G^{op}, \Delta, \pi \times \pi^{op}, \psi)\simeq \mathscr{Z}(\mathbf{2Vect}^{\pi}_G)$$ an equivalence of fusion 2-categories, where $\pi^{op}(g_1,g_2,g_3,g_4):=\pi(g_4,g_3,g_2,g_1)$ for every $g_1,g_2,g_3,g_4\in G^{op}$, and some 3-cochain $\psi$.
\end{Example}

As another immediate corollary of theorem \ref{thm:grouptheoreticalrecognition}, we also obtain the following result, which includes a categorification of proposition 8.44 of \cite{ENO1}.

\begin{Corollary}
Let $F:\mathfrak{C}\rightarrow \mathfrak{D}$ be a monoidal 2-functor between two fusion 2-categories. If $\mathfrak{D}$ is group-theoretical, then so is $\mathfrak{C}$.
\end{Corollary}
\begin{proof}
The variant of Deligne's theorem given in theorem 9.9.22 of \cite{EGNO} asserts that a braided fusion 1-category $\mathcal{B}$ is equivalent to $\mathbf{Rep}(G)$ for some finite group $G$ if and only if there exists a braided monoidal functor $\mathcal{B}\rightarrow\mathbf{Vect}$. Now, the monoidal 2-functor $F$ induces a braided monoidal functor $\Omega F:\Omega\mathfrak{C}\rightarrow\Omega\mathfrak{D}$ (see proposition 2.4.7 of \cite{D2}). But, it follows from theorem \ref{thm:grouptheoreticalrecognition} and the first fact recalled above that there exists a braided monoidal functor $E:\Omega\mathfrak{D}\rightarrow\mathbf{Vect}$. Then the composite $E\circ \Omega F:\Omega\mathfrak{C}\rightarrow\mathbf{Vect}$ is a braided monoidal functor, and the variant of Deligne's theorem therefore implies that $\mathfrak{C}$ satisfies the criterion of theorem \ref{thm:grouptheoreticalrecognition}, which concludes the proof.
\end{proof}

\subsection{The Underlying 2-Category}

We now give an explicit description of the underlying 2-category of a group-theoretical fusion 2-category. In particular, we establish that the connected components of a group-theoretical fusion 2-category are given by double cosets. In fact, more generally, we describe the 2-category of bimodules over two algebras $\mathbf{Vect}_H^{\psi}$ and $\mathbf{Vect}_K^{\omega}$ in $\mathbf{2Vect}_G^{\pi}$. This categorifies proposition 3.1 of \cite{O2}. We also note that a version of the statement below with $H=K$ and $\psi=\omega$ has already appeared at a physical level of rigor in section 3 of \cite{BBFP}.

\begin{Proposition}\label{prop:2categorycosetdecomposition}
Let $G$ be a finite group, and $\pi\in Z^4_{gp}(G;\mathds{k}^{\times})$. Given $H,K\subseteq G$ two subgroups, and $\psi\in C^3_{gp}(H;\mathds{k}^{\times})$, $\omega\in C^3_{gp}(K;\mathds{k}^{\times})$ such that $d\psi = \pi|_H$ and $d\omega = \pi|_K$, then there is an equivalence of 2-categories $$\mathbf{Bimod}_{\mathbf{2Vect}_G^{\pi}}(\mathbf{Vect}_H^{\psi}, \mathbf{Vect}_K^{\omega})\simeq \boxplus_{\lbrack g\rbrack \in H\backslash G/ K} \mathbf{Mod}(\mathbf{Vect}_{H\cap gKg^{-1}}^{\xi_g}),$$ where the sum is taken over representatives in $G$ for the double cosets in $H\backslash G/ K$, and $\xi_g$ is the 3-cocycle for $H\cap gKg^{-1}$ given on $h_1,h_2,h_3\in H\cap gKg^{-1}$ by\footnote{Since the target is $\mathds{k}^{\times}$, we use multiplication to denote the combination of cochains.} \begin{align}\label{eq:xig}\xi_g(h_1,h_2&,h_3):=\frac{\psi(h_1,h_2,h_3)}{\omega(g^{-1}h_3^{-1}g,g^{-1}h_2^{-1}g,g^{-1}h_1^{-1}g)}\\
&\frac{\pi(h_3,h_3^{-1}g,g^{-1}h_2^{-1}g,g^{-1}h_1^{-1}g)\pi(g,g^{-1}h_3^{-1}g,g^{-1}h_2^{-1}g,g^{-1}h_1^{-1}g)}{\pi(h_2,h_3,h_3^{-1}h_2^{-1}g,g^{-1}h_1^{-1}g)\pi(h_1,h_2,h_3, h_3^{-1}h_2^{-1}h_1^{-1}g)}.\notag\end{align}
\end{Proposition}

We begin by recalling the following well-known lemma as well as its proof in a special case.

\begin{Lemma}\label{lem:freeforget}
Let $\mathfrak{C}$ be a monoidal 2-category, and be $\mathcal{A}$, $\mathcal{B}$ two algebras in $\mathfrak{C}$. For any object $V$ of $\mathfrak{C}$, and $\mathcal{A}$-$\mathcal{B}$-bimodule $\mathcal{P}$, there is an equivalence of 1-categories \begin{equation}\label{eq:freeforget}
    \begin{tabular}{r c c c} $U:$ & $Hom_{\mathcal{A}\mathrm{-}\mathcal{B}}(\mathcal{A}\Box C\Box \mathcal{B}, \mathcal{P})$ & $\xrightarrow{\sim}$ & $Hom_{\mathfrak{C}}(C,\mathcal{P}),$\\ & $f$ & $\mapsto$ & $f\circ (i^{\mathcal{A}}\Box C\Box i^{\mathcal{B}})$\end{tabular}
\end{equation} where $i^{\mathcal{A}}:I\rightarrow \mathcal{A}$ and $i^{\mathcal{B}}:I\rightarrow \mathcal{B}$ are the units of the algebras.
\end{Lemma}
\begin{proof}
For later use in the proof of proposition \ref{prop:2categorycosetdecomposition}, we need to spell out the pseudo-inverse. More precisely, for simplicity, we take $\mathfrak{C} = \mathbf{2Vect}_G^{\pi}$, and will therefore use the notations introduced in \S\ref{ap:detaileddefinitions} above.

Let us fix two subgroups $H\subseteq G$ and $K\subseteq G$ as well as two normalized 3-cochains $\psi$ for $H$ and $\omega$ for $K$ such that $d\psi = \pi|_H$ and $d\omega = \pi|_K$. We can then consider the two algebras $\mathcal{A} = \mathbf{Vect}_H^{\psi}$ and $\mathcal{B}=\mathbf{Vect}_K^{\omega}$. Further, we will only consider the case $C = \mathbf{Vect}_{g}$ for some fixed element $g\in G$.
Given any $\mathbf{Vect}_H^{\psi}$-$\mathbf{Vect}_K^{\omega}$-bimodule $\mathcal{P}$ in $\mathbf{2Vect}_G^{\pi}$, the pseudo-inverse to (\ref{eq:freeforget}) is given by sending the grading preserving functor $F:\mathbf{Vect}_{g}\rightarrow \mathcal{P}$ to the $\mathbf{Vect}_H^{\psi}$-$\mathbf{Vect}_K^{\omega}$-bimodule functor $\mathscr{F}(F):\mathbf{Vect}_H^{\psi}\boxtimes (\mathbf{Vect}_{g}\boxtimes \mathbf{Vect}_K^{\omega})\rightarrow \mathcal{P}$ given by $X\boxtimes V\boxtimes Y\mapsto X\otimes \big( F(V)\otimes Y\big)$, with $\mathbf{Vect}_H^{\psi}$-$\mathbf{Vect}_K^{\omega}$-bimodule structure inherited from $\mathcal{P}$. For legibility, given an element $x\in G$, we will simply write $x$ for the corresponding simple object of $\mathbf{Vect}_{x}$. For every $h_1,h_2\in H$, and $k_1,k_2\in K$, we have $$\chi^{\mathscr{F}(F)}_{h_1,h_2F(g)k_1}=\kappa_{h_1,h_2,F(g)k_1},$$ $$\phi^{\mathscr{F}(F)}_{h_2gk_1,k_2}=\beta^{-1}_{h_2,F(g)k_1,k_2}\cdot \nu^{-1}_{F(g),k_1,k_2}.$$ Further, the assignment $\mathscr{F}$ is functorial, and one verifies easily that it is indeed a pseudo-inverse for $U$. In particular, for any $\mathbf{Vect}_H^{\psi}$-$\mathbf{Vect}_K^{\omega}$-bimodule 1-morphism $E:\mathbf{Vect}_H^{\psi}\boxtimes (\mathbf{Vect}_{g}\boxtimes \mathbf{Vect}_K^{\omega})\rightarrow \mathcal{P}$, the counit is given by \begin{equation}\label{eq:counit}\epsilon_E:\mathscr{F}(U(E))\xrightarrow{\sim} E\ \ \ \text{with}\ \ \ (\epsilon_E)_{hgk}= \phi_{g,k}\cdot \chi_{h,gk}.\end{equation}
\end{proof}

\begin{proof}[Proof of Prop. \ref{prop:2categorycosetdecomposition}]
We start by introducing a particularly useful and managebale class of $\mathbf{Vect}_H^{\psi}$-$\mathbf{Vect}_K^{\omega}$-bimodule. Specifically, to any element $g\in G$, we associate a $\mathbf{Vect}_H^{\psi}$-$\mathbf{Vect}_K^{\omega}$-bimodule
\begin{equation}\label{eq:definitionVg}\mathbf{V}_{\lbrack g\rbrack} := \mathbf{Vect}_H^{\psi}\boxtimes \big(\mathbf{Vect}_g\boxtimes \mathbf{Vect}_K^{\omega}\big).\end{equation} For later use, we also spell out the $\mathbf{Vect}_H^{\psi}$-$\mathbf{Vect}_K^{\omega}$-bimodule structure constraints for $\mathbf{V}_{\lbrack g\rbrack}$ (using the notations of \S\ref{ap:detaileddefinitions}). Recall that, given an element $f\in G$, we will simply write $f$ for the simple object $\mathbf{Vect}_f$ of $\mathbf{2Vect}_G^{\pi}$. With this convention, we have that the bimodule structure on $\mathbf{V}_{\lbrack g\rbrack}=\mathbf{Vect}_H^{\psi}\boxtimes \mathbf{Vect}_{g}\boxtimes \mathbf{Vect}_K^{\omega}$ is given by \begin{align}\kappa_{h_1,h_2, h_3gk_1}=&  \psi(h_1, h_2, h_3)\cdot\pi^{-1}(h_1,h_2,h_3,gk_1),\notag \\ \nu_{h_3gk_1, k_2, k_3}=&  \omega(k_1,k_2,k_3)\cdot\pi^{-1}(h_3, gk_1,k_2,k_3)\cdot \pi^{-1}(g,k_1,k_2,k_3),\label{eq:explicitconstraints}\\ \beta_{h_2,h_3gk_1,k_2}=& \pi(h_2,h_3,gk_1,k_2),\notag \end{align} for any $h_1,h_2,h_3\in H$, and $k_1,k_2,k_3\in K$. 

We proceed to analyze the finite set of connected components of the finite semisimple 2-category $\mathbf{Bimod}_{\mathbf{2Vect}_G^{\pi}}(\mathbf{Vect}_H^{\psi}, \mathbf{Vect}_K^{\omega})$. To this end, $\mathcal{P}$ be an arbitrary simple $\mathbf{Vect}_H^{\psi}$-$\mathbf{Vect}_K^{\omega}$-bimodule. Given any simple summand $\mathbf{Vect}_g$ of $\mathcal{P}$ in $\mathbf{2Vect}_G^{\pi}$, there is a non-zero $\mathbf{Vect}_H^{\psi}$-$\mathbf{Vect}_K^{\omega}$-bimodule 1-morphism $\mathbf{V}_{\lbrack g\rbrack}\rightarrow \mathcal{P}$. Moreover, from lemma \ref{lem:freeforget}, we find
\begin{equation}\label{eq:freeforgetbody}End_{\mathbf{Vect}_H^{\psi}\mathrm{-}\mathbf{Vect}_K^{\omega}}(\mathbf{V}_{\lbrack g\rbrack})\simeq Hom_{\mathbf{2Vect}_G^{\pi}}(\mathbf{Vect}_g, \mathbf{V}_{\lbrack g\rbrack}),\end{equation} with the identity bimodule 1-morphism on
$\mathbf{V}_{\lbrack g\rbrack}$ corresponding to the canonical inclusion $\mathbf{Vect}_g\hookrightarrow \mathbf{Vect}_H^{\psi}\boxtimes\mathbf{Vect}_g\boxtimes \mathbf{Vect}_K^{\omega}$. In particular, $\mathbf{V}_{\lbrack g\rbrack}$ is a simple $\mathbf{Vect}_H^{\psi}$-$\mathbf{Vect}_K^{\omega}$-bimodule, and its equivalence class only depends on the double coset $\lbrack g\rbrack \in H\backslash G/ K$. This shows that $$\pi_0(\mathbf{Bimod}_{\mathbf{2Vect}_G^{\pi}}(\mathbf{Vect}_H^{\psi}, \mathbf{Vect}_K^{\omega}))\cong H\backslash G/ K.$$

We now begin to study the internal structure of the connected components of $\mathbf{Bimod}_{\mathbf{2Vect}_G^{\pi}}(\mathbf{Vect}_H^{\psi}, \mathbf{Vect}_K^{\omega})$. By the results of the previous paragraph, it is enough to describe $End_{\mathbf{Vect}_H^{\psi}\mathrm{-}\mathbf{Vect}_K^{\omega}}(\mathbf{V}_{\lbrack g\rbrack})$ for the double cosets $\lbrack g\rbrack \in H\backslash G/ K$. As a first observation, note that for any fixed $g$ in $G$, the simple 1-morphisms $\mathbf{Vect}_g\hookrightarrow\mathbf{Vect}_H\boxtimes\mathbf{Vect}_g\boxtimes \mathbf{Vect}_K$ in $\mathbf{2Vect}_G^{\pi}$ are parameterized by $H\cap g K g^{-1}$. More precisely, for every $h$ in $H\cap g K g^{-1}$, there corresponds a unique simple 1-morphism $$J_h:\mathbf{Vect}_g= \mathbf{Vect}_h\boxtimes\mathbf{Vect}_g\boxtimes \mathbf{Vect}_{(g^{-1}h^{-1}g)}\hookrightarrow \mathbf{Vect}_H\boxtimes\mathbf{Vect}_g\boxtimes \mathbf{Vect}_K.$$ Thence, it follows from the equivalence given in (\ref{eq:freeforgetbody}) that there is an equivalence \begin{equation}\label{eq:componentsfusionrules}End_{\mathbf{Vect}_H^{\psi}\mathrm{-}\mathbf{Vect}_K^{\omega}}(\mathbf{V}_{\lbrack g\rbrack})\simeq \mathbf{Vect}_{H\cap gKg^{-1}}\end{equation} of finite semisimple 1-categories, which is compatible with the fusion rules on both sides.

Finally, in order to complete our description of the finite semisimple 2-category $\mathbf{Bimod}_{\mathbf{2Vect}_G^{\pi}}(\mathbf{Vect}_H^{\psi}, \mathbf{Vect}_K^{\omega})$, it only remains to understand how to upgraded the equivalence given in equation \eqref{eq:componentsfusionrules} to a monoidal equivalence. More precisely, $End_{\mathbf{Vect}_H^{\psi}\mathrm{-}\mathbf{Vect}_K^{\omega}}(\mathbf{V}_{\lbrack g\rbrack})$ is a monoidal 1-category under composition, which is in fact strict. The monoidal structure on $End_{\mathbf{Vect}_H^{\psi}\mathrm{-}\mathbf{Vect}_K^{\omega}}(\mathbf{V}_{\lbrack g\rbrack})$ can be transferred via (\ref{eq:freeforget}) to endow $Hom_{\mathbf{2Vect}_G^{\pi}}(\mathbf{Vect}_g, \mathbf{V}_{\lbrack g\rbrack})\simeq \mathbf{Vect}_{H\cap gKg^{-1}}$ with a monoidal structure, which we denote by $\bigstar$. In more detail, given $h_1,h_2,h_3\in H\cap gKg^{-1}$, the monoidal product of $J_{h_1}$ and $J_{h_2}$ in $Hom_{\mathbf{2Vect}_G^{\pi}}(\mathbf{Vect}_g, \mathbf{V}_{\lbrack g\rbrack})$ is given by $$J_{h_1}\bigstar J_{h_2} := U(\mathscr{F}(J_{h_2})\circ \mathscr{F}(J_{h_1})))= J_{h_1h_2}.$$ Further, the associator of $Hom_{\mathbf{2Vect}_G^{\pi}}(\mathbf{Vect}_g, \mathbf{V}_{\lbrack g\rbrack})$, $$\alpha_{J_{h_1},J_{h_2},J_{h_3}}: J_{h_1}\bigstar(J_{h_2}\bigstar J_{h_3}) \cong (J_{h_1}\bigstar J_{h_2})\bigstar J_{h_3}$$ is given by the composite $$\alpha_{J_{h_1},J_{h_2},J_{h_3}}:U\Big(\mathscr{F}U\big(F(J_{h_3})\circ \mathscr{F}(J_{h_2})\big)\circ \mathscr{F}(J_{h_1})\Big)\xrightarrow{\epsilon}U\Big(\big(\mathscr{F}(J_{h_3})\circ \mathscr{F}(J_{h_2})\big)\circ \mathscr{F}(J_{h_1})\Big)$$ $$=U\Big(\mathscr{F}(J_{h_3})\circ \big(\mathscr{F}(J_{h_2})\circ \mathscr{F}(J_{h_1})\big)\Big)\xrightarrow{\epsilon^{-1}} U\Big(\mathscr{F}(J_{h_3})\circ \mathscr{F}U\big(\mathscr{F}(J_{h_2})\circ \mathscr{F}(J_{h_1})\big)\Big).$$ Tracing through the definitions above, in particular equations \eqref{eq:counit} and \eqref{eq:explicitconstraints}, and by writing $k_i:= g^{-1}h_i^{-1}g$ for $i=1,2,3$, we find that 

$$\alpha_{J_{h_1},J_{h_2},J_{h_3}}=\big(\epsilon_{\mathscr{F}(J_{h_3})\circ \mathscr{F}(J_{h_2})}\big)_{h_1gk_1}\cdot \big(\epsilon^{-1}_{\mathscr{F}(J_{h_2})\circ \mathscr{F}(J_{h_1})}\big)_g$$


$$=\frac{\psi(h_1,h_2,h_3)\pi(h_3,gk_3, k_2,k_1)\pi(g,k_3,k_2,k_1)}{\omega(k_3,k_2,k_1)\pi(h_2,h_3,gk_3k_2,k_1)\pi(h_1,h_2,h_3, gk_3k_2k_1)}.$$

\noindent This is precisely the expression of the 3-cocycle $\xi_g$ for $H\cap gKg^{-1}$ given in the statement of the proposition. The proof is therefore complete.
\end{proof}

\begin{Remark}
The following homotopical interpretation of proposition 
\ref{prop:2categorycosetdecomposition} was suggested by one of the referees. We have decided to include it not only because it is conceptually interesting, but also because it may be useful in establishing higher categorical versions of the proposition above.

Let $G$ be a finite group, and $H$, $K$ two subgroups of $G$. We may consider the corresponding diagram of spaces $\mathrm{B}H\rightarrow \mathrm{B}G\leftarrow\mathrm{B}K$, and write $P$ for its homotopy pullback. The space $P$ admits a concrete model given by the 1-groupoid whose objects are elements in $G$ and whose 1-morphisms $g_1\rightarrow g_2$ are pairs $(h,K)\in H\times K$ such that $hg_1 = g_2k$. It follows that the connected components of $P$ are in bijective correspondence with the double cosets $H\backslash G/K$, and, given such a double coset $[g]$, the corresponding group of automorphisms is $H\cap gKg^{-1}$. On one hand, observe that the collection of 3-cocycles $\xi_g$ from \eqref{eq:xig} may then be interpreted as a map of spaces $\xi: P\rightarrow \mathrm{B}^3\mathds{k}^{\times}$. On the other hand, it is possible to construct such a map of spaces from purely homotopical considerations. In detail, the 4-cocycle $\pi$ may be interpreted as a map of spaces $\mathrm{B}G\rightarrow \mathrm{B}^4\mathds{k}^{\times}$, and the 3-cocycles $\psi$ and $\omega$ as homotopies in the diagram of spaces depcited below: $$\begin{tikzcd}[sep = small]
\mathrm{B}H \arrow[rr] \arrow[dd] &  & \mathrm{B}G \arrow[dd, "\pi"] \arrow[lldd, Rightarrow, "\psi"', shorten >=3ex, shorten <=3ex] \arrow[rrdd, Rightarrow, "\omega", shorten >=3ex, shorten <=3ex] &  & \mathrm{B}K \arrow[dd] \arrow[ll] \\
                                  &  &                                                                            &  &                                   \\
* \arrow[rr]                      &  & \mathrm{B}^4\mathds{k}^{\times}                                            &  & *. \arrow[ll]                     
\end{tikzcd}$$ Upon taking homotopy pullbacks along both rows, this yields a map of spaces $\zeta:P\rightarrow \mathrm{B}^3\mathds{k}^{\times}$. We expect that $\xi$ and $\zeta$ are homotopy equivalent maps. We leave the details of this verification to the interested reader.
\end{Remark}

\begin{Corollary}
For $G$ a finite group, $H$ a subgroup, $\pi$ a 4-cocycle on $G$ and $\psi$ a trivialization on $H$, the underlying finite semisimple 2-category of $\mathfrak{C}(G,H,\pi,\psi) $ is $\boxplus_{[g]\in H \backslash G/H} \Mod(\Vect^{\xi_g}_{H \cap gHg^{-1}})$, where $\xi_g$ is the 3-cocycle for $H\cap gHg^{-1}$ defined in \eqref{eq:xig}.
\end{Corollary}

For later use, we also record the following related result.

\begin{Lemma}\label{lem:connectedcomponentunitgrouptheoretical}
For any finite groups $H\subseteq G$, 4-cocycle $\pi$ for $G$, and trivialization $\psi$ for $\pi|_H$, we have an equivalence $$\mathfrak{C}(G,H,\pi,\psi)^0\simeq \mathbf{2Rep}(H)$$ of fusion 2-categories.
\end{Lemma}
\begin{proof}
By definition, we have $\mathfrak{C}(G,H,\pi,\psi)=\mathbf{Bimod}_{\mathbf{2Vect}_G^{\pi}}(\mathbf{Vect}_H^{\psi})$ with monoidal unit given by $\mathbf{Vect}_H^{\psi}$. We wish to identify $End_{\mathbf{Vect}_H^{\psi}\mathrm{-}\mathbf{Vect}_H^{\psi}}(\mathbf{Vect}_H^{\psi})$ as a braided fusion 1-category. Given that this computation only involves bimodule morphisms, it is enough to carry it out in the special case $G=H$. As $d\psi=\pi|_H$, we may therefore without loss of generality assume that $\pi$ is trivial. Thus, we wish to describe $\mathfrak{C}(H,H,triv,\psi)$ the Morita dual to the fusion 2-category $\mathbf{2Vect}_H$ with respect to the finite semisimple module 2-category $\mathbf{2Vect}$ with $\mathbf{2Vect}_H$-module structure corresponding to $\psi$. This was done in proposition 4.3.3.10 of \cite{D:thesis}, and we have $$\mathfrak{C}(H,H,triv,\psi)\simeq \mathbf{2Rep}(H).$$ This concludes the proof of the result.
\end{proof}

\begin{Remark}\label{rem:alternative3cocycledescription}
We note that if $\pi$ is the trivial 4-cocycle for $G$, then for any $g\in G$, the 3-cocycle $\xi_g$ for $H\cap gKg^{-1}$ is cohomologous to the 3-cocycle $\zeta_g$ given on $h_1,h_2,h_3\in H\cap gKg^{-1}$ by $$\zeta_g(h_1,h_2,h_3):=\frac{\psi(h_1,h_2,h_3)}{\omega(g^{-1}h_1g,g^{-1}h_2g,g^{-1}h_3g)}.$$
\end{Remark}

Proposition \ref{prop:2categorycosetdecomposition} can be used to recover known descriptions of the underlying finite semisimple 2-category of certain fusion 2-categories of interest.

\begin{Example}\label{ex:2category2rep2group}
Let $\mathcal{G} = A[1]\rtimes H[0]$ be a finite split 2-group with Postnikov data $\beta \in H^3(\mathrm{B}\mathcal{G};\mathds{k}^{\times})$. We will establish in corollary below \ref{cor:grouptheoretical2group} that $\mathfrak{C}(H\ltimes \widehat{A}, \widehat{A}, \pi, triv)\simeq \mathbf{2Rep}(\mathcal{G})$. Noting that $H\backslash (\widehat{A}\rtimes H) / H\cong (\widehat{A})_H$, the set of orbits in $\widehat{A}$ under the action of $H$, and writing $\mathrm{Stab}_H(\alpha)$ for the subgroup of $H$ fixing $\alpha\in \widehat{A}$, the above proposition gives $$\mathbf{2Rep}(\mathcal{G})\simeq \boxplus_{\lbrack\alpha\rbrack\in(A)_H}\mathbf{Mod}(\mathbf{Vect}^{\alpha\circ\beta}_{\mathrm{Stab}_H(\alpha)}).$$ This result first appeared in \cite{El}.
\end{Example}

\begin{Example}
Let us fix $G$ a finite group. We now explain how to recover the description of the underlying 2-category of $\mathscr{Z}(\mathbf{2Vect}_G)$ obtained in \cite{KTZ}. Write $\Delta\hookrightarrow G\times G^{op}$ for the diagonal subgroup, we have $\Delta\backslash (G\times G^{op})/\Delta \cong \mathrm{Cl}(G)$, the set of conjugacy classes in $G$. Thus, combining the above proposition with example \ref{ex:center} yields $$\mathscr{Z}(\mathbf{2Vect}_G)\simeq \boxplus_{\lbrack g\rbrack\in\mathrm{Cl}(G)}\mathbf{Mod}(\mathbf{Vect}_{\mathrm{Z}(g)}),$$ where $\mathrm{Z}(g)$ denotes the centralizer of $g\in G$. Given $\pi\in Z^4_{gp}(G;\mathds{k}^{\times})$, a similar description was obtained for $\mathscr{Z}(\mathbf{2Vect}_G^{\pi})$ in \cite{KTZ}. Up to matching conventions and coboundaries, this has to agree with the description that is obtained by applying proposition \ref{prop:2categorycosetdecomposition} to the general case of example \ref{ex:center}.
\end{Example}

\subsection{A Partial Description of the Fusion Rules}\label{sub:fusionrules}

We now give a partial description of the fusion rules of the group-theoretical fusion 2-category $\mathfrak{C}(G,H,\pi,\omega)$. More precisely, given a class $\lbrack g\rbrack \in H\backslash G/ H$, we write $\mathbf{V}_{\lbrack g\rbrack}$ for the simple object of $\mathfrak{C}(G,H,\pi,\omega)$ given by $\mathbf{Vect}_{H\cap gHg^{-1}}^{\xi_g}$ in the connected component $\mathbf{Mod}(\mathbf{Vect}_{H\cap gHg^{-1}}^{\xi_g})$ as described in proposition \ref{prop:2categorycosetdecomposition}. We emphasize that $\mathbf{V}_{\lbrack e\rbrack}$ is in general not equivalent to the monoidal unit $I$ of $\mathfrak{C}(G,H,\pi,\omega)$, though they are related by a 2-condensation.

\begin{Proposition}\label{prop:partialfusionrule}
For any $f,g\in G$, we have $$\mathbf{V}_{\lbrack f\rbrack}\Box \mathbf{V}_{\lbrack g\rbrack}\simeq \boxplus_{k\in fHg}\mathbf{V}_{\lbrack k\rbrack}$$ in $\mathfrak{C}(G,H,\pi,\omega)$.
\end{Proposition}
\begin{proof}
Recall from equation (\ref{eq:definitionVg}) that $\mathbf{V}_{\lbrack g\rbrack} = \mathbf{Vect}_H^{\omega}\boxtimes \mathbf{Vect}_g\boxtimes \mathbf{Vect}_H^{\omega}$ is a simple object in $\mathfrak{C}(G,H,\pi,\omega)=\mathbf{Bimod}_{\mathbf{2Vect}_G^{\pi}}(\mathbf{Vect}_H^{\omega})$. Furthermore, the monoidal structure on the fusion 2-category $\mathbf{Bimod}_{\mathbf{2Vect}_G^{\pi}}(\mathbf{Vect}_H^{\omega})$ is given by $\boxtimes_{\mathbf{Vect}_H^{\omega}}$ the relative tensor product over $\mathbf{Vect}_H^{\omega}$ as defined in section 3.1 of \cite{D8}. Now, it is clear that $\mathbf{Vect}_H^{\omega}\boxtimes_{\mathbf{Vect}_H^{\omega}}\mathbf{Vect}_H^{\omega} = \mathbf{Vect}_H^{\omega}$ in $\mathbf{Bimod}_{\mathbf{2Vect}_G^{\pi}}(\mathbf{Vect}_H^{\omega})$. Thus, for any $f,g\in G$, we find that \begin{align*}\big(\mathbf{Vect}_H^{\omega}\boxtimes \mathbf{Vect}_f &\boxtimes \mathbf{Vect}_H^{\omega}\big)\boxtimes_{\mathbf{Vect}_H^{\omega}}\big(\mathbf{Vect}_H^{\omega}\boxtimes \mathbf{Vect}_g\boxtimes \mathbf{Vect}_H^{\omega}\big)\simeq \\ & \mathbf{Vect}_H^{\omega}\boxtimes \big(\mathbf{Vect}_f\boxtimes\mathbf{Vect}_H\boxtimes \mathbf{Vect}_g\big)\boxtimes \mathbf{Vect}_H^{\omega}\simeq \\ & \mathbf{Vect}_H^{\omega}\boxtimes \mathbf{Vect}_{fHg}\boxtimes \mathbf{Vect}_H^{\omega}.\end{align*} This gives the desired result.
\end{proof}

\begin{Remark}
For groups $H\subseteq G$ of low order, the above proposition can be used to determine the fusion rules completely. However, this is in general not enough. In fact, it would already be very interesting to derive the fusion rules of $\mathbf{2Rep}(\mathcal{G})$ for a general finite 2-group $\mathcal{G}$.
\end{Remark}

\subsection{Classification}

We investigate when two group-theoretical fusion 2-categories are monoidally equivalent. In order to do so, we first establish a result of independent interest, which characterizes the Morita equivalences of bosonic strongly fusion 2-categories.

\begin{Proposition}\label{prop:Morita2vectgroupequivalence}
Let $G_1$ and $G_2$ be two finite groups, and $\pi_1\in Z^4_{gp}(G_1;\mathds{k}^{\times})$, $\pi_2\in Z^4_{gp}(G_2;\mathds{k}^{\times})$. Then, $\mathbf{2Vect}_{G_1}^{\pi_1}$ and $\mathbf{2Vect}_{G_2}^{\pi_2}$ are Morita equivalent if and only if they are monoidally equivalent, that is there exists a group isomorphism $f:G_1\cong G_2$ such that $f^*\pi_2 = \pi_1$ in $H^4_{gp}(G_1;\mathds{k}^{\times})$.
\end{Proposition}
\begin{proof}
By definition, monoidally equivalent fusion 2-categories are Morita equivalent, so it is enough to prove the converse. Assume that $\mathbf{2Vect}_{G_1}^{\pi_1}$ and $\mathbf{2Vect}_{G_2}^{\pi_2}$ are Morita equivalent. In particular, there exists a connected separable algebra $\mathcal{A}$ in $\mathbf{2Vect}_{G_1}^{\pi_1}$ and a monoidal equivalence $\mathbf{Bimod}_{\mathbf{2Vect}_{G_1}^{\pi_1}}(\mathcal{A})\simeq \mathbf{2Vect}_{G_2}^{\pi_2}$ of fusion 2-categories. Recall from the proof of theorem 4.2.1 of \cite{D9} that $\mathcal{A}\Box \mathcal{A}$ is a simple object in the connected component of the identity of $\mathbf{Bimod}_{\mathbf{2Vect}_{G_1}^{\pi_1}}(\mathcal{A})$. But, we have $\mathbf{Bimod}_{\mathbf{2Vect}_{G_1}^{\pi_1}}(\mathcal{A})^0\simeq\mathbf{2Vect}$ by assumption, so that $\mathcal{A}\Box \mathcal{A}\simeq \mathcal{A}$ as $\mathcal{A}$-$\mathcal{A}$-bimodules. This implies that $\mathcal{A}\cong I$, as desired.
\end{proof}

The next result follows by inspecting the definitions. We recall that we follow the definitions of \cite{GPS} (see also section 2.3 of \cite{SP}) for monoidal 2-categories and monoidal 2-functors, and those of section 1.7 of \cite{EGNO} for group cohomology.

\begin{Proposition}\label{prop:monoidalequivalences}
Let $G_1$ and $G_2$ be two finite groups, and $\pi_1\in Z^4_{gp}(G_1;\mathds{k}^{\times})$, $\pi_2\in Z^4_{gp}(G_2;\mathds{k}^{\times})$. Then, monoidal equivalences up to monoidal 2-natural equivalence between $\mathbf{2Vect}_{G_1}^{\pi_1}$ and $\mathbf{2Vect}_{G_2}^{\pi_2}$ are classified by a group isomorphism $f:G_1\cong G_2$ together with a 3-cochain $\xi\in C^3_{gp}(G_1;\mathds{k}^{\times})$ considered up to 3-coboundary such that $\pi_1 = f^*{\pi_2}\cdot d\xi$.
\end{Proposition}
\begin{proof}
Observe that monoidal autoequivalences between $\mathbf{2Vect}_{G_1}^{\pi_1}$ and $\mathbf{2Vect}_{G_2}^{\pi_2}$ are uniquely determined by their restriction to the monoidal sub-2-categories $(\mathbf{2Vect}_{G_1}^{\pi_1})^{\times}$ and $(\mathbf{2Vect}_{G_2}^{\pi_2})^{\times}$ of invertible objects, invertible 1-morphisms, and invertible 2-morphisms. The result then follows from the homotopy hypothesis. Alternatively, and more concretely, let $\mathbf{F}:(\mathbf{2Vect}_{G_1}^{\pi_1})^{\times}\simeq (\mathbf{2Vect}_{G_2}^{\pi_2})^{\times}$ be a monoidal equivalence. It induces an isomorphism $f:G_1\cong G_2$ between the groups of invertible objects. By inspecting the definition of a monoidal 2-functor given in definition 3.1 of \cite{GPS}, we find that $\mathbf{F}$ also comes with the coherence data $\xi\in C^3_{gp}(G_1;\mathds{k}^{\times})$, and the equation $\pi_1 = f^*{\pi_2}\cdot d\xi$ corresponds to equation (HTA1) of \cite{GPS}. Conversely, the data of $f$ and $\xi$ as above yields a monoidal equivalence $(\mathbf{2Vect}_{G_1}^{\pi_1})^{\times}\simeq (\mathbf{2Vect}_{G_2}^{\pi_2})^{\times}$ and therefore also $\mathbf{2Vect}_{G_1}^{\pi_1}\simeq\mathbf{2Vect}_{G_2}^{\pi_2}$.
\end{proof}

The next result is essentially a categorification of \cite[Example 2.1 and Theorem 3.1]{O2} (see also \cite[Theorem 1.1]{Nat}), which classified Morita equivalence classes of algebras in pointed fusion 1-categories. We will use the following notation:\ For any $g\in G$ of a finite group $G$, we use $\mathrm{Ad}_g$ to denote the conjugation action by $g$ on $G$, that is $\mathrm{Ad}_g(x) = gxg^{-1}$ for every $x\in G$.

\begin{Lemma}\label{lem:Moritaequivalentalgebras}
Let $G$ be a finite group, and $\pi\in Z_{gp}^4(G;\mathds{k}^{\times})$. Let $H_1$ and $H_2$ be two subgroups of $G$, and $\psi_1\in C^3_{gp}(H_1;\mathds{k}^{\times})$, $\psi_2\in C^3_{gp}(H_2;\mathds{k}^{\times})$ such that $d\psi_1 = \pi|_{H_1}$ and $d\psi_2 = \pi|_{H_2}$. The connected rigid algebras $\mathbf{Vect}_{H_1}^{\psi_1}$ and $\mathbf{Vect}_{H_2}^{\psi_2}$ in $\mathbf{2Vect}_G^{\pi}$ are Morita equivalent, that is, have equivalent module 2-categories of modules, if and only if there exists $g\in G$ such that $\mathrm{Ad}_g(H_1) = H_2$ and $(\mathrm{Ad}_{g}^*\psi_2)/(\psi_1\cdot\Omega_g|_{H_1})$ is a 3-coboundary for $H_1$, where \begin{equation}\label{eq:Omegag}\Omega_g(g_1,g_2,g_3):=\frac{\pi(\mathrm{Ad}_g(g_1),\mathrm{Ad}_g(g_2),g,g_3)\pi(g,g_1,g_2,g_3)}{\pi(\mathrm{Ad}_g(g_1),\mathrm{Ad}_g(g_2),\mathrm{Ad}_g(g_3),g)\pi(\mathrm{Ad}_g(g_1),g,g_2,g_3)}\end{equation} for every $g_1,g_2,g_3\in G$.
\end{Lemma}
\begin{proof}
We will consider the $\mathbf{2Vect}_G^{\pi}$ module 2-category $\mathbf{Mod}_{\mathbf{2Vect}_G^{\pi}}(\mathbf{Vect}_{H_1}^{\psi_1})$. By theorem 4.2.2 of \cite{D4}, this endows the 2-category $\mathbf{Mod}_{\mathbf{2Vect}_G^{\pi}}(\mathbf{Vect}_{H_1}^{\psi_1})$ with a $\mathbf{2Vect}_G^{\pi}$-enriched structure. For any object $M$ of $\mathbf{Mod}_{\mathbf{2Vect}_G^{\pi}}(\mathbf{Vect}_{H_1}^{\psi_1})$, we write $\underline{End}(M)$ for the corresponding algebra in $\mathbf{2Vect}_G^{\pi}$. We note that this algebra is rigid by theorem 5.2.7 of \cite{D4}, and that it is connected if and only if $M$ is simple. Now, it follows from theorem 5.1.2 of \cite{D8} that every connected rigid algebra in $\mathbf{2Vect}_G^{\pi}$ that is Morita equivalent to $\mathbf{Vect}_{H_1}^{\psi_1}$ must be of the form $\underline{End}(M)$  for some simple object $M$ of $\mathbf{Mod}_{\mathbf{2Vect}_G^{\pi}}(\mathbf{Vect}_{H_1}^{\psi_1})$. It will therefore suffice to identify these algebras.

The simple objects of $\mathbf{Mod}_{\mathbf{2Vect}_G^{\pi}}(\mathbf{Vect}_{H_1}^{\psi_1})$ are all of the form $\mathbf{Vect}_g\boxtimes\mathbf{Vect}_{H_1}^{\psi_1}$ for some $g\in G$. Let us also record that $$\underline{End}(\mathbf{Vect}_{H_1}^{\psi_1})\simeq \mathbf{Vect}_{H_1}^{\psi_1}$$ as algebras in $\mathbf{Mod}_{\mathbf{2Vect}_G^{\pi}}(\mathbf{Vect}_{H_1}^{\psi_1})$. Now, fix an arbitrary element $g\in G$. The 2-functor $\mathbf{L}_g$ of left multiplication by the invertible object $\mathbf{Vect}_g$ of $\mathbf{2Vect}_G^{\pi}$ induces an equivalence of finite semisimple 2-categories $$\mathbf{L}_g := \mathbf{Vect}_g\boxtimes(-):\mathbf{Mod}_{\mathbf{2Vect}_G^{\pi}}(\mathbf{Vect}_{H_1}^{\psi_1})\xrightarrow{\sim}\mathbf{Mod}_{\mathbf{2Vect}_G^{\pi}}(\mathbf{Vect}_{H_1}^{\psi_1}).$$ This equivalence can be upgraded to an equivalence of left $\mathbf{2Vect}_G^{\pi}$-module 2-categories. Namely, let us write $\mathbf{F}_g$ for the monoidal autoequivalence of $\mathbf{2Vect}_G^{\pi}$ given by conjugation by the invertible object $\mathbf{Vect}_{g}$, that is $$\mathbf{F}_g:=\big(\mathbf{Vect}_g\boxtimes (-)\big)\boxtimes \mathbf{Vect}_{g^{-1}}:\mathbf{2Vect}_G^{\pi}\xrightarrow{\sim}\mathbf{2Vect}_G^{\pi}.$$
It follows by direct inspection that, under the bijection of proposition \ref{prop:Morita2vectgroupequivalence}, $\mathbf{F}_g$ corresponds to the group isomorphism $\mathrm{Ad}_g$ and the 3-cochain $\Omega_g\in C^3_{gp}(G;\mathds{k}^{\times})$ defined above. The key point being the relation $d\Omega_g = \pi/\mathrm{Ad}_g^*\pi$. It follows by construction that $\mathbf{L}_g$ promotes to an equivalence  of left $\mathbf{2Vect}_G^{\pi}$-module 2-categories $$\mathbf{L}_g:\mathbf{Mod}_{\mathbf{2Vect}_G^{\pi}}(\mathbf{Vect}_{H_1}^{\psi_1})\xrightarrow{\sim}{_{\mathbf{F}_g}\mathbf{Mod}_{\mathbf{2Vect}_G^{\pi}}(\mathbf{Vect}_{H_1}^{\psi_1})}.$$ We emphasize that the left action on the right hand-side is twisted by $\mathbf{F}_g$ in the manner discussed at the end of section \S\ref{sub:MoritaEquivalence}.

Finally, as the enrichment is defined using right adjoints \cite{D4}, it follows from their uniqueness that $$\underline{End}(\mathbf{Vect}_{g}\boxtimes\mathbf{Vect}_{H_1}^{\psi_1})\simeq \mathbf{F}_{g}\big(\underline{End}(\mathbf{Vect}_{H_1}^{\psi_1})\big)$$ as connected rigid algebras in $\mathbf{2Vect}_G^{\pi}$. Moreover, we have that $$\mathbf{F}_{g}\big(\underline{End}(\mathbf{Vect}_{H_1}^{\psi_1})\big)\simeq \mathbf{Vect}_{\mathrm{Ad}_g(H_1)}^{(\mathrm{Ad}_{g}^{-1})^*(\psi_1\cdot \Omega_g)}$$ as algebras by tracing through the definitions, so that the proof is complete.
\end{proof}

\begin{Proposition}\label{prop:classificationgrouptheoretical}
Two group-theoretical fusion 2-categories $\mathfrak{C}(G_1,H_1,\pi_1,\psi_1)$ and $\mathfrak{C}(G_2,H_2,\pi_2,\psi_2)$ are monoidally equivalent if and only if there exists an isomorphism $f:G_1\cong G_2$ and a class $\xi\in C^3_{gp}(G_1;\mathds{k}^{\times})$ such that 
\begin{enumerate}[label=(\alph*)]
    \item We have $\pi_1 = f^*\pi_2\cdot d\xi$ in $Z^4_{gp}(G_1;\mathds{k}^{\times})$;
    \item We have $\mathrm{Ad}_g(f(H_1)) = H_2$ for some $g\in G_2$;
    \item The class $\big(f^*(\mathrm{Ad}_g)^*\psi_2\big)/(\xi|_{H_1}\cdot \psi_1\cdot f^*\Omega_g|_{H_1})$ is a 3-coboundary for $H_1$ with $\Omega_g\in C^3(G_2;\mathds{k}^{\times})$ as defined in \eqref{eq:Omegag}.
\end{enumerate}
\end{Proposition}
\begin{proof}
Let $\mathbf{E}:\mathfrak{C}(G_1,H_1,\pi_1,\psi_1)\simeq\mathfrak{C}(G_2,H_2,\pi_2,\psi_2)$ be a monoidal equivalence of fusion 2-categories. Let us also consider the canonical invertible $\mathbf{2Vect}_{G_1}^{\pi_1}$-$\mathfrak{C}(G_2,H_2,\pi_2,\psi_2)$-bimodule 2-category $\mathbf{Mod}_{\mathbf{2Vect}_{G_1}^{\pi_1}}(\mathbf{Vect}_{H_1}^{\psi_1})$ and the invertible $\mathbf{2Vect}_{G_1}^{\pi_1}$-$\mathfrak{C}(G_2,H_2,\pi_2,\psi_2)$-bimodule 2-category $\mathbf{Mod}_{\mathbf{2Vect}_{G_2}^{\pi_2}}(\mathbf{Vect}_{H_2}^{\psi_2})$. These two Morita equivalences supply the horizontal arrows in the diagram below. Furthermore, thanks to the transitivity of Morita equivalence, which follows from theorem 5.4.3 of \cite{D8} and is recalled above in \S\ref{sub:MoritaEquivalence}, there exists a Morita equivalence, represented by the dashed arrow, that makes the square below commute $$\begin{tikzcd}
{\mathfrak{C}(G_1,H_1,\pi_1,\psi_1)} \arrow[dd, "\mathbf{E}"'] &  & \mathbf{2Vect}_{G_1}^{\pi_1} \arrow[dd, "\mathbf{F}"', dashed] \arrow[ll, leftrightarrow, "\sim"'] \\
&  & \\
{\mathfrak{C}(G_2,H_2,\pi_2,\psi_2)}                       &  & \mathbf{2Vect}_{G_2}^{\pi_2}. \arrow[ll,leftrightarrow, "\sim"]\arrow[uu, "\mathbf{F}^{-1}"', dashed, bend right]                     
\end{tikzcd}$$ It follows from proposition \ref{prop:Morita2vectgroupequivalence} above that the dashed Morita equivalence is actually implemented by a monoidal equivalence $\mathbf{F}:\mathbf{2Vect}_{G_1}^{\pi_1}\simeq\mathbf{2Vect}_{G_2}^{\pi_2}$. Under proposition \ref{prop:Morita2vectgroupequivalence}, $\mathbf{F}$ corresponds to an isomorphism $f:G_1\cong G_2$ and a 3-cochain $\xi\in C^3(G_1;\mathds{k}^{\times})$ such that $\pi_1 = f^*\pi_2\cdot d\xi$. Its pseudo-inverse $\mathbf{F}^{-1}$ then corresponds to the isomorphism $f^{-1}$ with 3-cochain $(f^{-1})^*(\xi^{-1})\in C^3(G_2;\mathds{k}^{\times})$. 

Commutativity of the above square implies that there is an equivalence of $\mathbf{2Vect}_{G_2}^{\pi_2}$-$\mathfrak{C}(G_2,H_2,\pi_2,\psi_2)$-bimodule 2-category between $\mathbf{Mod}_{\mathbf{2Vect}_{G_2}^{\pi_2}}(\mathbf{Vect}_{H_2}^{\psi_2})$ and $_{\mathbf{F}^{-1}}\mathbf{Mod}_{\mathbf{2Vect}_{G_1}^{\pi_1}}(\mathbf{Vect}_{H_1}^{\psi_1})_{\mathbf{E}}$. The subscripts $\textbf{F}^{-1}$ and $\textbf{E}$ are used to indicate that the actions are twisted. This notation was introduced at the end of section \S\ref{sub:MoritaEquivalence}. But, with $\varpi:= (f^{-1})^*(\psi_1\cdot \xi|_{H_1})$, we have $\mathbf{F}^{-1}(\mathbf{Vect}_{f^{-1}(H_1)}^{\varpi})\simeq \mathbf{Vect}_{H_1}^{\psi_1}$ by inspection, so that $_{\mathbf{F}^{-1}}\mathbf{Mod}_{\mathbf{2Vect}_{G_1}^{\pi_1}}(\mathbf{Vect}_{H_1}^{\psi_1})\simeq \mathbf{Mod}_{\mathbf{2Vect}_{G_2}^{\pi_2}}(\mathbf{Vect}_{f^{-1}(H_2)}^{\varpi})$ as left $\mathbf{2Vect}_{G_2}^{\pi_2}$-module 2-categories. The forward direction of the result therefore follows from lemma \ref{lem:Moritaequivalentalgebras} above.

Conversely, it follows from proposition \ref{prop:monoidalequivalences} that the data of $f:G_1\cong G_2$ and $\xi\in C^3_{gp}(G_1;\mathds{k}^{\times})$ with $\pi_1 = f^*\pi_2\cdot d\xi$ provides us with a monoidal equivalence $\mathbf{F}:\mathbf{2Vect}_{G_1}^{\pi_1}\rightarrow \mathbf{2Vect}_{G_2}^{\pi_2}$. The 3-cochain $\xi$ prescribes the coherence data for $\mathbf{F}$, we therefore find that $\mathbf{F}(\mathbf{Vect}_{H_1}^{\psi_1})\simeq \mathbf{Vect}_{f(H_1)}^{\varpi}$ as algebras in $\mathbf{2Vect}_{G_2}^{\pi_2}$ with $\varpi=(f^{-1})^*\big(\psi_1\cdot\xi|_{H_1}\big)$ as in the preceeding paragraph. But, if $\big(f^*(\mathrm{Ad}_g)^*\psi_2\big)/(\xi|_{H_1}\cdot \psi_1\cdot f^*\Omega_g|_{H_1})$ is a 3-coboundary for $H_1$, we find that $\mathbf{Vect}_{f(H_1)}^{\varpi}$ and $\mathbf{Vect}_{H_2}^{\psi_2}$ are Morita equivalent connected rigid algebras in $\mathbf{2Vect}_{G_2}^{\pi_2}$ thanks to lemma \ref{lem:Moritaequivalentalgebras}. Putting this discussion together, we have monoidal equivalences $$\mathfrak{C}(G_1,H_1,\pi_1,\psi_1) \simeq\mathbf{Bimod}_{\mathbf{2Vect}_{G_2}^{\pi_2}}(\mathbf{Vect}_{f(H_1)}^{\varpi})
\simeq\mathfrak{C}(G_2,H_2,\pi_2,\psi_2).$$ This gives the backward direction of the result, and thereby concludes the proof.
\end{proof}

\begin{Example}
Let $G=\mathbb{Z}/4\oplus\mathbb{Z}/2$, and $H=\mathbb{Z}/2\oplus\mathbb{Z}/2$. Using the K\"unneth formula, we find that the canonical map $H^3_{gp}(G;\mathds{k}^{\times})\rightarrow H^3_{gp}(H;\mathds{k}^{\times})$ is isomorphic to $P=pr\oplus 0\oplus Id:\mathbb{Z}/4\oplus\mathbb{Z}/2\oplus\mathbb{Z}/2\rightarrow\mathbb{Z}/2\oplus\mathbb{Z}/2\oplus\mathbb{Z}/2$. This computation follows from the naturality of the K\"unneth formula. Let us write $\psi$ for a 3-cocycle for $H$ representing the non-zero element in the cokernel, which is isomorphic to $\mathbb{Z}/2$. Then, $\mathfrak{C}(G,H,triv,triv)$ and $\mathfrak{C}(G,H,triv,\psi)$ are inequivalent fusion 2-categories.
\end{Example}

\begin{Example}\label{ex:cocyclesgivesdistinctf2cs}
Let $A$ be a finite abelian group together with an action by $\mathbb{Z}/2$, i.e.\ the data of an automorphism $\rho:A\cong A$ with $\rho^2 = Id$. In addition, let $\psi$ be a 3-cocycle for $A$ such that $\omega:=\psi/\rho^*\psi$ is a non-trivializable 3-cocycle. It follows from proposition \ref{prop:2categorycosetdecomposition} and remark \ref{rem:alternative3cocycledescription} that the underlying finite semisimple 2-categories of $\mathfrak{C}(A\rtimes\mathbb{Z}/2,A,triv,triv)$ and $\mathfrak{C}(A\rtimes\mathbb{Z}/2,A,triv,\psi)$ are distinct. In fact, they do not even have the same number of simple objects. Namely, it follows from proposition \ref{prop:2categorycosetdecomposition} that there are equivalences of 2-categories $$\mathfrak{C}(A\rtimes\mathbb{Z}/2,A,triv,triv)\simeq \mathbf{Mod}(\mathbf{Vect}_A)\boxplus\mathbf{Mod}(\mathbf{Vect}_A),$$ $$\mathfrak{C}(A\rtimes\mathbb{Z}/2,A,triv,\psi)\simeq \mathbf{Mod}(\mathbf{Vect}_A)\boxplus\mathbf{Mod}(\mathbf{Vect}_A^{\omega}),$$ and it follows from the main theorem of \cite{Nat} that $\mathbf{Mod}(\mathbf{Vect}_A^{\omega})$ has fewer equivalence classes of simple objects than $\mathbf{Mod}(\mathbf{Vect}_A)$. As a specific example, one can take $A=\mathbb{Z}/2\oplus\mathbb{Z}/2$, $\rho$ the automorphism given by swapping the two summand, and $\psi$ a 3-cocycle representing the pullback of the non-trivial class in $H^3_{gp}(\mathbb{Z}/2;\mathds{k}^{\times})\cong \mathbb{Z}/2$ along the projection $\mathbb{Z}/2\oplus\mathbb{Z}/2\twoheadrightarrow\mathbb{Z}/2$ to the first factor.
\end{Example}

\begin{Example}\label{ex:exotic2groupequivalence}
In section 3.4 of \cite{BSNW} (see also example 3.5.2 of \cite{BBFP}), the authors observed that there should be an equivalence of fusion 2-categories $$\mathbf{2Rep}(\mathbb{Z}/4\lbrack 1\rbrack \rtimes \mathbb{Z}/2\lbrack 0\rbrack )\simeq \mathbf{2Rep}((\mathbb{Z}/2\oplus \mathbb{Z}/2)\lbrack 1\rbrack \rtimes \mathbb{Z}/2\lbrack 0\rbrack ),$$ which might a priori look surprising. However, using the classification of proposition \ref{prop:classificationgrouptheoretical}, these two fusion 2-categories correspond to the group-theoretical fusion 2-categories $\mathfrak{C}(D_8, \langle s\rangle, triv, triv)$ and $\mathfrak{C}((\mathbb{Z}/2\oplus \mathbb{Z}/2)\rtimes \mathbb{Z}/2, 0\rtimes \mathbb{Z}/2, triv, triv)$, with $D_8$ the dihedral group of order 8, and $\langle s\rangle$ a reflection subgroup. It is well-known that $(\mathbb{Z}/2\oplus \mathbb{Z}/2)\rtimes \mathbb{Z}/2\cong D_8$, and it is not hard to show that such an isomorphism can be chosen to preserve the two $\mathbb{Z}/2$ subgroups under consideration. Thus, the above equivalence of fusion 2-categories is recovered by the proposition.
\end{Example}

\section{Fusion 2-Categories with a Fiber 2-Functor}\label{section:Fib2Functor}

\subsection{Classification}
\begin{Definition}
Let $\mathfrak{C}$ be a fusion 2-category. A fiber 2-functor for $\mathfrak{C}$ is a monoidal 2-functor $\mathfrak{C}\rightarrow \mathbf{2Vect}$.
\end{Definition}

The following definition from classical group theory \cite{Za, Sz} will play a key role in our classification result.

\begin{Definition}
Let $G$ be a finite group. An exact factorization of $G$ consists of two subgroups $H$ and $K$ of $G$ such that $H\cap K = \{e\}$ and $HK = G$.\footnote{This type of factorization of $G$ is also referred to as a Zappa-Sz\'ep product of $H$ and $K$.} If $H$ is a fixed subgroup, we say that a subgroup $K$ as above is a complement for $H$ in $G$.
\end{Definition}

The following result is a categorification of corollary 3.1 of \cite{O2}. Physically speaking, it completely describes what kind of fusion 2-category the topological operators in a $(2+1)$-dimensional theory must form in order for all the (noninvertible) symmetries to be gauged.

\begin{Theorem}\label{thm:fiber2functor}
A fusion 2-category admits a fiber 2-functor to $\mathbf{2Vect}$ if and only if it is equivalent to a group-theoretical fusion 2-category $\mathfrak{C}(G,H,\pi,\psi)$ such that $H$ admits a complement $K$ with $\pi|_K$ trivializable.
\end{Theorem}
\begin{proof}
Let $\mathfrak{C}$ be a fusion 2-category admitting a fiber 2-functor $\mathbf{F}:\mathfrak{C}\rightarrow\mathbf{2Vect}$. In particular, $\mathbf{2Vect}$ may be viewed as a finite semisimple left $\mathfrak{C}$-module 2-category, so that $\mathfrak{C}$ is group-theoretical by corollary \ref{cor:grouptheoreticaltechnical}. By definition, it therefore follows that $\mathfrak{C}$ is a fusion 2-category of the form $\mathfrak{C}(G,H,\pi,\psi)$ for some finite group $G$ equipped with a subgroup $H\subseteq G$ together with a 4-cocycle $\pi$ on $G$ and a trivialization $\psi$ for $\pi|_H$.

Now, the data of the 2-functor $\mathbf{F}:\mathfrak{C}\rightarrow\mathbf{2Vect}$ is equivalent to the data of a left $\mathfrak{C}(G,H,\pi,\psi)$-module structure on $\mathbf{2Vect}$. By the transitivity of Morita equivalence, there exists a connected rigid algebra $\mathcal{A}$ in $\mathbf{2Vect}_G^{\pi}$ such that $$\mathbf{2Vect}\simeq \mathbf{Bimod}_{\mathbf{2Vect}_G^{\pi}}(\mathbf{Vect}_H^{\psi}, \mathcal{A})$$ as left $\mathfrak{C}(G,H,\pi,\psi)=\mathbf{Bimod}_{\mathbf{2Vect}_G^{\pi}}(\mathbf{Vect}_H^{\psi})$-module 2-categories. But, the dual fusion 2-category is then given by $\big(\mathbf{Bimod}_{\mathbf{2Vect}_G^{\pi}}(\mathcal{A})\big)^{mop}$, which is necessarily group-theoretical by corollary \ref{cor:grouptheoreticaltechnical}. In particular, as we have assumed that $\mathcal{A}$ is connected, we must have that $\mathcal{A} \simeq \mathbf{Vect}_K^{\omega}$ for some subgroup $K\subseteq G$ and 3-cochain $\omega$ for $K$ such that $d\omega = \pi|_K$ by lemma \ref{lem:rigidalgebragrouptheoretical}. It then follows from proposition \ref{prop:2categorycosetdecomposition} that $$\mathbf{2Vect}\simeq \mathbf{Bimod}_{\mathbf{2Vect}_G^{\pi}}(\mathbf{Vect}_H^{\psi}, \mathbf{Vect}_K^{\omega})$$ as 2-categories if and only if both $HK= G$ and $H\cap K = \{e\}$, i.e.\ $H$ and $K$ give an exact factorization of $G$. This last argument not only finishes the proof of the forward direction, but also establishes the backward direction. The proof of the result is therefore complete.
\end{proof}

With respect to $\Vect^\psi_H$, the algebra $\mathcal{A}$ used in the above proof can be viewed as the ``screened objects". If one decides instead to condense the algebra $\mathcal{A}\simeq\Vect^\omega_K$, then $\Vect^{\pi}_H$ would give the screened objects.

\begin{Corollary}\label{cor:fiber2functor}
Let $\mathfrak{C}(G,H,\pi,\psi)$ be a group-theoretical fusion 2-category. Fiber 2-functors for $\mathfrak{C}(G,H,\pi,\psi)$ up to monoidal 2-natural equivalence are classified by pairs $(K,\omega)$ such that $K$ is a complement for $H$ in $G$ and $\omega$ is a 3-cochain for $K$ such that $\pi|_K=d\omega$, where two such pairs $(K_1,\omega_1)$ and $(K_2,\omega_2)$ are equivalent if and only if there exists an element $g\in G$ such that $\mathrm{Ad}_g(K_1) = K_2$, and $\mathrm{Ad}_g^*\omega_2/(\omega_1\cdot \Omega_g|_{K_1})$ is a 3-coboundary for $K_1$.
\end{Corollary}
\begin{proof}
By construction, we have that $\mathfrak{C}(G,H,\pi,\psi)$ and $\mathbf{2Vect}_G^{\pi}$ are Morita equivalent fusion 2-categories via $\mathbf{Mod}_{\mathbf{2Vect}_G^{\pi}}(\mathbf{Vect}_H^{\psi})$. In particular, it follows from theorem 5.4.3 of \cite{D8} that the two finite semisimple $\mathfrak{C}(G,H,\pi,\psi)$-module 2-categories $$\mathbf{Bimod}_{\mathbf{2Vect}_G^{\pi}}(\mathbf{Vect}_H^{\psi},\mathbf{Vect}_{K_1}^{\omega_1})\ \mathrm{and}\ \mathbf{Bimod}_{\mathbf{2Vect}_G^{\pi}}(\mathbf{Vect}_H^{\psi},\mathbf{Vect}_{K_2}^{\omega_2})$$ are equivalent if and only if the two finite semisimple $\mathbf{2Vect}_G^{\pi}$-module 2-categories $$\mathbf{Mod}_{\mathbf{2Vect}_G^{\pi}}(\mathbf{Vect}_{K_1}^{\omega_1})\ \mathrm{and}\ \mathbf{Mod}_{\mathbf{2Vect}_G^{\pi}}(\mathbf{Vect}_{K_2}^{\omega_2})$$ are equivalent. The statement therefore follows from lemma \ref{lem:Moritaequivalentalgebras}.
\end{proof}

The next corollary follows from the definitions recalled in section \S\ref{sub:MoritaEquivalence}.
\begin{Corollary}\label{cor:F2Fdual}
The dual to $\mathfrak{C}(G,H,\pi,\psi)$ with respect to the fiber 2-functor given by the complement $K\subseteq G$, and $\omega\in C_{gp}^3(K;\mathds{k}^{\times})$, is $\mathfrak{C}(G,K,\pi,\omega)$.
\end{Corollary}

In the perspective of figure \ref{fig:bimodAB}, the fiber 2-functor that we have constructed corresponds to $\Bimod_{\textbf{2}\Vect^\pi_G}(\Vect^\psi_H,\Vect^\omega_K)$, which gives the gapped boundary between the two dual phases. 

\subsection{Applications}

We now apply our classification of fiber 2-functors to offer complete and general mathematical proofs of some results in  \cite{BBFP,BBSNT,DelT} regarding gauging in fusion 2-categories associated to 2-groups. Let $\mathcal{G}=A[1]\boldsymbol{\cdot}H[0]$ be a finite 2-group. In particular, there is an action $\rho:H\rightarrow \mathrm{Aut}(A)$, and a 3-cocycle $\beta\in Z^3_{gp}(H;A)$. Without loss of generality, we can assume that $\beta$ is normalized. The group $H$ acts on $\widehat{A}$, the Pontryagin dual of $A$, by precomposition, so that we can consider the finite group $H\ltimes \widehat{A}$.\footnote{The assignment $(h,\alpha)\mapsto (\alpha \circ \rho(h^{-1}), h)$ defines an isomorphism of groups $H\ltimes \widehat{A}\cong\widehat{A}\rtimes H$. We use $H\ltimes \widehat{A}$ because it yields a very simple description of $\pi$.} We define a 4-cocycle $\pi$ for the finite group $H\ltimes \widehat{A}$ with coefficients in $\mathds{k}^{\times}$ by \begin{equation}\pi((h_1,\alpha_1), (h_2,\alpha_2),(h_3,\alpha_3),(h_4,\alpha_4)):=\alpha_1(\beta(h_2,h_3,h_4)),\label{eq:definitionpi}\end{equation} for every $h_1,...,h_4\in H$, and $\alpha_1,...,\alpha_4\in \widehat{A}$. The 4-cocycle condition for $\pi$ follows in a straightforward fashion from the 3-cocycle condition for $\beta$. We also record that $\pi$ is not only normalized, but is furthermore trivial if any of $\alpha_1$, $h_2$, $h_3$, or $h_4$ is trivial.

Physically speaking, this 4-cocycle has the expected form: We may regard $\pi$ as the mixed anomaly in the gauged theory and consider the Serre spectral sequence used to compute such an anomaly. The anomaly is a class in $H^4_{gp}(H\ltimes \widehat{A}; \mathds{k}^\times)$, and the $E_2$ page in total degree 4 has two nonvanishing pieces given by $H^1_{gp}(H;H^3_{gp}(\widehat{A};\mathds{k}^\times))$ and $H^3_{gp}(H,H^1_{gp}(
\widehat{A};\mathds{k}^\times))$. But, using $H^1_{gp}(
\widehat{A};\mathds{k}^\times) = \widehat{(\widehat{A})} = A$, we find that the data of $\beta$ gives a class in $H^3_{gp}(H;A) = H^3_{gp}(H; H^1_{gp}(
\widehat{A};\mathds{k}^\times))$. As a consequence, we see that the contribution 
to the anomaly naturally lives in $H^3_{gp}(H,H^1_{gp}(\widehat{A};\mathds{k}^\times))$ on the $E_2$ page and so the formula for $\pi$ should reflect this fact. Since $\beta$ is valued in $A$, and $\alpha_1 \in \widehat{A}$, the formula in equation \eqref{eq:definitionpi} is valued in $\mathds{k}^\times$.

\begin{Proposition}\label{prop:2groupgraded2vectorspaces}
For any finite 2-group $\mathcal{G} = A[1]\boldsymbol{\cdot}H[0]$ with Postnikov class $\beta$,
we have an equivalence of monoidal 2-categories $$\mathfrak{C}(H\ltimes \widehat{A}, \widehat{A}, \pi, triv)\simeq \mathbf{2Vect}_{\mathcal{G}}\,,$$
where $\pi$ is the 4-cocycle on $H\ltimes \widehat{A}$ defined in \eqref{eq:definitionpi} whose restriction to $\widehat{A}$ is trivial.
\end{Proposition}
\begin{proof}
We wish to describe the 2-category of $\mathbf{Vect}_{\widehat{A}}$-$\mathbf{Vect}_{\widehat{A}}$-bimodules in the fusion 2-category $\mathbf{2Vect}_{H\ltimes \widehat{A}}^{\pi}$. We will identify this fusion 2-category with the fusion 2-category of 2-vector spaces graded by the finite 2-group $\mathcal{G}$. In remark 2.1.17 of \cite{DR}, an effective characterization of such fusion 2-categories is given. More precisely, in order to prove that a fusion 2-category $\mathfrak{C}$ is equivalent to the fusion 2-category of (potentially twisted) 2-vector spaces graded by a 2-group, it is enough to check that every simple 1-morphism in $\Omega\mathfrak{C}$ is invertible, and that there exists a finite set of invertible objects in $\mathfrak{C}$, exactly one in each connected component of $\mathfrak{C}$, which is closed under the monoidal product.

We begin by introducing an auxiliary class of simple bimodules, which is distinct from those given in (\ref{eq:definitionVg}). For any $h\in H$, we consider the object $\mathbf{Vect}_{h\widehat{A}}$ of $\mathbf{2Vect}_{H\ltimes \widehat{A}}^{\pi}$ with its canonical $\mathbf{Vect}_{\widehat{A}}$-$\mathbf{Vect}_{\widehat{A}}$-bimodule structure. In the notation of definition \ref{def:bimodule}, the structure 2-isomorphisms $\kappa$, $\nu$, and $\beta$ for $\mathbf{Vect}_{h\widehat{A}}$ are all trivial. This follows from the specific form of the 4-cocycle $\pi$ given in (\ref{eq:definitionpi}). In fact, for almost all of the computations that we will carry out below, the 4-cocycle $\pi$ will not appear.

We claim that \begin{equation}\label{eq:1moprhisms}End_{\mathbf{Vect}_{\widehat{A}}\mathrm{-}\mathbf{Vect}_{\widehat{A}}}(\mathbf{Vect}_{h\widehat{A}})\simeq \mathbf{Vect}_{A}.\end{equation} Namely, it is clear that $End_{\mathbf{Vect}_{\widehat{A}}}(\mathbf{Vect}_{h\widehat{A}})\simeq \mathbf{Vect}$, so that there is a unique simple right $\mathbf{Vect}_{\widehat{A}}$-module 1-endomorphism $F_h:\mathbf{Vect}_{h\widehat{A}}\rightarrow \mathbf{Vect}_{h\widehat{A}}$ in $\mathbf{2Vect}_{H\ltimes \widehat{A}}^{\pi}$. The corresponding structure 2-isomorphism $\phi_h$ can be safely assumed to be trivial. On the other hand, by inspecting definition \ref{def:bimodulemorphism}, we find that endowing the right $\mathbf{Vect}_{\widehat{A}}$-module 1-morphism $F_h$ with a compatible left $\mathbf{Vect}_{\widehat{A}}$-module structure is equivalent to the data of an element in $\widehat{\big(\widehat{A}\big)}$.

We give explicit representatives for the simple $\mathbf{Vect}_{\widehat{A}}$-bimodule 1-morphisms $\mathbf{Vect}_{h\widehat{A}}\rightarrow \mathbf{Vect}_{h\widehat{A}}$ in $\mathbf{2Vect}_{H\ltimes \widehat{A}}^{\pi}$ of the previous paragraph. Given $a\in\widehat{(\widehat{A})}=A$, we write $${^{a}F_h}:\mathbf{Vect}_{h\widehat{A}}\rightarrow \mathbf{Vect}_{h\widehat{A}}$$ for the simple $\mathbf{Vect}_{\widehat{A}}$-$\mathbf{Vect}_{\widehat{A}}$-bimodule 1-morphism whose right $\mathbf{Vect}_{\widehat{A}}$-action $^a\phi_h$ is the trivial one, and whose left $\mathbf{Vect}_{\widehat{A}}$-action $^a\chi_h$ is twisted by evaluating at $a$. Explicitly, we have $$(^a\chi_h)_{\alpha_1,(\alpha_2,h)}=\alpha_1(a):\mathbf{Vect}_{\alpha_1}\otimes {^aF_h}(\mathbf{Vect}_{(\alpha_2,h)})\cong {^aF_h}(\mathbf{Vect}_{\alpha_1}\otimes \mathbf{Vect}_{(\alpha_2,h)})$$ for every $\alpha_1,\alpha_2\in\widehat{A}$ and $h\in H$, and $(^a\phi_h)_{(\alpha_1,h),\alpha_2}=1$. It is also useful to consider the simple $\mathbf{Vect}_{\widehat{A}}$-bimodule 1-morphisms $F^a_h:\mathbf{Vect}_{h\widehat{A}}\rightarrow \mathbf{Vect}_{h\widehat{A}}$ whose left $\mathbf{Vect}_{\widehat{A}}$-action $\chi^a_h$ is the trivial one, and whose right $\mathbf{Vect}_{\widehat{A}}$-action $\phi^a_h$ is twisted by evaluating at $a$. It is straightforward to check using the definitions of \S\ref{ap:detaileddefinitions} that $$F^a_h\simeq {^{\rho(a)}F_h}:\mathbf{Vect}_{h\widehat{A}}\rightarrow \mathbf{Vect}_{h\widehat{A}}$$ as $\mathbf{Vect}_{\widehat{A}}$-bimodule 1-morphisms in $\mathbf{2Vect}_{H\ltimes \widehat{A}}^{\pi}$

As our next step, we analyze in which connected components of $\mathfrak{C}(H\ltimes \widehat{A}, \widehat{A}, \pi, triv)$ the bimodules $\mathbf{Vect}_{h\widehat{A}}$ lie.
It follows from the above paragraph that the $\mathbf{Vect}_{\widehat{A}}$-bimodule 1-morphism $F^0_h$, which is the unit bimodule 1-morphism on $\mathbf{Vect}_{h\widehat{A}}$, is a simple bimodule 1-morphism for every $h\in H$. By definition, this shows that the bimodule $\mathbf{Vect}_{h\widehat{A}}$ is simple. Furthermore, there is exactly one bimodule of this form in every connected component of $\mathfrak{C}(H\ltimes \widehat{A}, A, \pi, triv)$. This is a consequence of the observation that there is no non-zero 1-morphism in $\mathbf{2Vect}_{H\ltimes\widehat{A}}^{\pi}$ between $\mathbf{Vect}_{h_1\widehat{A}}$ and $\mathbf{Vect}_{h_2\widehat{A}}$ if $h_1\neq h_2$.

We now compute the fusion rules of the bimodules $\mathbf{Vect}_{h\widehat{A}}$. Given any $h_1,h_2\in H$, we find that \begin{equation}\label{eq:fusionrules}\mathbf{Vect}_{h_1\widehat{A}}\boxtimes_{\mathbf{Vect}_{\widehat{A}}}\mathbf{Vect}_{h_2\widehat{A}}\simeq \mathbf{Vect}_{h_1h_2\widehat{A}}\end{equation} Namely, the canonical $\mathbf{Vect}_{\widehat{A}}$-bimodule 1-morphism $\mathbf{Vect}_{h_1\widehat{A}}\boxtimes\mathbf{Vect}_{h_2\widehat{A}}\simeq \mathbf{Vect}_{h_1h_2\widehat{A}}$ satisfies the desired 2-universal property. Again, note that, thanks to the specific form of $\pi$ given in (\ref{eq:definitionpi}), the 4-cocycle $\pi$ does not appear. As a consequence, we find that $\mathbf{Vect}_{h\widehat{A}}$ is an invertible $\mathbf{Vect}_{\widehat{A}}$-$\mathbf{Vect}_{\widehat{A}}$-bimodule for any $h\in H$. Similarly, we have \begin{equation}\label{eq:action1}^{a}F_{h_1}\boxtimes_{\mathbf{Vect}_{\widehat{A}}}\mathbf{Vect}_{h_2\widehat{A}}\simeq {^{a}F_{h_1h_2}}.\end{equation} Likewise, we have $\mathbf{Vect}_{h_1\widehat{A}}\boxtimes_{\mathbf{Vect}_{\widehat{A}}}F^{a}_{h_2}\simeq F^{a}_{h_1h_2}$ as $\mathbf{Vect}_{\widehat{A}}$-bimodule 1-morphisms, which readily implies that \begin{equation}\label{eq:action2}\mathbf{Vect}_{h_1\widehat{A}}\boxtimes_{\mathbf{Vect}_{\widehat{A}}}{^{a}F_{h_2}}\simeq {^{\rho(h_1)(a)}F_{h_1h_2}}.\end{equation}

Furthermore, we claim that the associator 1-morphism
\begin{align*}(\mathbf{Vect}_{h_1{\widehat{A}}} \boxtimes_{\mathbf{Vect}_{\widehat{A}}}\mathbf{Vect}_{h_2{\widehat{A}}})&\boxtimes_{\mathbf{Vect}_{\widehat{A}}}\mathbf{Vect}_{h_3{\widehat{A}}}\\ & \simeq\mathbf{Vect}_{h_1{\widehat{A}}}\boxtimes_{\mathbf{Vect}_{\widehat{A}}}(\mathbf{Vect}_{h_2{\widehat{A}}}\boxtimes_{\mathbf{Vect}_{\widehat{A}}}\mathbf{Vect}_{h_3{\widehat{A}}})\end{align*} is given by $F^{\beta(h_1,h_2,h_3)}_{h_1h_2h_3}$ for any $h_1,h_2,h_3\in H$. Namely, it follows from the construction of the monoidal structure on 2-categories of bimodules given in section 3.2 of \cite{D8} that the above associator is obtained via a 2-universal property from the canonical equivalence $$(\mathbf{Vect}_{h_1{\widehat{A}}} \boxtimes\mathbf{Vect}_{h_2{\widehat{A}}})\boxtimes\mathbf{Vect}_{h_3{\widehat{A}}}\simeq\mathbf{Vect}_{h_1{\widehat{A}}}\boxtimes(\mathbf{Vect}_{h_2{\widehat{A}}}\boxtimes\mathbf{Vect}_{h_3{\widehat{A}}})$$ of $\mathbf{Vect}_{\widehat{A}}$-$\mathbf{Vect}_{\widehat{A}}$-bimodules. The underlying 1-morphism in $\mathbf{2Vect}_{H\ltimes \widehat{A}}^{\pi}$ is the identity. Given the specific form of $\pi$, the right $\mathbf{Vect}_{\widehat{A}}$-module structure is the canonical one. In contrast, the presence of $\pi$ forces the left $\mathbf{Vect}_{\widehat{A}}$-module to be twisted by evaluating at $\beta(h_1,h_2,h_3)$. Since the $\mathbf{Vect}_{\widehat{A}}$-bimodule structure on the associator 1-morphism is induced from those we have just described, this finishes the proof of the claim.

To summarize the above discussion, we have seen in equation \eqref{eq:fusionrules} that the $\mathbf{Vect}_{\widehat{A}}$-bimodules $\mathbf{Vect}_{h\widehat{A}}$ are invertible and their tensor products are described by the group structure of $H$. In addition, we have found in \eqref{eq:1moprhisms} that $\mathbf{Vect}_A$ is the fusion 1-category of bimodule 1-endomorphisms on $\mathbf{Vect}_{h\widehat{A}}$. We have therefore identified a finite 2-group of the form $A[1]\boldsymbol{\cdot}H[0]$ as a monoidal sub-1-category of the monoidal homotopy 1-category $\mathrm{h}\mathfrak{C}(H\ltimes \widehat{A}, \widehat{A}, \pi, triv)$. But, it follows from equations \eqref{eq:action1} and \eqref{eq:action2} together with the preceding paragraph that $A[1]\boldsymbol{\cdot}H[0]\simeq \mathcal{G}$. We can therefore appeal to the criterion given in remark 2.1.17 of \cite{DR} to conclude that $\mathfrak{C}(H\ltimes \widehat{A}, \widehat{A}, \pi, triv)\simeq \mathbf{2Vect}_{\mathcal{G}}^{\varpi}$ for some 4-cocycle $\varpi$. Finally, using the explicit description of $\pi$ given in (\ref{eq:definitionpi}), theorem \ref{thm:fiber2functor} implies that there exists a fiber 2-functor $\mathfrak{C}(H\ltimes \widehat{A}, \widehat{A}, \pi, triv)\rightarrow \mathbf{2Vect}$. This proves that $\varpi$ must be a coboundary, and concludes the proof of the theorem.
\end{proof}

\begin{Remark}\label{rem:partialgaugingofextensions}
Let $H\boldsymbol{\cdot} A$ be an extension of $H$ by $A$ classified by an action of $H$ on $A$ and a class in $Z^2_{gp}(H;A_{\rho})$. The first part of the proof of the above argument continues to apply, and shows that $\mathfrak{C}(H\boldsymbol{\cdot} A, A, triv, triv)\simeq \mathbf{2Vect}_{\widehat{A}[1]\rtimes H[0]}^{\pi}$ for some $\pi\in Z^4_{gp}(\widehat{A}[1]\rtimes H[0];\mathds{k}^{\times})$. The 4-cocycle $\pi$ is an elusive quantity, and is difficult to determine in general. We do identify it for the extension $\mathbb{Z}/2\boldsymbol{\cdot}\mathbb{Z}/2 = \mathbb{Z}/4$ in example \ref{ex:gauginZ/2inZ/4} below.
\end{Remark}

The following corollary has already appeared in the Physics literature in \cite{BBFP} (see also \cite{BBSNT, DelT} for the case of split 2-groups).

\begin{Corollary}\label{cor:grouptheoretical2group}
For any finite 2-group $\mathcal{G} = A[1]\boldsymbol{\cdot}H[0]$ with Postnikov class $\beta$,
we have an equivalence of monoidal 2-categories $$\mathfrak{C}(H\ltimes \widehat{A}, H, \pi, triv)\simeq \mathbf{2Rep}(\mathcal{G})\,,$$
where $\pi$ is the 4-cocycle on $H\ltimes \widehat{A}$ defined in \eqref{eq:definitionpi} whose restriction to $H$ is trivial.
\end{Corollary}
\begin{proof}
It was shown in proposition 4.3.3.10 of \cite{D:thesis} that the Morita dual of $\mathbf{2Vect}_{\mathcal{G}}$ with respect to the module 2-category $\mathbf{2Vect}$ is $\mathbf{2Rep}(\mathcal{G})$. As observed in remark 4.3.3.11 therein, this identification does not depend on the choice of $\mathbf{2Vect}_{\mathcal{G}}$-module structure on $\mathbf{2Vect}$. Alternatively, this also follows from the fact that all of these module structures are related by monoidal autoequivalences of $\mathbf{2Vect}_{\mathcal{G}}$. The result therefore follows by combining proposition \ref{prop:2groupgraded2vectorspaces} with corollary \ref{cor:F2Fdual}.
\end{proof}

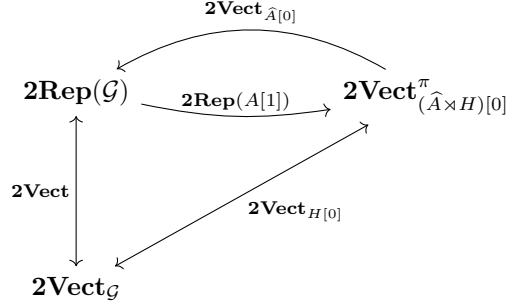
\begin{figure}[htb]
$$\begin{tikzcd}
\mathbf{2Rep}(\mathcal{G}) \arrow[rrr, "{\mathbf{2Rep}(A[1])}", bend right=10] &  &  & {\mathbf{2Vect}^{\pi}_{(\widehat{A}\rtimes H)[0]}} \arrow[lllddd, leftrightarrow, "{\mathbf{2Vect}_{H[0]}}"] \arrow[lll, "{\mathbf{2Vect}_{\widehat{A}[0]}}"', bend right] \\
 &  &  &  \\
 &  &  &  \\
\mathbf{2Vect}_{\mathcal{G}} \arrow[uuu, "\mathbf{2Vect}", leftrightarrow]                 &  &  &                                     
\end{tikzcd}$$
\caption{Morita equivalences between fusion 2-categories constructed from an arbitrary finite 2-group $\mathcal{G}={A}[1]\boldsymbol{\cdot}H[0]$. Recall that there is an equivalence $\mathbf{2Rep}(A[1])\simeq \mathbf{2Vect}_{\widehat{A}[0]}$ of finite semisimple 2-categories.}
    \label{fig:triangleRelation}
\end{figure}

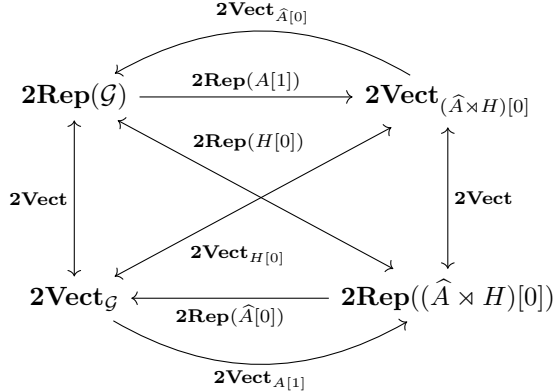
\begin{figure}[hbt]
$$\begin{tikzcd}
\mathbf{2Rep}(\mathcal{G}) \arrow[rrr, "{\mathbf{2Rep}(A[1])}"] \arrow[rrrddd, "{\mathbf{2Rep}(H[0])}", near start, leftrightarrow] &  &  & {\mathbf{2Vect}_{(\widehat{A}\rtimes H)[0]}} \arrow[lllddd, "{\mathbf{2Vect}_{H[0]}}", near end, leftrightarrow] \arrow[lll, "{\mathbf{2Vect}_{\widehat{A}[0]}}"', bend right] \\
&  &  &  \\
&  &  & \\
\mathbf{2Vect}_{\mathcal{G}} \arrow[uuu, "\mathbf{2Vect}", leftrightarrow] \arrow[rrr, "{\mathbf{2Vect}_{A[1]}}"', bend right]      &  &  & {\mathbf{2Rep}((\widehat{A}\rtimes H)[0])} \arrow[uuu, "\mathbf{2Vect}"', leftrightarrow] \arrow[lll, "{\mathbf{2Rep}(\widehat{A}[0])}"]                     
\end{tikzcd}$$
\caption{Morita equivalences between fusion 2-categories constructed from a split finite 2-group $\mathcal{G}={A}[1]\rtimes H[0]$. The arrow going to the bottom left of the square represents gauging $\mathbf{2Vect}_{H[0]}$, and the arrow going to the bottom right of the square represents gauging $\mathbf{2Rep}(H[0])$.
Recall that there is an equivalence $\mathbf{2Rep}(\widehat{A}[0])\simeq \mathbf{2Vect}_{A[1]}$ of finite semisimple 2-categories.}
\label{fig:squareRelation}
\end{figure}

\begin{Remark}
Let $\mathcal{G}= {A}[1]\boldsymbol{\cdot}H[0]$ be a finite 2-group. Our results yield the diagram of Morita equivalences depicted in figure \ref{fig:triangleRelation}. Namely, proposition \ref{prop:2groupgraded2vectorspaces} gives the top equivalence of the triangle, and corollary \ref{cor:grouptheoretical2group} gives the right side of the triangle. This diagram first appeared in the physics literature in section 3.5 of \cite{BBFP} (see also \cite{BBSNT} for particular cases).

When $\mathcal{G}$ is a split 2-group, the diagram of figure \ref{fig:triangleRelation} can be extended to the one in figure \ref{fig:squareRelation}. The missing diagonal was obtained in example 5.1.9 of \cite{D8}. The other missing equivalences are obtained by appealing to corollary \ref{cor:F2Fdual} above. This square has first appeared at a physical level of rigor in \cite{BBFP:I} (see also \cite{BBSNT} for particular cases). A more careful derivation was also given in \cite{DelT}.
\end{Remark}

\subsection{Examples}


We give some noteworthy examples of group-theoretical fusion 2-categories that illustrate the existence and behaviour of fiber 2-functors.

\begin{Example}
Take $G=S_4$, $H=S_3\subseteq S_4$, and the trivial 4-cocycle. Then, the fusion 2-category $\mathfrak{C}(S_4,S_3,triv,triv)$ admits a fiber 2-functor. However, the complement $K$ to $S_3$ in $S_4$ is not uniquely determined! For instance, one could take $K$ to be $\mathbb{Z}/4$, but also $\mathbb{Z}/2\oplus \mathbb{Z}/2$. Thence, the two associated fiber 2-functors are distinct. Moreover, different choices of trivialization for $\pi|_K$ yield inequivalent fiber 2-functors. These are parameterized by $H^3_{gp}(K;\mathds{k}^{\times})$.

This example is also particularly interesting for the following reason: The exact factorization of $S_4$ given by $H$ and $K$ is not a semi-direct product decomposition. Namely, neither $H$ nor $K$ is a normal subgroup of $S_4$. In particular, $\mathfrak{C}(S_4,S_3,triv,triv)$ is a fusion 2-category admitting a fiber 2-functor to $\mathbf{2Vect}$ that is not equivalent to the fusion 2-category of 2-representations of a 2-group, nor the fusion 2-category of twisted 2-group graded vector spaces. The picture below represents the underlying finite semisimple 2-category of $\mathfrak{C}(S_4,S_3,triv,triv)$. The simple objects are described using the main result of \cite{Nat}.
$$\begin{tikzcd}
& I \arrow["\mathbf{Rep}(\mathbb{Z}/2)"', ld, bend right, leftrightarrow] \arrow["\mathbf{Vect}" description, near start, dd, leftrightarrow] \arrow["\mathbf{Rep}(S_3)"', loop, distance=2em, in=125, out=55] &\\
X \arrow[rd, bend right, leftrightarrow] \arrow[rr, leftrightarrow] \arrow[loop, distance=2em, in=215, out=145] & & Y \arrow["\mathbf{Rep}(\mathbb{Z}/3)"',lu, bend right, leftrightarrow] \arrow[loop, distance=2em, in=35, out=325] \\
& Z \arrow[ru, bend right, leftrightarrow] \arrow["\mathbf{Vect}_{S_3}"',loop, distance=2em, in=305, out=235]           &\\ \\
A \arrow["\mathbf{Vect}_{\mathbb{Z}/2}"', loop, distance=2em, in=215, out=145] & & B \arrow["\mathbf{Vect}", ll, leftrightarrow] \arrow["\mathbf{Vect}_{\mathbb{Z}/2}"', loop, distance=2em, in=35, out=325]            
\end{tikzcd}$$
\end{Example}

\begin{Example}
Take $G=S_3$ and $H=\mathbb{Z}/3$. Then, thanks to corollary \ref{cor:grouptheoretical2group}, we have that $\mathfrak{C}(S_3, \mathbb{Z}/3, triv,triv)\simeq \mathbf{2Vect}_{\mathbb{Z}/3\lbrack 1\rbrack\rtimes\mathbb{Z}/2\lbrack 0\rbrack}$. On the other hand, if $H=\mathbb{Z}/2$, then, according to corollary \ref{cor:grouptheoretical2group}, we find $\mathfrak{C}(S_3, \mathbb{Z}/2, triv,triv)\simeq \mathbf{2Rep}({\mathbb{Z}/3\lbrack 1\rbrack\rtimes\mathbb{Z}/2\lbrack 0\rbrack})$.  In fact, as the subgroups $\mathbb{Z}/2,\ \mathbb{Z}/3\subseteq S_3$ form a Zappa-Sz\'ep decomposition of $S_3$, then not only do these two fusion 2-categories admit fiber 2-functors, but they are Morita dual to one another. This example was first considered in section 4 of \cite{BBSNT}.
\end{Example}

\begin{Example}\label{ex:gauginZ/2inZ/4}
We now study in detail an example of a group-theoretical fusion 2-category that does not admit a fiber 2-functor.
More precisely, we will consider the fusion 2-category $\mathfrak{C}(\mathbb{Z}/4,\mathbb{Z}/2,triv,triv)$, which was also studied in \cite{BBSNT}. Thanks to corollary \ref{cor:fiber2functor}, we see that $\mathfrak{C}(\mathbb{Z}/4,\mathbb{Z}/2,triv,triv)$ has no fiber 2-functor, as $\mathbb{Z}/2\subseteq \mathbb{Z}/4$ has no complement. On the other hand, it follows from remark \ref{rem:partialgaugingofextensions} that we have $\mathfrak{C}(\mathbb{Z}/4,\mathbb{Z}/2,triv,triv)\simeq \mathbf{2Vect}^{\varpi}_{\mathbb{Z}/2\lbrack 1\rbrack\times\mathbb{Z}/2\lbrack 0\rbrack}$ for some $\varpi$ in $Z^4(\mathbb{Z}/2\lbrack 2\rbrack\times\mathbb{Z}/2\lbrack 1\rbrack;\mathds{k}^{\times})$.

Thanks to the K\"unneth formula, we have that $H^4(\mathbb{Z}/2\lbrack 2\rbrack\times\mathbb{Z}/2\lbrack 1\rbrack;\mathds{k}^{\times})\cong \mathbb{Z}/4\oplus \mathbb{Z}/2$. More precisely, the first summand corresponds to $H^4(\mathbb{Z}/2\lbrack 2\rbrack;\mathds{k}^{\times})\cong \mathbb{Z}/4$, but $\varpi$ can not have a non-zero summand in $H^4(\mathbb{Z}/2\lbrack 2\rbrack;\mathds{k}^{\times})$ as we have $\Omega\mathfrak{C}(\mathbb{Z}/4,\mathbb{Z}/2,triv,triv)\simeq \mathbf{Rep}(\mathbb{Z}/2)$. Further, we see that $\varpi$ can not be trivial, as this would imply that $\mathfrak{C}(\mathbb{Z}/4,\mathbb{Z}/2,triv,triv)$ admits a fiber 2-functor. Thus, $\varpi$ must represent the non-trivial cohomology class in $H^4(\mathbb{Z}/2\lbrack 2\rbrack\times\mathbb{Z}/2\lbrack 1\rbrack;\mathds{k}^{\times})$ corresponding to the non-trivial element in the $\mathbb{Z}/2$ factor. In particular, $\varpi$ is non-trivial as was observed in \cite{BBSNT} (see figure 2 therein).
\end{Example}

\begin{Example}\label{ex:interestingfiber2functor2rep}
Let $\mathcal{G}:=\mathbb{Z}/2\lbrack 1\rbrack\boldsymbol{\cdot}\mathbb{Z}/2\lbrack 0\rbrack$ be the finite 2-group with non-trivial Postnikov data. We have $\mathbf{2Vect}_{\mathcal{G}}\simeq \mathfrak{C}(\mathbb{Z}/2\oplus\mathbb{Z}/2, \mathbb{Z}/2\oplus 0, \kappa, triv)$ with $\kappa$ a 4-cocycle for $\mathbb{Z}/2\oplus\mathbb{Z}/2$. The class of $\kappa$ in $H^4_{gp}(\mathbb{Z}/2\oplus\mathbb{Z}/2,\mathds{k}^{\times})\cong \mathbb{Z}/2\oplus\mathbb{Z}/2$ must be non-trivial because $\mathfrak{C}(\mathbb{Z}/2\oplus\mathbb{Z}/2, \mathbb{Z}/2\oplus 0, triv, triv)\simeq \mathbf{2Vect}_{\mathcal{H}}$ with $\mathcal{H}=\mathbb{Z}/2\lbrack 1\rbrack\times\mathbb{Z}/2\lbrack 0\rbrack$. Then, thanks to theorem \ref{thm:fiber2functor} and corollary \ref{cor:fiber2functor}, the fusion 2-category $\mathbf{2Vect}_{\mathcal{G}}$ has two inequivalent fiber 2-functors given by the diagonal subgroup $\Delta = \langle (1,1)\rangle$, and the subgroup $0\oplus \mathbb{Z}/2$ of $\mathbb{Z}/2\oplus\mathbb{Z}/2$. Nevertheless, these two fiber 2-functors may be identified upon precomposing with an autoequivalence. Explicitly, the group automorphism $f:\mathbb{Z}/2\oplus\mathbb{Z}/2\cong \mathbb{Z}/2\oplus\mathbb{Z}/2$ with $f(1,0)=(1,0)$ and $f(0,1)=(1,1)$ is such that $f(\mathbb{Z}/2\oplus 0)= \mathbb{Z}/2\oplus 0$
and $f^*\kappa = \kappa$ in $H^4_{gp}(\mathbb{Z}/2\oplus\mathbb{Z}/2,\mathds{k}^{\times})$, so that it induces a monoidal autoequivalence of $\mathfrak{C}(\mathbb{Z}/2\oplus\mathbb{Z}/2, \mathbb{Z}/2\oplus 0, \kappa, triv)$. Moreover, as $f(\Delta)=0\oplus\mathbb{Z}/2$, precomposition with this monoidal autoequivalence identifies the two fiber 2-functors of $\mathfrak{C}(\mathbb{Z}/2\oplus\mathbb{Z}/2, \mathbb{Z}/2\oplus 0, triv, triv)$.

Now, the dual to $\mathfrak{C}(\mathbb{Z}/2\oplus\mathbb{Z}/2, \mathbb{Z}/2\oplus 0, \kappa, triv)$ with respect to the finite semisimple module 2-category corresponding to $\Delta$ is $\mathfrak{C}(\mathbb{Z}/2\oplus\mathbb{Z}/2, \Delta, \kappa, triv)$. But, the dual fusion 2-category to $\mathbf{2Vect}_{\mathcal{G}}$ with respect to $\mathbf{2Vect}$ is $\mathbf{2Rep}(\mathcal{G})$, so we must have $\mathfrak{C}(\mathbb{Z}/2\oplus\mathbb{Z}/2, \Delta, \kappa, triv)\simeq \mathbf{2Rep}(\mathcal{G})$. Similarly, we find that $\mathfrak{C}(\mathbb{Z}/2\oplus\mathbb{Z}/2, 0\oplus \mathbb{Z}/2, \kappa, triv)\simeq \mathbf{2Rep}(\mathcal{G})$. Combining the above discussion together, we find that there exists a fiber 2-functor $\mathbf{2Rep}(\mathcal{G})\rightarrow \mathbf{2Vect}$ whose corresponding dual fusion 2-category is $\mathbf{2Rep}(\mathcal{G})$.
\end{Example}

\begin{Remark}
We have seen in the proof of corollary \ref{cor:grouptheoretical2group} that for any finite 2-group $\mathcal{G}$ and any fiber 2-functor $\mathbf{2Vect}_{\mathcal{G}}\rightarrow \mathbf{2Vect}$, the dual to $\mathbf{2Vect}_{\mathcal{G}}$ with respect to $\mathbf{2Vect}$ is $\mathbf{2Rep}(\mathcal{G})$. The last example above shows that there exists finite 2-groups $\mathcal{G}$ and fiber 2-functors $\mathbf{2Rep}(\mathcal{G})\rightarrow \mathbf{2Vect}$, such that the dual to $\mathbf{2Rep}(\mathcal{G})$ with respect to $\mathbf{2Vect}$ is not $\mathbf{2Vect}_{\mathcal{G}}$. It would be interesting to exhibit this duality in some physical $(2+1)$-dimensional system. This behaviour is well-known in the decategorified setting, and can be attributed to the fact that categories of projective representations for certain pairs of subgroup and 2-cocycle can yield additional fiber functors (see for instance corollary 7.12.20 of \cite{EGNO} and the subsequent remark).
\end{Remark}


\section{Fusion 2-Categories with a Tambara-Yamagami defect}\label{section:F2C_selfduality}
We define and study fusion 2-categories with a Tambara-Yamagami defect. More specifically, our objective is to categorify Tambara-Yamagami 1-categories, which possess a duality line. 
Going from Tambara-Yamagami 1-categories to the corresponding 2-categorical notion requires some care. Namely, it is too strong to require that the trivially graded factor be the fusion 2-category of 2-vector spaces graded by a finite abelian (1-)group, as this implies that the fusion 2-category is strongly fusion \cite{JFY}. Thus, we will instead require that the trivially graded factor be equivalent to the 2-category of 2-vector spaces graded by a 2-group. 

We start by reviewing the definition of Tambara-Yamagami 1-category. In \cite{TY}, Tambara and Yamagami introduced and studied fusion 1-categories associated to a finite group $A$, whose fusion rules are given by
\begin{equation}
  a \cdot b = ab\,, \quad a \cdot m = m\,, \quad m \cdot m = \bigoplus_{a\in A} a.
\end{equation}
Further, they show that $A$ is necessarily abelian. Reformulating the above properties, we obtain the following definition.

\begin{Definition}
A Tambara-Yamagami 1-category is a $\mathbb{Z}/2$-graded fusion 1-category $\mathcal{C}=\mathcal{C}_+\oplus \mathcal{C}_-$ such that $\mathcal{C}_+$ is pointed and $\mathcal{C}_-\simeq \mathbf{Vect}$.
\end{Definition}

\noindent Using direct explicit computations, Tambara and Yamgami obtained a complete classification of Tambara-Yamgami 1-categories.

We begin our analysis of the categorification of this notion by establishing some general properties of group-graded fusion 2-categories, and listing some noteworthy examples. Then, in \S\ref{sub:TYdefect}, we go on to study a general class of fusion 2-categories which have a simple object whose fusion rules are reminiscent of the duality line in Tambara-Yamagami 1-categories. After having established this general classification, we turn our attention in \S\ref{subsection:TYF2C} to Tambara-Yamagami 2-categories, that is, to the case when $\fC_{+} = \textbf{2}\Vect_{\cG}$, for some finite 2-group $\mathcal{G}$.

\subsection{Group-Graded Fusion 2-Categories}\label{sub:groupgradedF2C}

\begin{Definition}
Let $\mathfrak{C}$ be a fusion 2-category, and $G$ be a finite group. A $G$\textit{-grading} on $\mathfrak{C}$ is a decomposition $\mathfrak{C} = \boxplus_{g\in G}\mathfrak{C}_g$ of $\mathfrak{C}$ into a direct sum of finite semisimple full sub-2-categories such that given any $C$ in $\mathfrak{C}_g$ and $D$ in $\mathfrak{C}_h$, then $C\Box D$ is in $\mathfrak{C}_{gh}$. A $G$-grading on $\mathfrak{C}$ is \textit{faithful} if $\mathfrak{C}_g$ is non-zero for every $g$ in $G$.
\end{Definition}

\begin{Remark}
Let $\mathfrak{C}$ be a $G$-graded fusion 2-category. Then, as $\mathfrak{C}_g$ is a full sub-2-category, it is a union of connected components of $\mathfrak{C}$.
\end{Remark}

Just as for 1-categories, the defects for the symmetry $G$ define superselection sectors in the theory.

\begin{Definition}
Let $C$ be an object of the $G$-graded fusion 2-category $\mathfrak{C}$. The \textit{support} of $C$ is the subset of $G$ defined by $$\supp(C):=\{g\in G|C \textrm{ has a non-zero summand in } \mathfrak{C}_g\}.$$
\end{Definition}

It follows readily from the definitions and lemma 2.3.7 of \cite{D2} that the support is multiplicative in the sense that $\supp(C\Box D) = \supp(C)\supp(D)$ for any objects $C$ and $D$ in the $G$-graded fusion 2-category $\mathfrak{C}$. This property can be used to deduce the following lemma on how condensation can cut down on superselection sectors.

\begin{Lemma}\label{lemma:gradings}
Let $\mathfrak{C}$ be a faithfully $G$-graded fusion 2-category, and $\mathcal{A}$ an indecomposable separable algebra whose support is contained in $N$, a normal subgroup of $G$. Then, $\Bimod_{\mathfrak{C}}(\mathcal{A})$ is a faithfully $G/N$-graded fusion 2-category.
\end{Lemma}
\begin{proof}
Let $\mathcal{P}$ be any simple $\mathcal{A}$-$\mathcal{A}$-bimodule. Note that $\supp(\mathcal{P})$ is contained in $gN$ for some $g$ in $G$. Namely, for any simple summand $C$ of $\mathcal{P}$ in $\mathfrak{C}$, there is a non-zero 1-morphism $\mathcal{A}\Box C\Box \mathcal{A}\rightarrow \mathcal{P}$ of $\mathcal{A}$-$\mathcal{A}$-bimodules. In particular, as $\mathcal{P}$ is simple, it can be obtained by splitting a 2-condensation monad, as defined in \cite{GJF}, which is supported on $\mathcal{A}\Box C\Box \mathcal{A}$. But, we have $$\supp(\mathcal{A}\Box C\Box \mathcal{A})=\supp(\mathcal{A})\supp(C)\supp(\mathcal{A}) \subseteq N\supp(C)N= \supp(C) N,$$ which establishes the claim as the splitting of a 2-condensation monad induces an inclusion on supports.

It remains to show that this grading is compatible with the monoidal structure. Let $\mathcal{P}$ and $\mathcal{Q}$ be simple $\mathcal{A}$-$\mathcal{A}$-bimodules. Thanks to the argument above, there exists simple objects $C$ and $D$ in $\mathfrak{C}$ and non-zero 1-morphisms $\mathcal{A}\Box C\Box \mathcal{A}\rightarrow \mathcal{P}$ and $\mathcal{A}\Box D\Box \mathcal{A}\rightarrow \mathcal{Q}$. Recall from section 3.1 of \cite{D8} that $P\Box_{\mathcal{A}} Q$ is the splitting of a 2-condensation monad supported on $\mathcal{P}\Box \mathcal{Q}$, so that $$\supp(\mathcal{P}\Box_{\mathcal{A}} \mathcal{Q})\subseteq \supp(\mathcal{P})\supp(\mathcal{Q})\subseteq \supp(C)\supp(D) N.$$ This proves that $\Bimod_{\mathfrak{C}}(\mathcal{A})$ is a faithfully $G/N$-graded fusion 2-category.
\end{proof}

\begin{Definition}
Let $H\subseteq G$, we use $H^G$ to denote the minimal normal subgroup of $G$ containing $H$.
\end{Definition}
By definition, the fusion 2-category $\mathfrak{C}(G,H,\pi,\psi)$ is the fusion 2-category of bimodules over $\mathbf{Vect}_H^{\psi}$ in $\mathbf{2Vect}_G^{\pi}$. In particular, $\supp(\mathbf{Vect}_H^{\psi})= H$ is contained in  $H^G$. Lemma \ref{lemma:gradings} therefore yields the first part of the next result.

\begin{Corollary}\label{cor:grouptheoreticalgraded}
The fusion 2-category $\mathfrak{C}(G,H,\pi,\psi)$ is faithfully $G/H^G$-graded. Moreover, this grading is universal in the sense that any other faithful $\widetilde{G}$-grading on $\mathfrak{C}(G,H,\pi,\psi)$ by a finite group $\widetilde{G}$ is induced by a surjective group homomorphism $G/H^G\twoheadrightarrow \widetilde{G}$.
\end{Corollary}
\begin{proof}
Let $\widetilde{G}$ be a finite group, and let $\mathfrak{C}(G,H,\pi,\psi) = \boxplus_{\tilde{g}\in \widetilde{G}}\mathfrak{C}(G,H,\pi,\psi)^{\tilde{g}}$ be a faithful $\widetilde{G}$-grading. In the decomposition of proposition \ref{prop:2categorycosetdecomposition}, it follows from proposition \ref{prop:partialfusionrule} that the dual to the connected component indexed by $\lbrack g\rbrack\in H\backslash G/H$ is the connected component indexed by $\lbrack g^{-1}\rbrack$. Thence, for any $g\in G$, $\mathbf{V}_{\lbrack g\rbrack}\Box \mathbf{V}_{\lbrack g^{-1}\rbrack}$ is an object of $\mathfrak{C}(G,H,\pi,\psi)^{\tilde{e}}$. But, thanks to proposition \ref{prop:partialfusionrule}, for any $h\in H^G$, there exists $g\in G$ such that $\mathbf{V}_{\lbrack h\rbrack}$ is a summand of $\mathbf{V}_{\lbrack g\rbrack}\Box \mathbf{V}_{\lbrack g^{-1}\rbrack}$. This shows that the fusion sub-2-category $\mathfrak{C}(G,H,\pi,\psi)^{\tilde{e}}$ contains $\mathfrak{C}(G,H,\pi,\psi)_e$, the fusion sub-2-category associated to the canonical $G/H^G$-grading. But, using proposition \ref{prop:partialfusionrule} again, we have that for any $g\in G/H^G$, $\mathfrak{C}(G,H,\pi,\psi)_g$ is an indecomposable $\mathfrak{C}(G,H,\pi,\psi)_e$-module 2-category. Thus, we find that $\mathfrak{C}(G,H,\pi,\psi)_g$ is contained in $\mathfrak{C}(G,H,\pi,\psi)^{\tilde{g}}$ for a unique $\tilde{g}\in\widetilde{G}$ as desired.
\end{proof}

Let $\mathfrak{C}$ be a faithfully $G$-graded fusion 2-category. For any $g\in G$, the full sub-2-category $\mathfrak{C}_g$ inherits a $\mathfrak{C}_e$-$\mathfrak{C}_e$-bimodule structure. In particular, there is a canonical monoidal 2-functor $\mathfrak{C}_e^{mop}\rightarrow \mathbf{End}_{\mathfrak{C}_e}(\mathfrak{C}_g)$.

\begin{Proposition}\label{prop:gradedMoritaequivalence}
Let $\mathfrak{C}$ be a faithfully $G$-graded fusion 2-category $\mathfrak{C}$, and $g\in G$. Then, the finite semisimple $\mathfrak{C}_e$-$\mathfrak{C}_e$-bimodule 2-category $\mathfrak{C}_g$ witnesses a Morita equivalence between $\mathfrak{C}_e$ and itself, that is the canonical monoidal 2-functor  $\mathfrak{C}_e^{mop}\rightarrow \mathbf{End}_{\mathfrak{C}_e}(\mathfrak{C}_g)$ is an equivalence.
\end{Proposition}
\begin{proof}
Let us fix a simple object $D$ of $\mathfrak{C}_g$. We use $^{\sharp}D$ to denote the left dual to $D$. Then, $D\Box {^{\sharp}D}$ is an algebra in $\mathfrak{C}_e$. Further, it follows from example 4.1.3 and theorem 5.1.5 of \cite{D4} that the left $\mathfrak{C}$-module 2-functor $\mathfrak{C}\rightarrow\mathbf{Mod}_{\mathfrak{C}}(D\Box {^{\sharp}D})$ given by $C\mapsto C\Box {^{\sharp}D}$ is an equivalence. The proof of theorem 5.1.2 of \cite{D8} establishes that the canonical right action of $\mathfrak{C}$ on itself induces a monoidal 2-functor $F:\mathfrak{C}\rightarrow \mathbf{Bimod}_{\mathfrak{C}}(D\Box {^{\sharp}D})$ by $C\mapsto D\Box C\Box {^{\sharp}D}$. In addition, theorem 5.1.2 of \cite{D8} establishes that $F$ is an equivalence. In particular, using the fact that $\mathfrak{C}$ is $G$-graded, we find that restricting $F$ induces an equivalence $\mathfrak{C}_e\rightarrow \mathbf{Bimod}_{\mathfrak{C}_e}(D\Box {^{\sharp}D})$ of fusion 2-categories. 

Now, the left $\mathfrak{C}$-module 2-functor $\mathfrak{C}\rightarrow\mathbf{Mod}_{\mathfrak{C}}(D\Box {^{\sharp}D})$ induces a left $\mathfrak{C}_e$-module 2-functor $\mathfrak{C}_g\rightarrow\mathbf{Mod}_{\mathfrak{C}_e}(D\Box {^{\sharp}D})$, which is still an equivalence. Under this equivalence, the canonical right $\mathfrak{C}_e$-module structure on $\mathfrak{C}_g$ corresponds to the monoidal 2-functor $\mathfrak{C}_e\rightarrow \mathbf{Bimod}_{\mathfrak{C}_e}(D\Box {^{\sharp}D})$ given by $C\mapsto D\Box C\Box {^{\sharp}D}$. But we have seen above that this 2-functor is an equivalence. This concludes the proof of the result.
\end{proof}

\subsection{Examples of Group-Garded Fusion 2-Categories}\label{subsection:selfdualobj}

We now give examples of group-graded fusion 2-categories, along with descriptions of the underlying finite semisimple 2-categories, and their fusion rules.
\begin{Example}\label{ex:Z/2trivTambaraYamagami}
Let $H\cong \mathbb{Z}/2$ be a non-normal subgroup of $G=D_8$. Then, $\mathfrak{C}(D_8,H,triv, triv)$ is $\mathbb{Z}/2$-graded thanks to corollary \ref{cor:grouptheoreticalgraded} as $H^G\cong\mathbb{Z}/2\oplus\mathbb{Z}/2$. Writing $D_8/(\mathbb{Z}/2\oplus\mathbb{Z}/2)\cong \{+, -\}$, it follows from proposition \ref{prop:2categorycosetdecomposition} that $\mathfrak{C}(D_8,H,triv, triv)_-\simeq \mathbf{2Vect}$, and we write $D$ for the corresponding simple object. Further, $\mathfrak{C}(D_8,H,triv, triv)_+\simeq \mathfrak{C}(H^G,H,triv, triv) \simeq \mathbf{2Rep}(\mathbb{Z}/2\lbrack 1\rbrack \times \mathbb{Z}/2\lbrack 0\rbrack)$ has four simple objects, which we denote by $I$, $X$, $J$, and $Y$. The diagram below depicts the underlying 2-category $\mathfrak{C}(D_8,H,triv, triv)$, and the table gives its fusion rules.

$$\begin{tabular}{c c}
$\begin{tikzcd}
I \arrow[d, "\mathbf{Vect}"', bend right] \arrow["\mathbf{Vect}_{\mathbb{Z}/2}"', loop, distance=2em, in=125, out=55]  &                         & J \arrow[d, "\mathbf{Vect}"', bend right] \arrow["\mathbf{Vect}_{\mathbb{Z}/2}"', loop, distance=2em, in=125, out=55]  \\
X \arrow[u, "\mathbf{Vect}"', bend right] \arrow["\mathbf{Vect}_{\mathbb{Z}/2}"', loop, distance=2em, in=305, out=235] &                        & Y \arrow[u, "\mathbf{Vect}"', bend right] \arrow["\mathbf{Vect}_{\mathbb{Z}/2}"', loop, distance=2em, in=305, out=235] \\
& D \arrow["\mathbf{Vect}"', loop, distance=2em, in=305, out=235] &        
\end{tikzcd}$ & \begin{tabular}{ |c| c c c c c| } 
 \hline
 $\ $ & $I$ & $X$ & $J$ & $Y$ & $D$\\ 
 \hline
 $I$ & $I$ & $X$ & $J$ & $Y$ & $D$\\ 
 $X$ & $X$ & $2X$ & $Y$ & $2Y$ & $2D$\\ 
 $J$ & $J$ & $Y$ & $I$ & $X$ & $D$\\ 
 $Y$ & $Y$ & $2Y$ & $X$ & $2X$ & $2D$\\ 
 $D$ & $D$ & $2D$ & $D$ & $2D$ & $X\boxplus Y$\\ 
 \hline
\end{tabular}
\end{tabular}$$

\noindent We note that this example has already appeared in the physics literature \cite{BBSNT1,BBFP:I,BBFP,BBSNT}. 
\end{Example}

\begin{Example}
Let us consider $\mathbb{Z}/2\subseteq S_3$. It follows from proposition \ref{prop:2categorycosetdecomposition} that $\mathfrak{C}(S_3,\mathbb{Z}/2,triv, triv)\simeq \mathbf{2Rep}(\mathbb{Z}/2)\boxplus\mathbf{2Vect}$. However, this fusion 2-category is not $\mathbb{Z}/2$-graded as the smallest normal subgroup of $S_3$ containing $\mathbb{Z}/2$ is $S_3$. This last fact is also reflected in the fusion rules, which were first described in \cite{BBSNT}. Nevertheless, we may view $\mathfrak{C}(S_3,\mathbb{Z}/2,triv, triv)$ as an example of a near-group fusion 2-category. The diagram below gives the fusion rules.

$$\begin{tabular}{c c}
$\begin{tikzcd}[sep=small]
I \arrow[rr, "\mathbf{Vect}"', bend right] \arrow["\mathbf{Vect}_{\mathbb{Z}/2}"', loop, distance=2em, in=215, out=145] &  & X \arrow[ll, "\mathbf{Vect}"', bend right] \arrow["\mathbf{Vect}_{\mathbb{Z}/2}"', loop, distance=2em, in=35, out=325] \\ \\
 & D \arrow["\mathbf{Vect}"', loop, distance=2em, in=305, out=235] &
\end{tikzcd}$ & \begin{tabular}{ |c| c c c | } 
 \hline
 $\ $ & $I$ & $X$ & $D$\\ 
 \hline
 $I$ & $I$ & $X$ & $D$\\ 
 $X$ & $X$ & $2X$ & $2D$\\ 
 $D$ & $D$ & $2D$ & $X\boxplus D$\\ 
 \hline
\end{tabular}
\end{tabular}$$
\end{Example}

\begin{Example}
Let $Z\cong \mathbb{Z}/2$ be a normal subgroup of $D_8$. This forces $Z$ to be the center of $D_8$. Then, $\mathfrak{C}(D_8,Z,triv, triv)$ is a $\mathbb{Z}/2\oplus\mathbb{Z}/2$-graded fusion 2-category. Further, it follows from example \ref{ex:gauginZ/2inZ/4} above that $\mathfrak{C}(D_8,Z,triv, triv)\simeq \mathbf{2Vect}_{\mathcal{G}}^{\pi}$ with $\mathcal{G}=\mathbb{Z}/2\lbrack 1\rbrack\times(\mathbb{Z}/2\oplus\mathbb{Z}/2)\lbrack 0\rbrack$, and $\pi$ a 4-cocycle. This holds because every connected component of $\mathfrak{C}(D_8,Z,triv, triv)$ contains a unique invertible object. In addition, it also follows from example \ref{ex:gauginZ/2inZ/4} that the Postnikov data for $\mathcal{G}$ is trivial. This example has already appeared in \cite{BBFP, BBSNT}.
\end{Example}

\begin{Example}
Let us consider $H=\mathbb{Z}/2$ as a subgroup of the alternating group $A_4$. Then, the smallest normal subgroup containing $H$ is isomorphic to $\mathbb{Z}/2\oplus\mathbb{Z}/2$, so that $\mathfrak{C}(A_4,\mathbb{Z}/2,triv, triv)$ is a faithfully $\mathbb{Z}/3$-graded fusion 2-category. Further, it follows from proposition \ref{prop:2categorycosetdecomposition} that $\mathfrak{C}(A_4,\mathbb{Z}/2,triv, triv)\simeq \mathbf{2Rep}(\mathbb{Z}/2)\boxplus \mathbf{2Vect}\boxplus \mathbf{2Vect}$. We see that $\mathfrak{C}(A_4,\mathbb{Z}/2,triv, triv)$ has interesting fusion rules. More precisely, the simple object $T$ is a sort of ``self-triality'' object. Its connected components and fusion rules are given below.

$$\begin{tabular}{c c}
$\begin{tikzcd}[sep=small]
I \arrow[rr, "\mathbf{Vect}"', bend right] \arrow["\mathbf{Vect}_{\mathbb{Z}/2}"', loop, distance=2em, in=215, out=145] &  & X \arrow[ll, "\mathbf{Vect}"', bend right] \arrow["\mathbf{Vect}_{\mathbb{Z}/2}"', loop, distance=2em, in=35, out=325] \\
&  &\\
T \arrow["\mathbf{Vect}"', loop, distance=2em, in=215, out=145]&  & T' \arrow["\mathbf{Vect}"', loop, distance=2em, in=35, out=325]
\end{tikzcd}$ & \begin{tabular}{ |c| c c c c | } 
 \hline
 $\ $ & $I$ & $X$ & $T$ & $T'$\\ 
 \hline
 $I$ & $I$ & $X$ & $T$ & $T'$\\ 
 $X$ & $X$ & $2X$ & $2T$ & $2T'$\\ 
 $T$ & $T$ & $2T$ & $T'$ & $X$\\ 
 $T'$ & $T'$ & $2T'$ & $X$ & $T$\\
 \hline
\end{tabular}
\end{tabular}$$
\end{Example}

\subsection{Constructing Tambara-Yamagami defects}\label{sub:TYdefect}

\begin{Definition}
A simple object $D$ of a fusion 2-category $\mathfrak{C}$ is a \textit{Tambara-Yamagami defect} if $\mathfrak{C}$ admits a (necessarily faithful) $\mathbb{Z}/2$-grading $\mathfrak{C}=\mathfrak{C}_+\boxplus\mathfrak{C}_-$ such that $\mathfrak{C}_-\simeq \mathbf{2Vect}$ is generated by $D$.
\end{Definition}

By proposition \ref{prop:gradedMoritaequivalence}, the Tambara-Yamagami defect witnesses a self-duality of $\fC_{+}$, i.e.\ a Morita equivalence of $\fC_{+}$ with itself. We now examine the most elementary example of such a fusion 2-category.  

\begin{Example}
The most basic example of a fusion 2-category with a Tambara-Yamagami defect is $\mathfrak{C}(\mathbb{Z}/2,0,triv,triv)\simeq \mathbf{2Vect}_{\mathbb{Z}/2}$. In fact, it follows from \cite{JFY} that every fusion 2-category with an invertible Tambara-Yamgami defect is of this form.
\end{Example}

A priori a Tambara-Yamagami defect is additional data on a fusion 2-category. It turns out that this is not the case as we show in the next lemma, i.e.\ having a Tambara-Yamagami defect is really a property of a fusion 2-category.

\begin{Lemma}
If a fusion 2-category admits a Tambara-Yamagami defect, then it is unique.
\end{Lemma}
\begin{proof}
Let $\mathfrak{C}$ be a fusion 2-category. Let $D$ and $D'$ be two Tambara-Yamagami defect in $\mathfrak{C}$. In particular, $\mathfrak{C}$ admits two faithful $\mathbb{Z}/2$-gradings $\mathfrak{C}=\mathfrak{C}_+\boxplus\mathfrak{C}_-$and $\mathfrak{C}=\mathfrak{C}_0\boxplus\mathfrak{C}_1$with $\mathfrak{C}_-\simeq \mathbf{2Vect}$ generated by $D$ and $\mathfrak{C}_1\simeq \mathbf{2Vect}$ generated by $D'$. If $D$ and $D'$ are not equivalent, then $D$ lies in $\mathfrak{C}_0$ and $D'$ in $\mathfrak{C}_+$, so that $D\Box D'$ is a non-zero object in both $\mathfrak{C}_-$ and $\mathfrak{C}_1$. As $\mathfrak{C}_-$ and $\mathfrak{C}_1$ are both connected finite semisimple 2-categories, this forces $\mathfrak{C}_-=\mathfrak{C}_1$, so that $D\simeq D'$, a contradiction. 
\end{proof}

Let $G$ be a (finite) group. A $\Z/2$-grading $G=G_+\sqcup G_-$ on $G$ is precisely the data of a homomorphism $f:G\rightarrow\Z/2$. In particular, we have $G_+=f^{-1}(0)$ and $G_-=f^{-1}(1)$, and such a grading is faithful if and only if $f$ is surjective. Using our previous results together with this notion, we obtain the following proposition.

\begin{Proposition}\label{prop:selfdualityF2Cs}
Every fusion 2-category with a Tambara-Yamagami defect is equivalent to $\mathfrak{C}(G,H,\pi,\psi)$, with $G = G_+\sqcup G_-$ a faithfully $\mathbb{Z}/2$-graded group, $H$ a subgroup of $G_+$ such that $H\backslash G_-/H = *$, and such that, for any $g\in G_-$, $H$ and $gHg^{-1}$ form an exact factorization of $G_+$. In addition, $\pi$ is a 4-cocycle on $G$, and $\psi$ is a 3-cochain on $H$ such that $d\psi =\pi|_H$.
\end{Proposition}
\begin{proof}
It follows from the definition and theorem \ref{thm:grouptheoreticalrecognition} that any fusion 2-category with a Tambara-Yamagami defect is group-theoretical. Namely, a left $\mathfrak{C}(G,H,\pi,\psi)_+$-module structure on $\mathfrak{C}(G,H,\pi,\psi)_-\simeq\mathbf{2Vect}$ is precisely the data of a monoidal 2-functor $\mathfrak{C}(G,H,\pi,\psi)_+\rightarrow \mathbf{End}(\mathbf{2Vect})\simeq \mathbf{2Vect}$ (see lemma 2.1.5 of \cite{D4}). Further, corollary \ref{cor:grouptheoreticalgraded} gives that the group-theoretical fusion 2-category $\mathfrak{C}(G,H,\pi,\psi)$ is graded by $G/H^G$. Thence, we find that $G_+= H^G$ is a subgroup of $G$ of index 2. This implies that $G = G_+\sqcup G_-$ is a $\mathbb{Z}/2$-graded group. Finally, proposition \ref{prop:2categorycosetdecomposition} yields that $$\mathfrak{C}(G,H,\pi,\psi)_- \simeq \boxplus_{\lbrack g\rbrack\in H\backslash G_- /H} \mathbf{Mod}(\mathbf{Vect}^{\xi_g}_{H\cap gHg^{-1}}).$$ This finite semisimple 2-category is connected if and only if $H\backslash G_-/H = *$. Furthermore, it is equivalent to $\mathbf{2Vect}$ if and only if, in addition, $H\cap gHg^{-1} = \{e\}$ for any $g\in G_-$. But, provided that $H\backslash G_-/H = *$, we have $H\cap gHg^{-1} = \{e\}$ if and only if $H$ and $gHg^{-1}$ form an exact factorization of $G_+$.
\end{proof}

The most fundamental examples of fusion 2-categories with a Tambara-Yamagami defect are given by the construction that we describe in the next example. As we will shortly see, the fusion rules for the Tambara-Yama\-gami defect make it into a noninvertible operator. Therefore, we will have algebraically established the existence of a noninvertible operator arising as a condensation, without the need to write any background fields. For constructions that take this latter approach, see \cite{RSS22,IO23}. In \cite[Appendix D]{KOZ21}, the authors construct a specific Kramers-Wannier-like defect through gauging, and show that it gives a self-duality of the theory. From our point of view, the self-duality of the theory gives rise to a Tambara-Yamagami defect of the corresponding fusion 2-category.

\begin{Example}\label{ex:normalselfdualF2C}
Let $H$ be a finite group. There is an order 2 automorphism $\sigma:H\times H\cong H\times H$ given by swapping the two copies of $H$. We can therefore consider the semi-direct product $(H\times H)\rtimes\mathbb{Z}/2$, commonly known as the wreath product $H \wr \mathbb{Z}/2$. Slightly abusing notations, we use $H\times \{e\}$ to denote the subgroup $(H\times \{e\})\rtimes \{0\}\subset H \wr \mathbb{Z}/2$. The fusion 2-category $\mathfrak{C}(H \wr \mathbb{Z}/2,H\times \{e\},triv, triv)$ is an elementary example of a fusion 2-category with a Tambara-Yamagami defect. In this case, we have $$\mathfrak{C}(H \wr \mathbb{Z}/2,H\times \{e\},triv, triv)\simeq \Big(\mathbf{2Rep}(H)\boxtimes\mathbf{2Vect}_H \Big) \boxplus \mathbf{2Vect}.$$ So proposition \ref{prop:gradedMoritaequivalence} gives an object witnessing the Morita autoequivalence that exchanges $\textbf{2}\rRep(H)$ and $\textbf{2}\Vect_H$.
More generally, such examples can be twisted using a 4-cocycle $\pi$ in $Z^4_{gp}(H \wr \mathbb{Z}/2;\mathds{k}^{\times})$, together with a 3-cochain $\psi$ for $H\times \{e\}$ trivializing $\pi|_{H\times \{e\}}$.
\end{Example}

In fact, it is possible to recognize the fusion 2-categories with a Tambara-Yamagami defect that fall within the scope of example \ref{ex:normalselfdualF2C}.

\begin{Corollary}\label{cor:normalselfduality}
If we assume that $H$ is a normal subgroup of $G_+$ in the statement of proposition \ref{prop:selfdualityF2Cs}, then $G\cong H\wr \mathbb{Z}/2$.
\end{Corollary}
\begin{proof}
If $H$ is a normal subgroup of $G_+$, then so is $gHg^{-1}$ for any $g\in G_-$. This implies that $G_+\cong H\times gHg^{-1}$ as $H$ and $gHg^{-1}$ give an exact factorization of $G_+$. Thence, $G$ is an extension of $H\times H$ by $\mathbb{Z}/2$ with action $\mathbb{Z}/2\rightarrow Out(H\times H)$ given by the map $\sigma$ swapping the two factors. But, it was proven in \cite{EML} that such extensions, if they exist, are classified by a torsor over $H^2_{gp}(\mathbb{Z}/2; (Z(H)\times Z(H))_{\sigma})$, a group which is trivial as can be seen from an elementary direct computation. This proves that $G\cong H\wr \mathbb{Z}/2$.
\end{proof}

We now examine the fusion rule of the Tambara-Yamagami defect. We will do so using the notations introduced in section \ref{sub:fusionrules}.

\begin{Proposition}\label{prop:fusionD}
Let $\mathfrak{C}$ be a fusion 2-category with a Tambara-Yamagami defect. Writing $D$ for the simple object of $\mathfrak{C}_-$, we have $$D\Box D \simeq \boxplus_{f\in gHg^{-1}}\mathbf{V}_{\lbrack f\rbrack},$$ with $g$ any element of $G_-$.
\end{Proposition}
\begin{proof}
Let $g$ be any element of $G_-$. We necessarily have $D\simeq \mathbf{V}_{\lbrack g\rbrack}$. In particular, we also have $D\simeq \mathbf{V}_{\lbrack g^{-1}\rbrack}$. Thence, it follows from proposition \ref{prop:partialfusionrule} that $$D\Box D \simeq \mathbf{V}_{\lbrack g\rbrack}\Box \mathbf{V}_{\lbrack g^{-1}\rbrack}\simeq\boxplus_{f\in gHg^{-1}}\mathbf{V}_{\lbrack f\rbrack},$$ as desired.
\end{proof}

\subsection{Tambara-Yamagami 2-Categories}\label{subsection:TYF2C}

\begin{Definition}
A Tambara-Yamagami 2-category is a $\mathbb{Z}/2$-graded fusion 2-category $\mathfrak{C}:=\mathfrak{C}_+\boxplus\mathfrak{C}_-$ such that 
$\mathfrak{C}_+\simeq \mathbf{2Vect}_{\mathcal{G}}$ for some finite 2-group $\mathcal{G}$, and $\mathfrak{C}_-\simeq\mathbf{2Vect}$.
\end{Definition}

The next lemma gives strong restriction on which finite 2-group can occur as the trivially graded factor of a Tambara-Yamagami 2-category. This is a 2-categorical analogue of the well-known fact that Tamabara-Yamagami 1-categories only exists for finite abelian groups.

\begin{Lemma}
Let $\mathcal{G}$ be a finite 2-group. There exists an invertible $\mathbf{2Vect}_{\mathcal{G}}$-$\mathbf{2Vect}_{\mathcal{G}}$-bimodule 2-category whose underlying finite semisimple 2-category is $\mathbf{2Vect}$ if and only if $\mathcal{G}\simeq \widehat{A}\lbrack 1\rbrack\times A\lbrack 0\rbrack$ for a finite abelian group $A$.
\end{Lemma}
\begin{proof}
Proposition 4.3.3.10 of \cite{D:thesis} shows that the Morita dual of $\mathbf{2Vect}_{\mathcal{G}}$ with respect to $\mathbf{2Vect}$ (for any module structure) is $\mathbf{2Rep}(\mathcal{G})$. It is therefore sufficient to understand when there is a monoidal equivalence $\mathbf{2Vect}_{\mathcal{G}}\simeq \mathbf{2Rep}(\mathcal{G})$. Comparing the connected components of the identity, we find that $\mathbf{Vect}_{\pi_2(\mathcal{G})}\simeq \mathbf{Rep}(\pi_1(\mathcal{G}))$ as braided fusion 1-categories, so that $\pi_2(\mathcal{G})\simeq \pi_1(\mathcal{G})$. We denote this finite abelian group by $A$.

It remains to show that both the action of $\pi_1(\mathcal{G})$ on $\pi_2(\mathcal{G})$, and the 3-cocycle $\omega$ are trivial. But, it follows from the above discussion that the monoidal 1-category associated to $\mathcal{G}$ admits a symmetric structure. It therefore follows by inspection that the action of $\pi_1(\mathcal{G})$ on $\pi_2(\mathcal{G})$ is trivial. Further, by the existence of the symmetric structure above, theorem B of \cite{JO} applies and establishes that the the Postnikov data of $\mathcal{G}$ is trivial as claimed.
\end{proof}

To see this from the point of view of gauging, note that if we set $G = \widehat{A}[0] \times A[0]$ (which is non-canonically equivalent to ${A}[0] \times A[0]$) then equations $\eqref{eq:bim1}$ and $\eqref{eq:bim2}$ are equivalent. The fact that $\textbf{2}\Vect_{\cG}\simeq \textbf{2}\rRep(\cG)$ is what preserves a “self-duality” in the $\fC_+$ component.  Indeed, this is because both the 0-form and the 1-form symmetries are given by $A$ and therefore the gauging procedure that swaps the 0-form and 1-form symmetry leaves $\mathcal{G}$ invariant.

\begin{Proposition}\label{prop:2TYclassification}
There is a bijection between the set of equivalence classes of Tambara-Yamagami 2-categories and the quotient of the set of pairs consisting of a finite abelian group $A$ together with a 4-cocycle $\pi$ for $A\wr\mathbb{Z}/2$ with coefficients in $\mathds{k}^{\times}$, that is trivializable on $A\oplus A$, under the relation identifying two pairs $(A_1,\pi_1)$ and $(A_2,\pi_2)$ whenever there exists a group homomorphism $f:A_1\wr \mathbb{Z}/2\cong A_2\wr \mathbb{Z}/2$ such that $f(A_1\oplus 0) = A_2\oplus 0$, and $f^*\pi_2/\pi_1$ is a coboundary. We write $\mathbf{2TY}(A,\pi)$ for a Tambara-Yamagami fusion 2-category representing this equivalence class.\footnote{An explicit fusion 2-category representing this class is $\mathfrak{C}(A\wr\mathbb{Z}/2,A\oplus 0,\pi,\psi)$, where $\psi$ is a 3-cochain for $A\oplus 0$ with value in $\mathds{k}^{\times}$ such that $d\psi=\pi|_{A\oplus 0}$.} 
\end{Proposition}
\begin{proof}
Let $\mathfrak{T}$ be a Tambara-Yamagami 2-category so that $\mathfrak{T}_+\simeq \mathbf{2Vect}_{\mathcal{G}}$. Thence, the last lemma above yields that $\mathcal{G}\simeq A\lbrack 1\rbrack \times A \lbrack 0\rbrack$ for some finite abelian group $A$. But, by proposition \ref{prop:selfdualityF2Cs}, we also have that $\mathfrak{T}$ is equivalent as a fusion 2-category to $\mathfrak{C}(G,A,\pi,\psi)$ for some $\mathbb{Z}/2$-graded group $G=G_+\sqcup G_-$. As $\mathfrak{T}_+\simeq \mathbf{2Vect}_{A\lbrack 1\rbrack \times A \lbrack 0\rbrack}$, it follows from proposition \ref{prop:2categorycosetdecomposition} that $A$ must be a normal subgroup of $G_+$. Appealing to corollary \ref{cor:normalselfduality}, we therefore obtain that $\mathfrak{T}\simeq \mathfrak{C}(A\wr \mathbb{Z}/2,A\oplus 0,\pi,\psi)$. It then also follows that $$\mathbf{2Vect}_{A\lbrack 1\rbrack \times A \lbrack 0\rbrack}\simeq\mathfrak{T}_+\simeq\mathfrak{C}(A\oplus A,A\oplus 0,\pi|_{A\oplus A},\psi),$$ so that $\pi|_{A\oplus A}$ must be trivializable by propositions \ref{prop:2groupgraded2vectorspaces} and \ref{prop:classificationgrouptheoretical}. We may then, without loss of generality, assume that $\pi|_{A\oplus A} = triv$, which implies that $\psi$ is a 3-cochain for $A\oplus 0$. It only remains to show that the monoidal equivalence class of the fusion 2-category $\mathfrak{C}(A\wr \mathbb{Z}/2,A\oplus 0,\pi,\psi)$ does not depend on the choice of the 3-cochain $\psi$. Thanks to proposition \ref{prop:classificationgrouptheoretical}, in order to do so, it is enough to show that there exists $\xi \in Z^3(A\wr\mathbb{Z}/2;\mathds{k}^{\times})$ such that $\psi/(\xi|_{A\oplus 0})$ is a 3-coboundary for $A\oplus 0$. This is precisely the content of lemma \ref{lem:technical2TY} below, whence the proof is complete.
\end{proof}

\begin{Lemma}\label{lem:technical2TY}
For any finite abelian group $A$, the homomorphism $$H^3_{gp}(A\wr\mathbb{Z}/2;\mathds{k}^{\times})\rightarrow H^3_{gp}(A;\mathds{k}^{\times})$$ given by the pullback along the inclusion $A\oplus 0\hookrightarrow A\wr\mathbb{Z}/2$ is surjective.
\end{Lemma}
\begin{proof}
Equivalently, we show that the map $H_{gp}^4(A\wr\mathbb{Z}/2;\mathbb{Z})\rightarrow H_{gp}^4(A;\mathbb{Z})$ is surjective. Namely, as $\mathds{k}$ is an algebraically closed field of characteristic zero, there is a natural isomorphism $H^3_{gp}(G;\mathds{k}^{\times})\cong H^4_{gp}(G;\mathbb{Z})$. Let us begin by observing that the inclusion $A\oplus 0\hookrightarrow A\wr\mathbb{Z}/2$ factors as $A\oplus 0\hookrightarrow A\oplus A\hookrightarrow A\wr\mathbb{Z}/2$. Moreover, by the K\"unneth formula, we have $$H_{gp}^4(A\oplus A;\mathbb{Z})\cong \bigoplus_{i=0}^4 \big(H_{gp}^i(A;\mathbb{Z})\otimes H_{gp}^{4-i}(A;\mathbb{Z})\big)\oplus \bigoplus_{i=0}^5Tor_1^{\mathbb{Z}}\big(H_{gp}^i(A;\mathbb{Z})\otimes H_{gp}^{5-i}(A;\mathbb{Z})\big).$$ It was shown in section 4 of \cite{Lea} that, under the K\"unneth isomorphism above, the image of the map $$H_{gp}^4(A\wr\mathbb{Z}/2;\mathbb{Z})\rightarrow H_{gp}^4(A\oplus A;\mathbb{Z})$$ given by the pullback along the inclusion $A\oplus A\hookrightarrow A\wr\mathbb{Z}/2$ contains all the 4-cocycles of the from $\alpha \otimes \beta + \beta\otimes \alpha$ with $\alpha\in H_{gp}^i(A;\mathbb{Z})$ and $\beta\in H_{gp}^{4-i}(A;\mathbb{Z})$ with $i=0,...,4$. But, for any $\alpha\in H_{gp}^4(A;\mathbb{Z})$, it follows by inspection that the image of $\alpha\otimes 1 + 1\otimes \alpha$ under the map $H_{gp}^4(A\oplus A;\mathbb{Z})\rightarrow H_{gp}^4(A\oplus 0;\mathbb{Z})$ is just $\alpha$. This concludes the proof of the lemma.
\end{proof}

\begin{Example}\label{ex:twoTYF2Cs}
Let us take $A=\mathbb{Z}/2$, so that $G:=\mathbb{Z}/2\wr\mathbb{Z}/2\cong D_8$. We have $H^4_{gp}(D_8;\mathds{k}^{\times})\cong \mathbb{Z}/2\oplus \mathbb{Z}/2$. Further, the map $H^4_{gp}(D_8;\mathds{k}^{\times})\rightarrow H^4_{gp}(\mathbb{Z}/2\oplus \mathbb{Z}/2;\mathds{k}^{\times})$ is non-zero, and has kernel isomorphic to $\mathbb{Z}/2$. Let us write $\zeta$ for a 4-cochain representing the non-trivial element in this kernel. The Tambara-Yamagami 2-category $\mathbf{2TY}(\mathbb{Z}/2,triv)$ was examined in example \ref{ex:Z/2trivTambaraYamagami}. We find that there is another Tambara-Yamagami 2-category with $A=\mathbb{Z}/2$, namely $\mathbf{2TY}(\mathbb{Z}/2,\zeta)$.
\end{Example}

\begin{Example}\label{ex:beyondTamabraYamagami}
Let us take $A=\mathbb{Z}/2$, so that $G:=\mathbb{Z}/2\wr\mathbb{Z}/2\cong D_8$. Let $\kappa$ be a 4-cochain for $D_8$ whose image under $H^4_{gp}(D_8;\mathds{k}^{\times})\rightarrow H^4_{gp}(\mathbb{Z}/2\oplus \mathbb{Z}/2;\mathds{k}^{\times})$ is non-zero. We obtain a fusion 2-category $\mathfrak{C}(D_8,\mathbb{Z}/2,\kappa,triv)$ with a Tambara-Yamagami defect, whose underlying fusion 2-category and fusion rules are depicted below. We note that the trivially graded component is the fusion 2-category examined in example \ref{ex:cocyclesgivesdistinctf2cs} with $A=\mathbb{Z}/2$. In fact, there are two possibilities for the class of $\kappa$ in cohomology, but we suspect that the two associated fusion 2-categories are equivalent. In example \ref{ex:distinctkappa}, we will explain the relation between these fusion 2-categories and the one considered in example 3.5.2 of \cite{BBFP}.

$$\begin{tikzcd}
 & I \arrow[ld, "\mathbf{Vect}"', bend right] \arrow["\mathbf{Mod}(\mathbf{Vect}_{\mathbb{Z}/2})"', loop, distance=2em, in=125, out=55] &\\
X \arrow[ru, "\mathbf{Vect}"', bend right] \arrow["\mathbf{Mod}(\mathbf{Vect}_{\mathbb{Z}/2})"', loop, distance=2em, in=215, out=145] &                                    & Y \arrow["\mathbf{Mod}(\mathbf{Vect}_{\mathbb{Z}/2}^{\omega})"', loop, distance=2em, in=35, out=325] \\
 & D \arrow["\mathbf{Vect}"', loop, distance=2em, in=305, out=235]&
 \end{tikzcd}$$

 \begin{center}
\begin{tabular}{ |c| c c c c| } 
 \hline
 $\ $ & $I$ & $X$ & $Y$ & $D$\\ 
 \hline
 $I$ & $I$ & $X$ & $Y$ & $D$\\ 
 $X$ & $X$ & $2X$ & $2Y$ & $2D$\\ 
 $Y$ & $Y$ & $2Y$ & $2X$ & $2D$\\ 
 $D$ & $D$ & $2D$ & $2D$ & $X\boxplus Y$\\ 
 \hline
\end{tabular}
\end{center}
\end{Example}

\begin{Example}
Let us take $A=\mathbb{Z}/3$. We then have $\mathbb{Z}/3\wr\mathbb{Z}/2\cong S_3\times \mathbb{Z}/3$. It follows from the K\"unneth formula that $H^4_{gp}(S_3\times \mathbb{Z}/3;\mathds{k}^{\times})\cong \mathbb{Z}/3$. Further, the canonical map $H^4_{gp}(S_3\times \mathbb{Z}/3;\mathds{k}^{\times})\rightarrow H^4_{gp}(\mathbb{Z}/3\times \mathbb{Z}/3;\mathds{k}^{\times})$ is injective. Thus, $\mathbf{2TY}(\mathbb{Z}/3,triv)$ is the only Tambara-Yamagami fusion 2-category with $A=\mathbb{Z}/3$. The underlying 2-category and fusion rules are given below.

$$
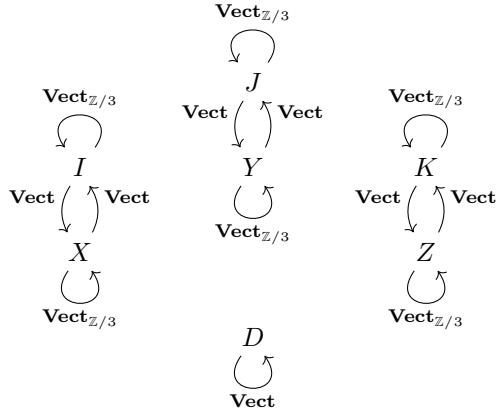
\begin{tikzcd}
 &  & J \arrow[d, "\mathbf{Vect}"', bend right] \arrow["\mathbf{Vect}_{\mathbb{Z}/3}"', loop, distance=2em, in=125, out=55]  &  &  \\
I \arrow[d, "\mathbf{Vect}"', bend right] \arrow["\mathbf{Vect}_{\mathbb{Z}/3}"', loop, distance=2em, in=125, out=55]  &  & Y \arrow[u, "\mathbf{Vect}"', bend right] \arrow["\mathbf{Vect}_{\mathbb{Z}/3}"', loop, distance=2em, in=305, out=235] &  & K \arrow["\mathbf{Vect}_{\mathbb{Z}/3}"', loop, distance=2em, in=125, out=55] \arrow[d, "\mathbf{Vect}"', bend right]  \\
X \arrow[u, "\mathbf{Vect}"', bend right] \arrow["\mathbf{Vect}_{\mathbb{Z}/3}"', loop, distance=2em, in=305, out=235] &  &   &  & Z \arrow["\mathbf{Vect}_{\mathbb{Z}/3}"', loop, distance=2em, in=305, out=235] \arrow[u, "\mathbf{Vect}"', bend right] \\
  &  & D \arrow["\mathbf{Vect}"', loop, distance=2em, in=305, out=235]   & &                                          
\end{tikzcd}$$

\begin{center}
\begin{tabular}{ |c| c c c c c c c| } 
 \hline
 $\ $ & $I$ & $X$ & $J$ & $Y$ & $K$ & $Z$ & $D$\\ 
 \hline
 $I$ & $I$ & $X$ & $J$ & $Y$ & $K$ & $Z$ & $D$\\ 
 $X$ & $X$ & $3X$ & $Y$ & $3Y$ & $3Z$ & $3Z$ & $3D$\\ 
 $J$ & $J$ & $Y$ & $K$ & $Z$ & $I$ & $X$ & $D$\\ 
 $Y$ & $Y$ & $3Y$ & $Z$ & $3Z$ & $X$ & $3X$ & $3D$\\ 
 $K$ & $K$ & $X$ & $I$ & $X$ & $J$ & $Y$ & $D$\\ 
 $Z$ & $Z$ & $3X$ & $X$ & $3X$ & $Y$ & $3Y$ & $3D$\\ 
 $D$ & $D$ & $3D$ & $D$ & $3D$ & $D$ & $3D$ & $X\boxplus Y\boxplus Z$\\ 
 \hline
\end{tabular}
\end{center}

On the other hand, if $\kappa$ denotes a 4-cocycle representing a non-trivial class in $H^4_{gp}(S_3\times \mathbb{Z}/3;\mathds{k}^{\times})$, then $\mathfrak{C}(\mathbb{Z}/3\wr\mathbb{Z}/2,\mathbb{Z}/2,\kappa,triv)$ is a fusion 2-category with a Tambara-Yamagami defect, which is not a Tambara-Yamagami 2-category. We have depicted its underlying 2-category below as well as its fusion rules.

$$
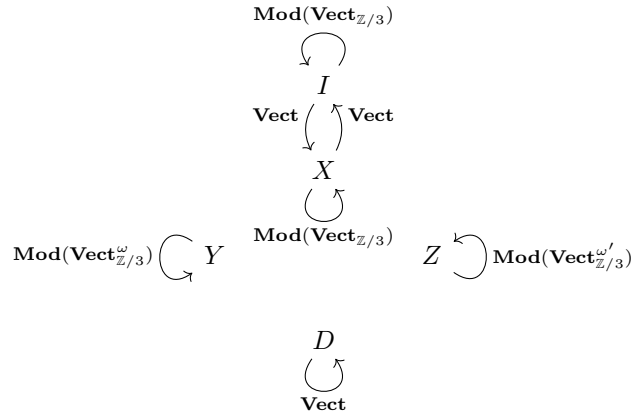
\begin{tikzcd}
 & I \arrow[d, "\mathbf{Vect}"', bend right] \arrow["\mathbf{Mod}(\mathbf{Vect}_{\mathbb{Z}/3})"', loop, distance=2em, in=125, out=55]  &\\
 & X \arrow[u, "\mathbf{Vect}"', bend right] \arrow["\mathbf{Mod}(\mathbf{Vect}_{\mathbb{Z}/3})"', loop, distance=2em, in=305, out=235] &\\
Y \arrow["\mathbf{Mod}(\mathbf{Vect}_{\mathbb{Z}/3}^{\omega})"', loop, distance=2em, in=215, out=145] & & Z \arrow["\mathbf{Mod}(\mathbf{Vect}_{\mathbb{Z}/3}^{\omega'})"', loop, distance=2em, in=35, out=325] \\
 & D \arrow["\mathbf{Vect}"', loop, distance=2em, in=305, out=235] &        
\end{tikzcd}$$

 \begin{center}
\begin{tabular}{ |c| c c c c c| } 
 \hline
 $\ $ & $I$ & $X$ & $Y$ & $Z$ & $D$\\ 
 \hline
 $I$ & $I$ & $X$ & $Y$ & $Z$ & $D$\\ 
 $X$ & $X$ & $3X$ & $3Y$ & $3Z$ & $3D$\\ 
 $Y$ & $Y$ & $3Y$ & $3Z$ & $3X$ & $3D$\\ 
 $Z$ & $Z$ & $3Z$ & $3X$ & $3Y$ & $3D$\\ 
 $D$ & $D$ & $3D$ & $3D$ & $3D$ & $X\boxplus Y\boxplus Z$\\ 
 \hline
\end{tabular}
\end{center}
\end{Example}

We now expand on the ideas of \S\ref{section:Fib2Functor} and examine when Tambara-Yamagami 2-categories admit fiber 2-functors. The decategorified question was considered in \cite{Tam}, and was used in \cite{Thorngren:2019iar,TW21}.

\begin{Proposition}
If $A$ is a finite abelian group of odd order, then there is exactly one Tambara-Yamagami fusion 2-category $\mathbf{2TY}(A,triv)$. Furthermore, $\mathbf{2TY}(A,triv)$ has a fiber 2-functor, and $\mathbf{2TY}(A,triv)$ is not equivalent to the 2-category of 2-representations of a 2-group if $A$ is non-trivial.
\end{Proposition}
\begin{proof}
As $A$ has odd order, it follows from corollary 6.1 of \cite{Lea} that the map $H^4_{gp}(A\wr \mathbb{Z}/2;\mathds{k}^{\times})\rightarrow H^4_{gp}(A\oplus A;\mathds{k}^{\times})$ is injective. This proves the first part of the claim. For the second part, note that, with $\Delta:=\{(a,a)|a\in A\}$, $\Delta\oplus\mathbb{Z}/2$ is a subgroup of $A\wr\mathbb{Z}/2$, and a complement of $A\oplus 0$. It therefore follows from theorem \ref{thm:fiber2functor} that $\mathbf{2TY}(A,triv)$ has a fiber 2-functor.

Observe that there is an isomorphism of groups $A\wr\mathbb{Z}/2\cong (A\rtimes \mathbb{Z}/2) \times \Delta$, where $\mathbb{Z}/2$ acts on $A$ by multiplication by $-1$. If $\mathbf{2TY}(A,triv)$ were the fusion 2-category of 2-representations of a finite 2-group, then $A\oplus 0$ would have an abelian normal complement $N$ in $A\wr \mathbb{Z}/2$. Now, as $A$ has odd order, $N$ must contain at least one element of order 2. Given that we have assumed that $N$ is normal, it follows from an elementary computation that $N$ must contain all the elements of order 2. In particular, $N$ must contain all of $A\rtimes \mathbb{Z}/2$, but this subgroup is not abelian.
\end{proof}

The above proposition gives a complete understanding of Tambara-Yamagami 2-categories associated to a finite group of odd order. By contrast, as we illustrate in the examples below, the situation is much more complicated when the finite group $A$ has order divisible by $2$.

\begin{Example}
Let $A=\mathbb{Z}/4$. Then, it follows from corollary 6.1 of \cite{Lea} that the kernel of the map $$H^4_{{gp}}(\mathbb{Z}/4\wr\mathbb{Z}/2;\mathds{k}^{\times})\rightarrow H^4_{{gp}}(\mathbb{Z}/4\oplus\mathbb{Z}/4;\mathds{k}^{\times})\oplus H^4_{{gp}}(\Delta\oplus\mathbb{Z}/2;\mathds{k}^{\times}),$$ whose components are the pullbacks along the inclusions, is $\mathbb{Z}/2$. Let $\xi$ denote a 4-cocycle for $\mathbb{Z}/4\wr\mathbb{Z}/2$ representing the non-trivial class in this kernel. Then, $\mathbf{2TY}(\mathbb{Z}/4,\xi)$ is a Tambara-Yamagami 2-category. Moreover, it admits a fiber 2-functor because $\xi|_{\Delta\oplus\mathbb{Z}/2}$ is trivial. Yet, the 4-cocycle $\xi$ for $\mathbb{Z}/4\wr\mathbb{Z}/2$ is not trivializable.
\end{Example}

\begin{Example}\label{ex:distinctkappa}
We now consider $A=\mathbb{Z}/2$, so that $G:=\mathbb{Z}/2\wr\mathbb{Z}/2\cong D_8$. Recall that we use $\zeta$ as in example \ref{ex:twoTYF2Cs} to denote a 4-cocycle representing the non-trivial class in the kernel of the map $H^4_{{gp}}(D_8;\mathds{k}^{\times})\rightarrow H^4_{{gp}}(\mathbb{Z}/2\oplus\mathbb{Z}/2;\mathds{k}^{\times})$. In this case, there are exactly two complements for $\mathbb{Z}/2\oplus 0$ given by the subgroups $\Delta\rtimes \mathbb{Z}/2\subset\mathbb{Z}/2\wr\mathbb{Z}/2$, which is isomorphic to $\mathbb{Z}/2\oplus\mathbb{Z}/2$, and $K:=\langle(1,0,1)\rangle$, which is isomorphic to $\mathbb{Z}/4$. It follows from corollary 6.1 of \cite{Lea} that the image of $\zeta$ under the map $H^4_{{gp}}(\mathbb{Z}/2\wr\mathbb{Z}/2;\mathds{k}^{\times})\rightarrow H^4_{gp}(\Delta\rtimes \mathbb{Z}/2;\mathds{k}^{\times})$ is non-zero. On the other hand, it follows from a direct computation that the map $H^4_{{gp}}(\mathbb{Z}/2\wr\mathbb{Z}/2;\mathds{k}^{\times})\rightarrow H^4_{gp}(K;\mathds{k}^{\times})$ is zero.

In particular, theorem \ref{thm:fiber2functor} establishes that $\mathbf{2TY}(\mathbb{Z}/2,\zeta)$ does admit a fiber 2-functor. However, we will argue that this last fusion 2-category is not equivalent to the 2-category of 2-representations of a finite 2-group. Namely, it follows from corollary \ref{cor:grouptheoretical2group} that the only possible 2-groups are $$\mathbb{Z}/4\lbrack 1\rbrack\rtimes\mathbb{Z}/2\lbrack 0\rbrack\ \ \ \mathrm{and}\ \ \ \mathbb{Z}/4\lbrack 1\rbrack\boldsymbol{\cdot}\mathbb{Z}/2\lbrack 0\rbrack$$ with the non-trivial actions by $\mathbb{Z}/2$. Namely, recall from example \ref{ex:exotic2groupequivalence} that $\big(\mathbb{Z}/2\oplus\mathbb{Z}/2\big)\lbrack 1\rbrack\rtimes\mathbb{Z}/2\lbrack 0\rbrack\cong \mathbb{Z}/4\lbrack 1\rbrack\rtimes\mathbb{Z}/2\lbrack 0\rbrack$. We have that $\mathbf{2TY}(\mathbb{Z}/2,triv)\simeq \mathbf{2Rep}(\big(\mathbb{Z}/2\oplus\mathbb{Z}/2\big)\lbrack 1\rbrack\times\mathbb{Z}/2\lbrack 0\rbrack)$. Moreover, it follows from example \ref{ex:2category2rep2group} that the underlying 2-category of $\mathbf{2Rep}(\mathbb{Z}/4\lbrack 1\rbrack\boldsymbol{\cdot}\mathbb{Z}/2\lbrack 0\rbrack)$, which is depicted in example \ref{ex:beyondTamabraYamagami} above, is not a Tambara-Yamagami fusion 2-category. This proves the claim.

Let us write $\kappa_1$ and $\kappa_2$ for two 4-cochains representing the two cohomology classes of $H^4_{gp}(D_8;\mathds{k}^{\times})$ whose image in $H^4_{gp}(\mathbb{Z}/2\oplus \mathbb{Z}/2;\mathds{k}^{\times})$ are non-zero. It follows from a variant of the above discussion that both of these fusion 2-categories have a Tamabra-Yamagami defect and admit a fiber 2-functor. Further, at least one of them must be equivalent to $\mathbf{2Rep}(\mathbb{Z}/4\lbrack 1\rbrack\boldsymbol{\cdot}\mathbb{Z}/2\lbrack 0\rbrack)$. We suspect that they are actually both equivalent to this last fusion 2-category.
\end{Example}

\begin{Remark}
Parts of this last example can be readily generalized to the case of an arbitrary finite cyclic group $\mathbb{Z}/n$. More precisely, every Tambara-Yamagami 2-category $\mathbf{2TY}(\mathbb{Z}/n,\pi)$ admits a fiber 2-functor. Namely, the subgroup of $\mathbb{Z}/n\wr \mathbb{Z}/2$ generated by $(1,0,1)$ is a complement of $\mathbb{Z}/n\oplus 0$, and is isomorphic to $\mathbb{Z}/2n$. If $n$ is not a power of $2$, then such fusion 2-categories are not equivalent to the 2-category of 2-representations of a finite 2-group. In general, we do not know whether every Tambara-Yamagami 2-category admits a fiber 2-functor.
\end{Remark}

\bibliography{references.bib}
\end{document}